\documentclass[10pt]{article}
\usepackage{amsmath}
\usepackage{amssymb}
\usepackage{authblk}
\usepackage{subcaption}
\usepackage{tikz}
\usepackage{soul}
\usepackage{color}

\usepackage[square,numbers]{natbib}
\usepackage[left=2.5cm,top=2.5cm,right=2.5cm, bottom=2.5cm]{geometry}

\usepackage{graphicx}
\usepackage{enumerate}

\usepackage{draftwatermark}
\SetWatermarkText{Submitted}
\SetWatermarkScale{5.5} 
\SetWatermarkColor[gray]{0.9} 

\newcommand{\bm}[1]{\boldsymbol{#1}}
\newcommand{\dd}{\mathrm{d}}
\newcommand{\pd}{\partial}
\newcommand{\Tr}{\intercal}

\newcommand{\lam}[1]{{\langle{#1}\rangle}}
\newcommand{\elm}[1]{{({#1})}}

\newcommand{\fpd}[2]{\frac{\pd #1}{\pd #2}}

\newcommand{\figref}[1]{\ref{#1}}

\newcommand{\jump}[1]{[[ #1 ]]}
\newcommand{\ave} [1]{\{ #1 \}}

\title{The Immersed Boundary Conformal Method for Kirchhoff-Love and Reissner-Mindlin shells}

\author[1]{Giuliano Guarino}
\author[2]{Alberto Milazzo}
\author[1]{Annalisa Buffa}
\author[1]{Pablo Antolin}
\affil[1]{Institute of Mathematics, École Polytechnique Fédérale de Lausanne, CH-1015, Lausanne, Switzerland}
\affil[2]{Department of Engineering, Universit\`a degli Studi di Palermo, 90128, Italy}

\begin{document}

\maketitle

\begin{abstract}
This work utilizes the Immersed Boundary Conformal Method (IBCM) to analyze Kirchhoff-Love and Reissner-Mindlin shell structures within an immersed domain framework. Immersed boundary methods involve embedding complex geometries within a background grid, which allows for great flexibility in modeling intricate shapes and features despite the simplicity of the approach. The IBCM method introduces additional layers conformal to the boundaries, allowing for the strong imposition of Dirichlet boundary conditions and facilitating local refinement. In this study, the construction of boundary layers is combined with high-degree spline-based approximation spaces to further increase efficiency. The Nitsche method, employing non-symmetric average operators, is used to couple the boundary layers with the inner patch, while stabilizing the formulation with minimal penalty parameters. High-order quadrature rules are applied for integration over cut elements and patch interfaces.     
Numerical experiments demonstrate the efficiency and accuracy of the proposed formulation, highlighting its potential for complex shell structures modeled through Kirchhoff-Love and Reissner-Mindlin theories. These tests include the generation of conformal interfaces, the coupling of Kirchhoff-Love and Reissner-Mindlin theories, and the simulation of a damaged shell.
\end{abstract}

\section{Introduction}\label{sec:INTRO}

Two-dimensional manifolds embedded in three-dimensional space are a fundamental topological structure frequently encountered in physics. Examples of phenomena involving surfaces include the distribution of soap bubbles minimizing surface tension, the growth patterns of leaves, and the propagation of defects and dislocations in crystals. While these phenomena can often be modeled using three-dimensional equations derived from conservation principles, the significant difference in scale between thickness and surface directions suggests considerably distinct behaviors for the problem variables. By making appropriate a-priori assumptions about the through-thickness response, it is possible to achieve dimension reduction with minimal loss in model accuracy.

Although such two-dimensional approaches facilitate faster numerical approximations, they introduce more complex formulations compared to their three-dimensional counterparts. This complexity necessitates borrowing concepts from surface differential geometry \cite{ciarlet2005introduction}. Specifically, the physical domain must be expressed as a surface where geometrical quantities, such as local bases, tangent and normal vectors, as well as differential operators, can be effortlessly computed. This often imposes additional requirements on the regularity of the surface itself. Simple surfaces can be constructed using analytical mappings from the unit square, but without altering the topology. However, in practical applications, surfaces of interest can exhibit significantly more complex geometrical features, including holes and kinks. For such objects, finding an effective mathematical representation is a non-trivial task.

A possible solution is to define the surface implicitly using a level-set function. Typically, the level-set depends on the Cartesian coordinates, and the set of points where is equal to a specific value corresponds to either the surface itself or the boundary of an existing surface. The main advantage of this method is that leads to efficient algorithms when paired with simple tensor-product background discretizations \cite{saye2015high, saye2022high}. However, although some implicitly-defined surfaces can reach high level of complexity, the flexibility of the level-set is quite limited, and not every surface can be represented in such way. 

A different approach, more commonly used in industrial CAD design, involves adopting multiple surfaces explicitly represented in curvilinear coordinates through spline technology \cite{piegl1997,rogers2001}. Splines offer extensive flexibility in representing curvature profiles; however, standard splines are based on the concept of patch that results from a mapping of a rectangular parametric domain. For complex topologies, multiple patches are juxtaposed at their edges to meet the required geometric continuity, whether $G^0$ (position continuity), $G^1$ (tangent continuity), or $G^2$ (curvature radius continuity) \cite{kapl2019isogeometric, farahat2023isogeometric, weinmuller2022}. Enforcing these continuity requirements becomes increasingly difficult when moving from $G^0$ to $G^2$. These mergings often lack exact solutions, resulting in geometrical inconsistencies at the interfaces. A branch of research on splines in this direction has been focused on developing new unstructured spline technologies capable of locally enriching the design space, abandoning the tensor-product structure of typical spline meshes. Examples of this trend can be found in works such as \cite{casquero2017arbitrary, casquero2020, thomas2022, wen2023}. To some extent, in these approaches the concept of designing by connecting patches is abandoned, demanding for a drastic change in practitioners' habits.

Another method that allows for modeling topologies equivalent to polygons with a finite number of cutouts is the trimming approach \cite{mantyla1987introduction, mortenson1997geometric}. In this method, an underlying surface is created using spline functions to achieve the desired curvature profile that meets design requirements. The surface is then trimmed to account for local features on the boundary or to model internal holes, retaining the tensor product nature of the surface map while drastically increasing design flexibility. However, multiple patches are still necessary if the desired geometry includes kinks.

B-splines are a powerful tool not only for the design phase but also for the numerical analysis of partial differential equations. The order and regularity of spline-based functions can be arbitrarily increased straightforwardly, allowing higher-order problems that require higher-continuity approximation spaces to be seamlessly addressed numerically. Additionally, spline-based approximations results in high convergence rates that typically grow with the polynomial order, achieving efficient and accurate computation without a significant increase in the number of degrees of freedom \cite{bressan2019}. Two-dimensional spline spaces can be easily constructed over explicitly defined surfaces. When the functional spaces used in defining the underlying surface and the approximation function coincide, the resulting method is referred to as Isogeometric Analysis (IGA) \cite{hughes2005isogeometric}. For trimmed surfaces, a tensor-product spline-based approximation space can still be constructed and adopted for analysis \cite{breitenberger2015, bauer2017embedded, antolin2019isogeometric}, but this comes with some complications: i) Dirichlet boundary conditions cannot be imposed in a strong sense on trimmed boundaries. ii) High-accuracy integration over cut elements requires ad-hoc quadrature rules. iii) Badly-cut elements can negatively affect the stability and conditioning of the method.

An innovative approach to circumvent some of the issues associated with trimming was recently proposed in \cite{wei2021}. This approach, called the Immersed Boundary Conformal Method (IBCM), involves creating an additional layer for each immersed boundary. These layers are generated in such a way to make the domain's boundary conformal, so it can be represented exactly by the computational model's basis functions, allowing for a strong imposition of Dirichlet boundary conditions. The boundary layer is then coupled with the internal patch through an interface that is conformal to the former but non-conformal to the latter. Similar to Dirichlet boundary conditions, coupling conditions along non-conforming edges also need to be enforced in a weak sense. To this end, various techniques can be used, such as Lagrange multipliers methods \cite{babuska1973lagrange,brivadis2015isogometric, dittman2019weak, chasapi2020patch, benvenuti2023isogeometric}, pure-penalty methods \cite{babuska1973penalty, schmidt2012, lei2015,breitenberger2015}, and Nitsche-based methods \cite{douglas2008interior,riviere1999improved,bauman1999discontinuous, nguyen2013nitsche, ruess2014weak, hu2018skew}. 

In particular, the last class relies on interface integrals that use the jump of the main variables and the average of the formulation's fluxes between two edges of an interface. Although variationally consistent, this method requires stabilization terms based on penalty parameters, which however become unbounded when the interfaces lie on cut elements. The IBCM method can leverage the conformal nature of the boundary layer at the interface to mitigate this issue. Specifically, the average operator can be defined considering fluxes only from the non-cut boundary layer, maintaining stability and effectively bounding the penalty parameters.

Among the advantages of this method, it is worth mentioning its capability to facilitate local refinement, a task that can be challenging with B-spline based approximation spaces. This benefit, along with the strong application of Dirichlet conditions on otherwise trimmed boundaries, were demonstrated in \cite{wei2021} through multiple numerical tests encompassing Laplace, two-dimensional elasticity, and advection-diffusion equations. Additionally, the IBCM can be utilized to model more accurate physical phenomena that occur locally. For instance, in \cite{lapina2024}, conformal layers are introduced between different regions of the domain to accurately model the behaviour at their interface. As a final remark, when dealing with trimmed cut-outs with size of the same scale as the B-spline based mesh elements, the problem of cross talk can occur \cite{coradello2020}. Since spline bases span across multiple elements, for small cut-outs the approximation space would still maintain some level of continuity across the sides of the cut-out although these are separated. The degree of freedom of a basis function affected by the cutouts should then be replicated on elements on different sides of the edge, which is a somewhat laborious task. As an alternative, the boundary layer introduces necessary discontinuities in the approximation space along the cut-outs, effectively addressing this issue.

When partial differential equations in the field of solid mechanics pertain to thin-walled domains, the resulting problems fall within the realm of shell analysis. The most emblematic example within this class is shell elasticity theory, which addresses inquiries spanning sectors such as transportation, biology, construction, and energy. Structural shell elements, in particular, are utilized to effectively withstand external loads. Due to the curvature effect, forces applied on the shell surface and its boundaries distribute more efficiently, maximizing the material's utilization \cite{reddy2006}. In this context, an innovative symbiosis that is increasingly becoming standard for high-performance structural applications is the integration of shell-like designs with laminated composite materials. These materials typically consist of layers of more resistant fibers oriented arbitrarily within a more ductile matrix, providing designers with flexibility to optimize strength in more solicited directions \cite{jones1998}.

Classical shell theories include the Kirchhoff-Love \cite{kirchhoff1850, love1888, schollhammer2019kirchhoff,kiendl2009} and the Reissner-Mindlin ones \cite{reissner1945, mindlin1951, schollhammer2019reissner,benson2010isogeometric}. Both theories are based on the assumption that a straight unit segment, perpendicular to the mid-surface, remains straight and unstretched after deformation. In the Reissner-Mindlin theory, in addition to the displacement of the mid-surface, a second primary variable representing the rotation of the segment is introduced, leading to a second-order problem. Conversely, the Kirchhoff-Love theory restricts primary variables to displacement alone, assuming that the segment remains perpendicular to the mid-surface even after deformation. However, this simplification transforms the problem into a fourth-order one, demanding for $C^1$-continuity basis functions, and limits its applicability to thin-shell scenarios. To capture out-of-plane stress components more accurately, higher-order theories employ refined through-the-thickness approximations of the displacement field \cite{carrera2003theories, patton2021efficient}.

In the context of shell structures, weak coupling between different parts of the domain is often achieved through Lagrange multipliers methods \cite{horger2019hybrid, sommerwerk2017reissner, hirschler2019dual, schu2019multi}, or pure-penalty methods \cite{herrema2019penalty, leonetti2020robust, zhao2022opensource, pash2021priori, coradello2021coupling, coradello2021projected, prosepio2022penalty, guarino2023discontinuous}. However, its variationally-consistent counterpart has also been extensively investigated using Nitsche-based coupling for Kirchhoff-Love plates and shells \cite{noels2008, noels2009discontinuous, zhang2016, bonito2020discontinuous}, Reissner-Mindlin plates and shells \cite{arnold2005, bosing2010, talamini2017, mu2018}, and higher-order plates and shells \cite{gulizzi2020implicit, guarino2021equivalent, guarino2021high, guarino2022accurate}. Particularly interesting is the application of weak shell coupling to IGA, with various examples regarding Kirchhoff-Love shells for both conforming \cite{guo2015, nguyen2017isogeometric}, and trimmed patches \cite{guo2017parameter, guo2018variationally, yu2023isogeometric, guo2021isogeometric, chasapi2023fast, guarino2024interior}, and Reissner-Mindlin plates \cite{song2024reissner}. However, to the best of the author's knowledge, there are no examples of weak coupling of Reissner-Mindlin IGA shells as well as for higher-order theories. 

In this work, the IBCM method is utilized to analyze Kirchhoff-Love and Reissner-Mindlin shell structures within an immersed domain framework, where ad-hoc layers are created to ensure conformal boundaries. The approximation spaces defined over the discretized domains are based on B-spline functions in order to take advantage of the higher degree. To stabilize the formulation with minimal penalty parameters, the Nitsche method employing non-symmetric average operators is employed. To perform integration over cut elements, high-order quadrature rules are obtained following the algorithm developed in \cite{antolin2022robust}.

In Section \ref{sec:IBCM}, we introduce the discretization technique based on the IBCM method, along with the definition of spline spaces over the discretized domain. Section \ref{sec:SHELLS} delves into the Kirchhoff-Love and the Reissner-Mindlin shell theories following a presentation of essential concepts from differential geometry. The discretized versions of the variational statements are formulated in Section \ref{sec:IP}, where we additionally describe the Nitsche-based coupling technique. Section \ref{sec:RESULTS} presents the results from numerical experiments demonstrating the efficiency and accuracy of the proposed formulation. Finally, Section \ref{sec:CONCLUSIONS} draws the conclusions of the present study, highlighting the key advantages recognized in the investigated method.

\section{The IBCM on surface manifolds} \label{sec:IBCM}

\subsection{Surface map} \label{ssec:IBCM surface map}
Let us denote as $\Omega \in \mathbb{E}^3$ the manifold surface, where $\mathbb{E}^3$ represents the three-dimensional Euclidean space, and as $\partial \Omega$ its boundary. Let us suppose, following the embedded concept, that $\Omega$ is constructed as a restriction of an untrimmed surface $\Pi_0 \in \mathbb{E}^3$, with boundary $\partial \Pi_0$. This untrimmed surface is obtained as the image of an untrimmed reference domain $\hat{\Pi}_0\in\mathbb{R}^2$ through the function $\hat{\bm{\mathcal{F}}}: \hat{\Pi}_0 \rightarrow \Pi_0$. Thus, a generic point $\bm{x} \in \Pi_0$ is obtained as:
\begin{equation} \label{eq:IBCM map}
    \bm{x}(\xi_1,\xi_2)= x_i(\xi_1,\xi_2)\bm{e}_i = \hat{\bm{\mathcal{F}}}(\xi_1,\xi_2) \;, \quad\mathrm{with} \; i=1,2,3\;,
\end{equation}
where the components of $\bm{x}$ refer to the standard Euclidean basis $\bm{e}_1,\bm{e}_2,\bm{e}_3$, and $\xi_1,\xi_2$ are the curvilinear coordinates such that $(\xi_1,\xi_2)\in\hat{\Pi}_0$. Consistently, the boundary of the untrimmed surface is obtained as $\pd\Pi_0=\hat{\bm{\mathcal{F}}}(\pd\hat{\Pi}_0)$, where $\pd\hat{\Pi}_0$ denotes the boundary of $\hat{\Pi}_0$. As typical in structured maps, $\hat{\Pi}_0$ is selected as the reference rectangle
\begin{equation}
    \hat{\Pi}_0 = [\xi_{1b},\xi_{1t}]\times[\xi_{2b},\xi_{2t}] \;,
\end{equation}
where $\xi_{\alpha b}$ and $\xi_{\alpha t}$ denote the bottom and top limits, respectively, over which the curvilinear variable $\xi_\alpha$ ranges, and are arbitrary parameters, tuned on the desired final shape. For instance, when adopting Bézier surface maps, $\hat{\Pi}_0$ corresponds to the unit square $[0,1]\times[0,1]$. 

It is pointed out that in the previous equations and throughout the remainder of the paper, Greek letter indices are adopted for curvilinear coordinates and take values in $\{1,2\}$, while Latin letter indices are adopted for Cartesian coordinates and take values in $\{1,2,3\}$, unless stated otherwise. Regarding superscripts and subscripts, the Einstein summation convention is employed over repeated indices. Additionally, quantities referring to the parametric space are identified by an accent $\hat{\bullet}$,  whereas corresponding quantities in the Euclidean space do not have one.

\subsection{Embedded domain}
By defining the map in Eq.\eqref{eq:IBCM map}, the underlying untrimmed surface is effectively molded with the desired curvature profile. The actual surface geometry $\Omega$ is then identified by outlining its boundary $\pd\Omega$ within $\Pi_0$. Through the inverse map, one can obtain $\pd\hat{\Omega}=\bm{\mathcal{F}}^{-1}(\pd\Omega)$ and use it to delimit the parametric domain $\hat{\Omega}$. More specifically, $\pd\hat{\Omega}$ and $\pd\Omega$ divide $\hat{\Pi}_0$ and $\Pi_0$, respectively, into two parts. One part for each set are consistently selected as active portions, coinciding with $\hat{\Omega}$ and $\Omega$, while the remaining ones are referred to as the non-active portions. In turn, once $\hat{\Omega}$ and $\Omega$ are constructed, the relationship $\Omega=\bm{\mathcal{F}}(\hat{\Omega})$ holds. It is worth mentioning that in some circumstances, it can be more convenient to define the boundary directly in the parametric domain $\hat{\Pi}_0$ rather than in the physical one ${\Pi}_0$, and obtain the latter through the mapping process. 

So far, no assumptions have been made regarding the shape of the boundary of the parametric domain and the technique to construct it. However, in this paper, the strategy adopted is based on trimming operations. Accordingly, $\hat{\Omega}$ is derived from the reference rectangle by subsequently removing regions that do not belong to the desired shape. Let us assume the number of regions to remove is $N_\Gamma$, and denote the $i$-th region as $\hat{\Pi}_i$, delimited by the simply-connected closed trimming curve $\pd\hat{\Pi}_i$. The former may be either internal or external to the latter, and $i$ ranges from $1$ to $N_\Gamma$. Without loss of generality let us also assume that $\pd\hat{\Pi}_i\cap\pd\hat{\Pi}_j=\emptyset$ if $i\ne j$, and that the first trimming curve $\pd\hat{\Pi}_i$ delimits the external boundary eventually coinciding with $\pd\hat{\Pi}_0$, while the remaining trimming curve outline internal cut-outs. Under these assumptions, $\hat{\Omega}$ is given by
\begin{subequations}
\begin{equation}
    \hat{\Omega}= \hat{\Pi}_0 \setminus \bigcup_{i=1}^{N_\Gamma}\mathrm{cl}(\hat{\Pi}_i) \;,
\end{equation}
and its boundary by
\begin{equation} \label{eq:IBCM domain boundary}
    \pd\hat{\Omega}= \bigcup_{i=1}^{N_\Gamma}\pd\hat{\Pi}_i \;,
\end{equation}
\end{subequations}
where cl$(\bullet)$ stands for the closure of $\bullet$, therefore being $\mathrm{cl}(\hat{\Pi}_i) = \hat{\Pi}_i \cup \pd\hat{\Pi}_i$. As such, general surfaces topologically equivalent to a two-dimensional polygon with an arbitrary number of holes can be easily constructed. 

To better understand these definitions, Figs.\figref{fig:IBCM - Construction a} to \figref{fig:IBCM - Construction c} illustrate the construction of $\hat{\Omega}$. Fig.\figref{fig:IBCM - Construction a} shows a reference rectangular parametric domain. Through the trimming operation shown in Fig.\figref{fig:IBCM - Construction b} the regions $\hat{\Pi}_1$ and $\hat{\Pi}_2$ are subtracted from $\Pi_0$. In Fig.\figref{fig:IBCM - Construction c}, the resulting parametric domain $\hat{\Omega}$ and the correspondent boundary $\pd\hat{\Omega}$ are depicted.

This construction offers notable flexibility in defining the surface, combining the generality of curvature profiles from the map function with the arbitrariness of boundary selection, making it a powerful tool for designing complex shapes. Nonetheless, this flexibility comes at a cost when trimming curves are defined in Euclidean space and only subsequently mapped back into $\mathbb{R}^2$. In fact, the inverse of Eq.\eqref{eq:IBCM map} does not have a closed-form solution in general, so the trimming curves in the parametric domain are usually approximations, accurate only up to a specific geometric tolerance.

\subsection{Boundary conformal layers} \label{ssec:IBCM boundary layers}
As expressed in Eq.\eqref{eq:IBCM domain boundary}, the boundary of the parametric domain upon completion of the trimming operations is expressed as the union of $N_\Gamma$ simply-connected oriented closed boundary curve $\pd\hat{{\Pi}}_i$, where two different curves are disjoint sets. Let us assume now that the $i$-th curve is parameterized through the curvilinear coordinate $\eta^i\in[\eta^i_{b},\eta^i_{t}]$, where $\eta^i_{b}$ and $\eta^i_{t}$ define the bottom and top limits, respectively, of the interval spanned by $\eta^i$. Then, for each $\pd\hat{\Pi}_i$, an offset curve $\pd\hat{\Pi}'_i$ is constructed in a region internal to $\hat{\Omega}$, adopting the same parameterization of the associated boundary curve, i.e., $\eta^i\in[\eta^i_b,\eta^i_t]$. Additionally, let us assume that these offset curves are non-intersecting with each other and with boundary curves, meaning
\begin{align*}
    &\pd\hat{\Pi}'_i \cap \pd\hat{\Pi}'_j= \emptyset  \quad \mathrm{if} \; i\ne j \;,\\
    &\pd\hat{\Pi}_i \cap \pd\hat{\Pi}'_j= \emptyset   \quad \forall (i,j) \;.
\end{align*}
Apart from these requirements, there remains relative arbitrariness in their construction. A simple approach is to maintain a uniform distance from the target boundary curve. However, if the local radius is smaller than the selected distance, this can alter the topology of the offset curve, making it no longer simply-connected. Furthermore, while satisfying the uniform distance criterion in  $\mathbb{R}^2$ is relatively simple, achieving uniform distance on a surface manifold in $\mathbb{E}^3$ poses a significantly more challenging task. A simpler strategy involves creating a curve that is internal to the parametric domain, even if it does not precisely follow the geometric features of the boundary. Although this may appear to be a critical choice during the development of a discretization to approximate the problem of interest, \cite{wei2021} demonstrated that the method's accuracy remains robust regardless of whether the curve adheres strictly to the boundary's geometry or not.

It is important to mention that in those cases where the domain presents local features tight to each other, constructing non-intersecting offset curves can be challenging. Indeed, such geometrical configurations impose limits on the local distance between the offset curve and its corresponding boundary one. These critical cases are not addressed in the present contribution, while being left for future developments.

The approach proposed here relies on a domain sub-structuring strategy where each pair of boundary and offset curves are used to construct a boundary layer and to restrict the main patch's extension through an additional trimming operation. More specifically, the offset curve $\pd\hat{\Pi}'_i$ divides $\mathbb{R}^2$ into two regions. The one that contains the boundary curve $\pd\hat{\Pi}_i$, that by construction is entirely either internal or external to $\pd\hat{\Pi}'_i$, and its complement to $\mathbb{R}^2$. The former, in particular, is denoted as $\hat{\Pi}'_i$. The region enclosed between $\pd\hat{\Pi}'_i$ and $\pd\hat{\Pi}_i$ is the boundary layer $\hat{\Omega}_i$ defined as
\begin{align}
    &\hat{\Omega}_{i}=\hat{\Omega}\cap \hat{\Pi}'_i     \;,\\
    &\pd\hat{\Omega}_i= \pd\hat{\Pi}_i \cup \pd\hat{\Pi}'_i  \;,
\end{align}
being $\pd\hat{\Omega}_i$ the boundary of $\hat{\Omega}_i$. Once the boundary layers are constructed, the remaining part of the domain is defined through
\begin{align}
    &\hat{\Omega}_0 =\hat{\Omega}\setminus\bigcup_{i=1}^{N_\Gamma}\mathrm{cl}(\hat{\Omega}_i) \;,\\
    &\pd\hat{\Omega}_0 =\bigcup_{i=1}^{N_\Gamma} \pd\hat{\Pi}'_i \;.
\end{align}
The corresponding surfaces and boundary curves in the physical space are then obtained through mapping, i.e., $\Omega_i=\hat{\bm{\mathcal{F}}}(\hat{\Omega}_i)$, and $\pd\Omega_i=\hat{\bm{\mathcal{F}}}(\pd\hat{\Omega}_i)$ for $i=0,\cdots, N_\Gamma$. The interfaces between the main domain and the boundary layers are introduced in the parametric domain as $\hat{\Gamma}_i=\pd\hat{\Omega}_0\cap\pd\hat{\Omega}_i$, as well as their counterpart in the physical space $\Gamma_i=\hat{\bm{\mathcal{F}}}(\hat{\Gamma}_i)$.

\begin{figure}	\centering
    \includegraphics[width = 0.99\textwidth]{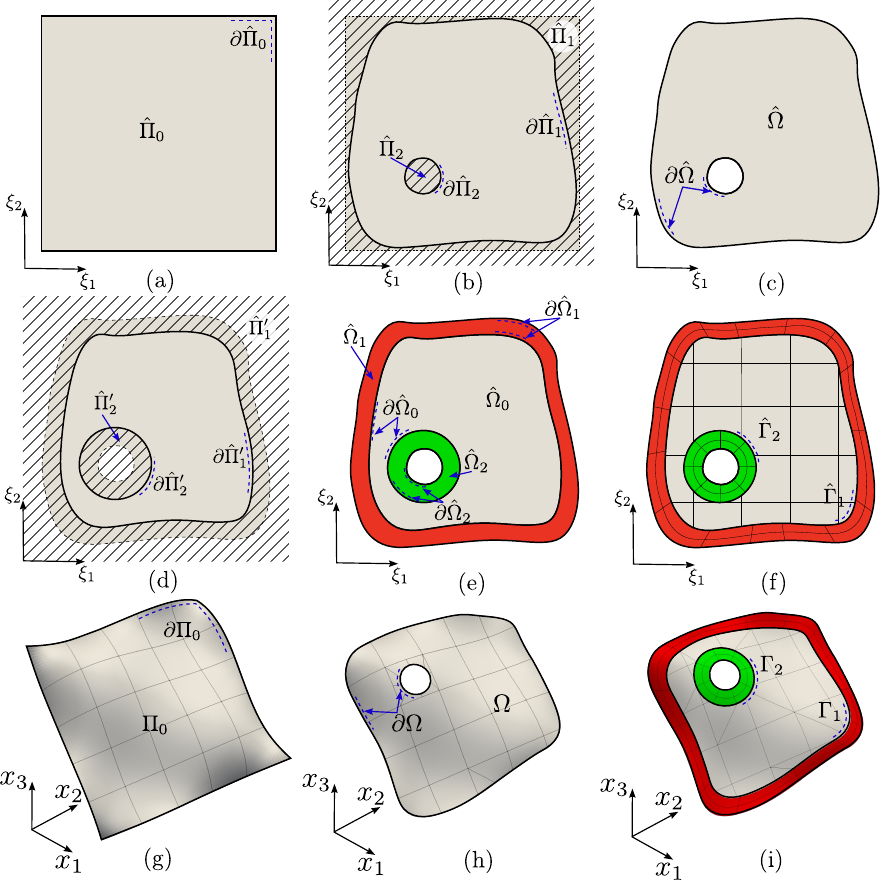}
    \begin{subfigure}{0\textwidth}\phantomcaption\label{fig:IBCM - Construction a}\end{subfigure}
    \begin{subfigure}{0\textwidth}\phantomcaption\label{fig:IBCM - Construction b}\end{subfigure}
    \begin{subfigure}{0\textwidth}\phantomcaption\label{fig:IBCM - Construction c}\end{subfigure}
    \begin{subfigure}{0\textwidth}\phantomcaption\label{fig:IBCM - Construction d}\end{subfigure}
    \begin{subfigure}{0\textwidth}\phantomcaption\label{fig:IBCM - Construction e}\end{subfigure}
    \begin{subfigure}{0\textwidth}\phantomcaption\label{fig:IBCM - Construction f}\end{subfigure}
    \begin{subfigure}{0\textwidth}\phantomcaption\label{fig:IBCM - Construction g}\end{subfigure}
    \begin{subfigure}{0\textwidth}\phantomcaption\label{fig:IBCM - Construction h}\end{subfigure}
    \begin{subfigure}{0\textwidth}\phantomcaption\label{fig:IBCM - Construction i}\end{subfigure}
    \caption{Step-by-step construction of the internal patch and boundary layers with corresponding discretizations. (a) Untrimmed parametric domain $\hat{\Pi}_0$ and associated boundary. (b) Construction of the parametric domain through trimming the regions $\hat{\Pi}_1$ and $\hat{\Pi}_2$ delimited by the corresponding trimming curves $\pd\hat{\Pi}_1$ and $\pd\hat{\Pi}_2$. (c) Resulting parametric $\Omega$ domain and associated boundary $\pd\Omega$. (d) Offset curves $\partial\hat{\Pi}'_1$ and $\partial\hat{\Pi}'_2$ and associated trimming regions $\hat{\Pi}'_1$ and $\hat{\Pi}'_2$. (e) Resulting internal domain $\hat{\Omega}_0$ and associated boundary layers $\hat{\Omega}_1$ and $\hat{\Omega}_2$ with correspondent boundaries $\pd\hat{\Omega}_0$, $\hat{\Omega}_1$, and $\pd\hat{\Omega}_2$, respectively. (f) Discretization of the internal patch and boundary layers, and interfaces $\hat{\Gamma}_1$ and $\hat{\Gamma}_2$. (g) Untrimmed surface $\Pi_0= \hat{\bm{\mathcal{F}}}(\hat{\Pi}_0)$ and correspondent boundary $\pd\Pi_0 = \hat{\bm{\mathcal{F}}}(\pd\hat{\Pi}_0)$ with superimposted discretization. (h) Trimmed surface $\Omega=\hat{\bm{\mathcal{F}}}(\hat{\Omega})$ and correspondent boundary $\pd\Omega=\hat{\bm{\mathcal{F}}}(\pd\hat{\Omega})$ with superimposted discretization. (i) Internal surface $\Omega_0 = \hat{\bm{\mathcal{F}}}(\hat{\Omega}_0)$ and correspondent boundary layers $\Omega_1 = \hat{\bm{\mathcal{F}}}(\hat{\Omega}_1)$ and $\Omega_2 = \hat{\bm{\mathcal{F}}}(\hat{\Omega}_2)$ with superimposed discretization and interfaces $\Gamma_1=\hat{\bm{\mathcal{F}}}(\hat{\Gamma}_1)$ and $\Gamma_2=\hat{\bm{\mathcal{F}}}(\hat{\Gamma}_2)$}.
\end{figure}

Figs.\figref{fig:IBCM - Construction d} and \figref{fig:IBCM - Construction e} show the construction of the internal domain and the corresponding boundary layers. In Fig.\figref{fig:IBCM - Construction d}, the parametric domain $\hat{\Omega}$ is restricted by trimming operations through the offset curves $\pd\hat{\Pi}'_1$ and $\pd\hat{\Pi}'_2$ that delimit the regions $\hat{\Pi}'_1$ and $\hat{\Pi}'_2$. The resulting internal domain $\hat{\Omega}_0$, and boundary layers $\hat{\Omega}_1$ and $\hat{\Omega}_2$ are shown in Fig.\figref{fig:IBCM - Construction e}, together with their boundaries $\pd\hat{\Omega}_0$, $\pd\hat{\Omega}_1$, and $\pd\hat{\Omega}_2$, respectively.

Provided the parameterization of $\pd\hat{\Pi}_i$ and $\pd\hat{\Pi}'_i$, that refers for construction to the same curvilinear coordinate $\eta^i\in[\eta^i_b,\eta^i_t]$, the boundary layer in the parametric domain can be described through an auxiliary map that relies on two curvilinear coordinates $\eta_1^i=\eta^i$ and $\eta_2^i$. The simplest choice consist in constructing this map as a ruled domain between the two curves, although non-convex curves might require more sophisticated algorithms. 

The boundary layer is therefore obtained in the parametric domain as the image of $\tilde{\Omega}_i=[\eta_{1b}^i,\eta_{1t}^i] \times [\eta_{2b}^i, \eta_{2t}^i]$ through the function $\tilde{\bm{\mathcal{F}}}_i: \tilde{\Omega}_i \rightarrow \hat{\Omega}_i$, where $\eta_{\alpha b}^i$ and $\eta_{\alpha t}^i$ are the limits for the interval of definition of $\eta^i_\alpha$. It is remarked that for construction $\eta^i_{1 b}$ and $\eta^i_{1 t}$ correspond to $\eta^i_b$ and $\eta^i_t$, respectively, while $\eta^i_2$ can be chosen, in absence of other constraints and without loss of generality, in the interval $[0,1]$. It is important pointing out that, in the notation introduced hereby, quantities referred to the auxiliary space are denoted by a tilde. Provided the parameterization of the $i$-th boundary layer, a generic point $\bm{\xi}=(\xi_1,\xi_2)$ in the parametric space is obtained as
\begin{equation}
    \bm{\xi}(\eta_1^i,\eta_2^i)= \xi_\alpha(\eta_1^i,\eta_2^i) = \tilde{\bm{\mathcal{F}}}_i(\eta_1^i,\eta_2^i) \;, \quad\mathrm{with} \; \alpha = 1,2 \; \mathrm{and} \; i=1,\cdots,N_\Gamma \;.
\end{equation}
with $(\eta_1^i,\eta_2^i)\in\tilde{\Omega}_i$. In turn, the domain of the boundary layer in the physical space is given by the composition of two maps $\Omega_i=\hat{\bm{\mathcal{F}}}\circ\tilde{\bm{\mathcal{F}}}_i(\tilde{\Omega}_i)$, and a generic point of the boundary layer in the physical space is given by 
\begin{equation} \label{eq:FOR - map composition}
    \bm{x}(\eta_1^i,\eta_2^i)= \hat{\bm{\mathcal{F}}} \circ \tilde{\bm{\mathcal{F}}}_i(\eta_1^i,\eta_2^i) \;, \quad\mathrm{with} \; i=1,\cdots,N_\Gamma \;.
\end{equation}
It is worth noting that the maps in the parametric domain $\tilde{\bm{\mathcal{F}}}_i(\pd\tilde{\Omega}_i)$ and in the physical domain $\hat{\bm{\mathcal{F}}}\circ \tilde{\bm{\mathcal{F}}}_i(\pd\tilde{\Omega}_i)$ do not correspond exactly to $\pd\hat{\Omega}_i$ and  $\pd\Omega_i$, respectively, since the curves identified by $\bm{\xi}(\eta^i_1=\eta^i_{1b},\eta^i_2)$ and $\bm{\xi}(\eta^i_1=\eta^i_{1t},\eta^i_2)$ are internal to $\hat{\Omega}_i$, and, similarly, the curves identified by $\bm{x}(\eta^i_1=\eta^i_{1b},\eta^i_2)$ and $\bm{x}(\eta^i_1=\eta^i_{1t},\eta^i_2)$ are internal to $\Omega_i$.

As detailed in the remainder of the paper, depending on the nature of the differential problem stated on the manifold, certain derivatives of the map in Eq.\eqref{eq:FOR - map composition} might be required. These derivatives can be straightforwardly obtained using the chain rule. For clarity, they are provided up to the third order in Appendix \ref{sec:APP - MAP COMP}.

\subsection{Discretization strategy} \label{ssec: discretization}
Each region resulting from the domain sub-structuring requires discretization to generate a suitable mesh for numerical investigations. The discretization algorithm employed here utilizes quadrilateral elements, with specific details varying depending on whether the region is a boundary layer or the internal domain.

\subsubsection{Discretization of the boundary layers}
For boundary layers, a structured, tensor-product-based, and quadrilateral grid is generated in each auxiliary domain $\tilde{\Omega}_i$. Within the interval of definition of $\eta_\alpha^i$,  $N_{e\alpha}^i+1$ increasing and non-repeating values $\theta^{i,j}_\alpha$ are identified collected in $\Theta^i_\alpha = \left\{\theta^{i,1}_\alpha,\cdots,\theta^{i,N_{e\alpha}^i+1}_\alpha\right\}$, where $\eta^i_{\alpha b}=\theta^{i,1}_\alpha$ and $\eta^i_{\alpha t}=\theta^{i,N_{e\alpha}^i+1}_\alpha$. As such, in the $i$-th auxiliary space the $e$-th element is defined as 
\begin{equation}
    \tilde{\Omega}_i^\elm{e} = \left[ \eta^{i,\elm{e}}_{1b}, \eta^{i,\elm{e}}_{1t} \right] \times \left[\eta^{i,\elm{e}}_{2b}, \eta^{i,\elm{e}}_{2t}\right] \;,
\end{equation}
where $\eta^{i,\elm{e}}_{\alpha b}$ and $\eta^{i,\elm{e}}_{\alpha t}$ are two consecutive values in $\Theta^i_\alpha$, and being $N_e^i=N_{e1}^i \times N_{e2}^i$ the number of elements of the $i$-th boundary layer. The corresponding element in the parametric and in the Euclidean spaces are given, respectively by 
\begin{subequations}
\begin{align}
    &\hat{\Omega}_i^\elm{e} = \tilde{\bm{\mathcal{F}}}_i\left(\tilde{\Omega}_i^\elm{e}\right) \;, \\
    &\Omega_i^\elm{e} = \hat{\bm{\mathcal{F}}}\circ\tilde{\bm{\mathcal{F}}}_i\left(\tilde{\Omega}_i^\elm{e}\right) \;.
\end{align}
\end{subequations}
It is remarked that in the previous equations $i$ spans $\{1,\cdots,N_\Gamma\}$, that is the set of indices associated with the boundary layers. A different discretization strategy based on trimmed elements is adopted for $\Omega_0$ as detailed in the subsequent section.

\subsubsection{Discretization of the internal domain}
Also for the internal domain, the discretization strategy starts from a background rectangular grid. In this case, however, the rectangular grid is constructed over the reference rectangular domain $\hat{\Pi}_0$. Within the interval of definition of $\xi_\alpha$, $N_{e\alpha}^0+1$ increasing and non-repeating values are identified and collected in $\Theta^0_\alpha = \left\{\theta^{0,1}_\alpha,\cdots,\theta^{0,N_{e\alpha}^0+1}_\alpha\right\}$, where $\xi_{\alpha b}=\theta^{0,1}_\alpha$ and $\xi^0_{\alpha t}=\theta^{0,N_{e\alpha}^0+1}_\alpha$. The $e$-th cell of the background grid is denoted as $\hat{\Pi}^\elm{e}_0$ and defined as
\begin{equation}
    \hat{\Pi}^\elm{e}_0 = \left[ \xi^{\elm{e}}_{1b}, \xi^{\elm{e}}_{1t} \right] \times \left[\xi^{\elm{e}}_{2b}, \xi^{\elm{e}}_{2t}\right] \;,
\end{equation}
where $\xi^{\elm{e}}_{\alpha b}$ and $\xi^{\elm{e}}_{\alpha t}$ are two consecutive values in $\Theta^0_\alpha$, and the total number of cells is $N_{e}^0=N_{e1}^0\times N_{e2}^0$. The background grid is then intersected with the actual internal domain $\hat{\Omega}_0$, and for the $e$-th cell the corresponding element is defined as 
\begin{align} 
	&\hat{\Omega}^\elm{e}_0 = \hat{\Pi}^\elm{e}_0\cap \hat{\Omega}_0 \;, \\
	&\partial\hat{\Omega}^\elm{e}_0 = (\partial\hat{\Pi}^\elm{e}_0\cap \hat{\Omega}_0 ) \cup (\hat{\Pi}^\elm{e}_0\cap \partial\hat{\Omega}_0 ) \;.
\end{align}
where $\pd\hat{\Pi}^\elm{e}_0$ is the boundary of the $e$-th background cell, and $\hat{\Omega}^\elm{e}_0 $ and $\pd\hat{\Omega}^\elm{e}_0$ are the domain and the boundary of the $e$-th element. It follows that, accordingly to the specificity of the intersection $\hat{\Pi}^\elm{e}_0\cap \hat{\Omega}_0$, the elements are classified into three groups: 
\begin{itemize}
    \item[-] The group of entire elements, for which $\hat{\Omega}^\elm{e}_0=\hat{\Pi}^\elm{e}_0$ ;
    \item[-] The group of partial elements, for which $\hat{\Omega}^\elm{e}_0\subset\hat{\Pi}^\elm{e}_0$ ;
    \item[-] The group of empty elements, for which $\hat{\Omega}^\elm{e}_0\cap\hat{\Pi}^\elm{e}_0=\emptyset$ .
\end{itemize}
The active elements are only the entire and the partial ones. Similarly to what stated for boundary layers, the elements' domain in the Euclidean space is obtained mapping their correspondents in the parametric space
\begin{equation}
    \Omega_i^\elm{e} = \hat{\bm{\mathcal{F}}}\left(\hat{\Omega}_0^\elm{e}\right) \;.
\end{equation}

\subsubsection{Quadrature rule}
The formulation considered in this work requires numerical integration over the domain and the boundary of elements, along with interfaces between contiguous patches. For elements in the boundary layers, as well as in the entire elements of the internal patch, integrals are computed using high-order tensor-product Gauss-Legendre quadrature rules in the corresponding non-mapped domain. On the other hand, integration on trimmed elements necessitates more refined strategies to accurately capture only their active parts.

In this work, the algorithm presented in \cite{wei2021,antolin2022robust} is adopted. In essence, this algorithm reparameterizes each partial element as a set of tiles, allowing standard Gauss-Lagrange quadrature to be applied in the auxiliary space of each tile. These tiles can be observed, e.g., in Figs.\figref{fig:IBCM - Construction h} and \figref{fig:IBCM - Construction i} for $\Omega$ and $\Omega_0$, respectively. The reparameterization starts by determining the topology of the element through the information regarding its boundary. Accordingly, to approximate the element a decomposition strategy is selected and high-order tiles are created as Bezier two-dimensional surfaces. It is essential to emphasize that the tiling process is solely used to obtain high-fidelity integration. From an analysis perspective, no new elements are introduced.

Integration over interfaces requires the subdivision of the common curve into segments, with each segment corresponding to a pair of elements belonging to distinct yet adjacent patches connected by the interface.  In cases where one or both sides of the interfaces are trimmed, this subdivision must consider individual tiles rather than the active elements. For further details on this segmentation, interested readers are referred to \cite{antolin2021overlapping}. 

Ultimately, as the integration rules are constructed within the non-mapped domain of the elements, appropriate terms are multiplied to the integration weights to perform transformation of infinitesimal area and length in surface and line integration, respectively.

\subsubsection{Strategy to address the ill-conditioning}
For the discretization employed herein, partial elements might be characterized by an unboundedly small ratio between the trimmed and the original area. When this situation occurs the condition number of the stiffness matrix critically increases, affecting the numerical stability of the solution.

In the construction of the boundary layers, while no control is left over trimming curves, whose definition depends on the desired shape, the construction of offset curves is to some extent arbitrary. This flexibility could be leveraged by an algorithm that ensures that the cutting process of the internal domain does not produce critical elements, given a certain homogeneity requirement the elements' size. While for a given discretization level the criteria could be met, on subsequent refinements motivated by the need for improved accuracy of the analysis, the homogeneity requirement might not be satisfied any longer.

As such, the approach preferred here is to avoid complicating the construction of offset curves and using instead a diagonal scaling of the stiffness matrix through Jacobi preconditioning. Although relatively simple, this method has proven to be very efficient in addressing the ill-conditioning issue caused by size difference on the support of shape functions. Nonetheless, in some critical cutting scenarios, shape functions having support only on the small cut element might become linearly dependent, making the matrix numerically singular. When this situation occurs, a dedicated preconditioner, e.g., the Schwarz preconditioner  \cite{deprenter2019,deprenter2023stability}, should be adopted. However, the simplicity of the Jacobi is preferred in this work, leaving more robust techniques to future studies. 

\subsection{B-spline based approximation spaces}\label{ssec:splines}
The space chosen for approximating the solution in the variational problems investigated hereafter leverages B-spline technology \cite{hughes2005isogeometric,cottrell2009}. This allows for constructing higher-continuity basis functions over the discretized domain. To this end, let us first define the B-spline space over non-mapped one-dimensional and two-dimensional domains.

Let us suppose that a uni-dimensional interval  $\mathrm{I}=[\xi_b,\xi_t]$, referred to the  coordinate $\xi$, is discretized using the breaks in the vector $\Theta=\left\{\theta^1,\dots,\theta^{N_e+1}\right\}$, with $\xi_b=\theta^1$ and $\xi_t=\theta^{N_e+1}$, being $N_e$ the number of elements of the discretization. To construct a spline-based approximation space over this discretized interval it is first defined the knot vector by taking adequate repetitions of the elements in $\Theta$, thus resulting in the vector of non-decreasing values $\Xi=\left\{\xi^1,\dots,\xi^{n+p+1}\right\}$, where $n$ is the number of shape functions and $p$ is the degree of the spline. The details of the construction of the knot vector are not reported here for the sake of conciseness. The interested reader is referred to \cite{piegl1997} for further explainations. It is worth mentioning, however, that the spline space is $C^{\infty}$ everywhere except at the internal knots, where the continuity is at most $p-1$ and it is reduced by one for every repetition of the knot. A generic univariate B-spline over $\mathrm{I}$ is constructed recursively using the Cox de Boor formula as:
\begin{subequations}
    \begin{align}
    N_{i,0}(\xi) &= \left\{\begin{array}{ll}
    1 &\mathrm{for}\ \xi^i\le \xi \le \xi^{i+1}\\
    0 &\mathrm{otherwise}\end{array}\right., \\
    N_{i,p}(\xi) &= \frac{(\xi-\xi^i)N_{i,p-1}(\xi)}{\xi^{i+p-1}-\xi^i}+\frac{(\xi^{i+p}-\xi)N_{i+1,p-1}(\xi)}{\xi^{i+p}-\xi^{i+1}} \;,
    \end{align}
\end{subequations}
where $N_{i,p}(\xi)$ is the $i$-th spline basis with $i\in\{1,\dots,n\}$ and degree $p$. In order to extend this construction to two dimension, let us now suppose that a rectangular domain $\Pi=[\xi_{1b},\xi_{1t}]\times[\xi_{2b},\xi_{2t}]$ is discretized using a rectangular grid based on a tensor product structure into $N_{e} = N_{e1}\times N_{e2}$ elements, being $N_{e\alpha}$ the number of elements in the direction $\alpha$. As such, $\Theta_\alpha = \left\{\theta^1_\alpha,\dots,\theta^{N_{e\alpha}+1}_\alpha\right\}$ is the list of breaks in the $\alpha$-th direction and $\Xi_1=\left\{\xi^1_\alpha,\dots,\xi^{n_\alpha+p+1}_\alpha\right\}$ is the correspondent knot vector. It is noted that, following this construction, the rectangular domain is also obtained as  $\Pi=[\theta_1^1,\theta^{N_{e1}+1}_1] \times[\theta_1^1,\theta^{N_{e2}+1}_1]$. With these definitions at hand a generic bi-variate B-spline is defined as 
\begin{equation} \label{eq:IBCM generic spline}
    B_{ij}(\xi_1,\xi_2)=N_i(\xi_1)N_j(\xi_2) \quad \mathrm{with}\;i\in\{1,\cdots,n_1\},j\in\{1,\cdots,n_2\}\;,
\end{equation}
where the secondary indices denoting the degree of $N_i(\xi_1)$ and $N_j(\xi_2)$ were dropped assuming they are both equal to $p$.

The approximation space for the $k$-th boundary layer is then obtained by first constructing bi-variate B-splines over $\tilde{\Omega}_k$ following the construction that leads to Eq.\eqref{eq:IBCM generic spline} and taking $\Pi=\tilde{\Omega}_k$ and $\Theta_\alpha=\Theta_\alpha^k$. The approximation space is therefore defined as
\begin{equation}
\mathcal{S}_{h,\Omega_k}=\mathrm{span}\{B_{ij}^k\circ\bm{\mathcal{F}}^{-1}_k: i\in\{1,...,n_1^k\} ,\; j\in\{1,...,n_2^k\}\}   \end{equation}
where $B_{ij}^k$ is the $ij$-th spline defined over $\tilde{\Omega}_k$, $\bm{\mathcal{F}}_k=\hat{\bm{\mathcal{F}}}\circ\tilde{\bm{\mathcal{F}}}_k$, and $n_\alpha^k$ is the number of spline functions in the $\alpha$-th direction for the boundary layer $k$.

A similar approach is adopted in constructing the approximation space for the internal patch $\Omega_0$. In this case the bi-variate B-spline are constructed over $\hat{\Pi}_0$, therefore selecting $\Pi=\hat{\Pi}_0$ and $\Theta_\alpha=\Theta_\alpha^0$. However, due to the trimmed definition of $\hat{\Omega}_0$ some of the basis functions $B_{ij}^0$ may result in a support consisting only on non-active elements. To remove this superfluous degrees of freedom, the approximation space is defined as
\begin{equation}
\mathcal{S}_{h,\Omega_0}=\mathrm{span}\{B_{ij}^0\circ\bm{\mathcal{F}}^{-1}_0: i\in\{1,...,n_1^0\} ,\; j\in\{1,...,n_2^0\} ,\; \mathrm{supp}(B_{ij}^0)\cap\Omega_0\neq \emptyset\} \;
\end{equation}
where $\mathrm{supp}(B_{ij}^0)$ is the union of the elements where $B_{ij}^0\ne 0$,  $\hat{\bm{\mathcal{F}}}_0=\hat{\bm{\mathcal{F}}}$, and $n_\alpha^0$ is the number of spline functions in the $\alpha$-th direction.

The approximation space over the entire domain is simply obtained as the vector spaces sum of the approximation spaces defined over each patch, in other words the internal one and the boundary layers, meaning
\begin{equation}
\mathcal{S}_{h}=\mathcal{S}_{h,\Omega_0}\oplus\mathcal{S}_{h,\Omega_1}\oplus\cdots\oplus\mathcal{S}_{h,\Omega_{N_\Gamma}} \;.
\end{equation}
where $\oplus$ denotes the vector spaces sum operation.

\section{Classical shell theories} \label{sec:SHELLS}

The problems investigated in this work pertain to the linear elastic analysis of laminated shell structures. Specifically, the starting point are the Kirchhoff-Love and Reissner-Mindlin shell equations formulated in a weak sense over the shell mid-surface. The following section presents both shell kinematics within a unified framework.

\subsection{Differential geometry}
In this section, we briefly introduce the essential concepts of differential geometry pertinent to the definition of the proposed problem. Given the map $\bm{\mathcal{F}}$ of the surface, the vectors of the local covariant basis are defined as
\begin{equation}
    \bm{a}_\alpha = \bm{\mathcal{F}}_{,\alpha} \;.
\end{equation}
A unit normal vector is associated to each point on the surface, and is computed from the covariant basis as
\begin{equation}
    \bm{a}_3 = \frac{\bm{a}_1 \times \bm{a}_2}{|\bm{a}_1 \times \bm{a}_2|} \;,
\end{equation}
where $\times$ denotes the cross product and $|\bullet|$ the standard Euclidean norm. The covariant components of the surface metric tensor are defined as $a_{\alpha\beta}=\bm{a}_\alpha\cdot\bm{a}_\beta$, where $\cdot$ denotes the dot product. The determinant of the surface metric tensor is denoted as $a$. The contravariant components of the surface metric tensor are obtained from the covariant counterparts as 
\begin{equation}
    [a^{\alpha\beta}] = [a_{\alpha\beta}]^{-1}\;
\end{equation}
and allow us to compute the vectors of the local contravariant basis as
\begin{equation}
    \bm{a}^\alpha= a^{\alpha\beta}\bm{a}_\beta \;,
\end{equation}
that satisfy the property $\bm{a}_\alpha\cdot\bm{a}^\beta = \delta_\alpha^\beta$, being $\delta_\alpha^\beta$ the Kronecker delta. Additionally, we introduce the covariant components of the curvature tensor, defined respectively as
\begin{equation}
        b_{\alpha\beta} = \bm{a}_3\cdot\bm{a}_{\alpha,\beta} \;,
\end{equation}
while its mixed components are obtained through the index raising operation
\begin{equation}
        b^\alpha_\beta = a^{\alpha\gamma}b_{\gamma\beta} \;.
\end{equation}

Leveraging the surface map and the definition of the normal unit vector, it is possible to map the volume constructed by extruding in the direction of $\bm{a}_3$. To achieve this, a third curvilinear coordinate $\xi_3$ is introduced. This coordinate represents the signed distance between a point in the volume and its normal projection onto the mid-surface, spanning the interval $\hat{\mathrm{I}}_3=[\xi_{3b},\xi_{3t}]$. Consequently, the map for the entire volume $\bm{\mathcal{G}}:\hat{V}\rightarrow V$ is defined as
\begin{equation}
    \bm{X}(\xi_1,\xi_2,\xi_3) = \bm{\mathcal{G}}(\xi_1,\xi_2,\xi_3) =  \bm{x}(\xi_1,\xi_2)+\xi_3 \bm{a}_3(\xi_1,\xi_2) \;,
\end{equation}
where $\hat{V}=\hat{\Omega}\times\hat{\mathrm{I}}_3$. Similarly to the surface mapping, for the volume mapping the vectors of the local covariant basis are defined as $\bm{A}_i=\bm{\mathcal{G}}_{,i}$. More specifically
\begin{subequations}\label{eq:DIFF volume covariant basis}
\begin{align}
    & \bm{A}_\alpha = \bm{a}_\alpha + \xi_3 \bm{a}_{3,\alpha}  \;, \\
    & \bm{A}_3 = \bm{a}_{3} \;.
\end{align}
\end{subequations}
It is worth noting that, following the convention utilized in this work, for mappings associated with boundary layers, the curvilinear coordinates $\eta_1^i\eta_2^i$ are employed for defining geometrical quantities pertaining the surface and the volume in lieu of $\xi_1\xi_2$. This clarification also applies in the treatment of shell theories.

\subsection{Problem setting}
Let us assume that the mid-surface of the shell is a two-dimensional manifold in $\mathbb{E}^3$, as described in Section \ref{ssec:IBCM surface map}. The shell material is a laminate obtained juxtaposing $N_\ell$ layers, assumed homogeneous, orthotropic, and perfectly bonded. The superscript $\lam{\ell}$ denotes quantities specific to the $\ell$-th layer. Thus, $V^\lam{\ell}$ represents the volume of the $\ell$-th layer, which is assumed to have a uniform thickness $\tau^\lam{\ell}$, resulting in the total thickness of the laminate $\tau$ being also uniform and given by
\begin{equation}
    \tau = \sum\limits_{\ell=1}^{N_\ell}\tau^{\lam{\ell}} \;,
\end{equation}
and the overall volume of the shell given by $V=\cup_{\ell=1}^{N_\ell} V^\lam{\ell}$. Each layer conforms to the curvature of the shell, meaning that given $\xi_{3b}^\lam{\ell}$ and $\xi_{3t}^\lam{\ell} = \xi_{3b}^\lam{\ell}+\tau^\lam{\ell}$ as the coordinates of the layer's bottom and top reference surfaces, respectively, the volume $V^\lam{\ell}$ comprises the points such that $\bm{X} = \bm{X}(\xi_1,\xi_2, \xi_{3b}^\lam{\ell}\leq\xi_3\leq\xi_{3t}^\lam{\ell})$. Concerning the material parameters such as Young moduli, Poisson ratios, shear moduli, and lamination angle, these are discussed in detail in Appendix \ref{sec:APP constitutive}, where the constitutive relationships for the classical shell theories are derived. Here, it is assumed that the shell mid-surface refers to $\xi_3=0$,  so the coordinates of the bottom and top surfaces of the laminate are $\xi_{3b}=-\tau/2$ and $\xi_{3t}=+\tau/2$.

To what regards the external loads, distributed forces and moments denoted by $\bar{\bm{f}}$ and $\bar{\bm{m}}$, respectively, are assumed applied directly on the mid-surface $\Omega$. Its boundary, $\pd\Omega$ is bipartite in two distinct ways: as $\pd\Omega = \pd\Omega_{D_1}\cup\pd\Omega_{N_1}$ for one set of conditions, and as $\pd\Omega = \pd\Omega_{D_2}\cup\pd\Omega_{N_2}$ for the other one:
\begin{itemize} 
    \item[-] $\pd\Omega_{D_1}$ and $\pd\Omega_{D_2}$ represent the portions of the boundary where Dirichlet displacement and rotation boundary conditions are applied, respectively.
    \item[-] $\pd\Omega_{N_1}$ and $\pd\Omega_{N_2}$ represent the portions of the boundary where Neumann force and moment boundary conditions are applied, respectively.
\end{itemize}
It is important to mention that $\pd\Omega_{D_1}$ and $\pd\Omega_{N_1}$ are disjoint sets whose union coincide with $\pd\Omega$, and the same applies to $\pd\Omega_{D_2}$ and $\pd\Omega_{N_2}$. It is also defined the set of corners $\chi\in\pd\Omega$ which is further divided into 
\begin{itemize} 
    \item[-] $\chi_D\in \mathrm{cl}({\pd\Omega_{D_1}})$ where Dirichlet displacement boundary conditions are applied.
    \item[-] $\chi_N\in\pd\Omega_{N_1}$ where Neumann force boundary conditions are applied. 
\end{itemize}
For the notation of boundary conditions, applied boundary force and moment are denoted as $\bar{\bm{F}}$ and $\bar{\bm{M}}$, respectively, while applied displacement and rotation are denoted as $\bar{\bm{u}}$ and $\bar{\bm{\theta}}$, respectively.

It is important to mention that while $\bar{\bm{f}}$, $\bar{\bm{F}}$, and $\bar{\bm{u}}$ can have three independent components in the Euclidean space, $\bar{\bm{m}}$, $\bar{\bm{M}}$, and $\bar{\bm{\theta}}$ lie in the plane locally tangent to the mid-surface. Additionally, although theoretically the division between different types of Dirichlet and Neumann boundaries could be component-based, meaning that different components of the displacement or rotation vectors could belong to different boundary types on the same edge, for simplicity, the classification of the boundary conditions is kept as described.

\subsection{Kinematic hypothesis, generalized strain, and generalized stress}
Shell theories are developed based on an axiomatic assumption about how the displacement varies along the coordinate $\xi_3$, which spans the thickness direction. The classical shell theories, namely  Kirchhoff-Love and Reissner-Mindlin, are founded on the kinematic hypothesis that a unit vector normal to the undeformed mid-surface remains straight and unextended after deformation. This assumption can be expressed as follows
\begin{equation}
    \bm{U}(\xi_1,\xi_2,\xi_3) = \bm{u}(\xi_1,\xi_2) +\xi_3\bm{\theta}(\xi_1,\xi_2) \;,
\end{equation}
where $\bm{U}$ denotes displacement of a generic point within the shell volume, $\bm{u}$ denotes displacement of its projection onto the mid-surface, and $\bm{\theta}$ the rotation vector. The assumption of the non-extensibility of the perpendicular unit vector implies that the rotation must lie within the tangent plane of the mid-surface, which leads to the condition that $\bm{a}_3\cdot\bm{\theta}=0$. The covariant components of the linear strain tensor for a three-dimensional solid \(V\) are derived from
\begin{equation} \label{eq:RM volume strain}
    \epsilon_{ij} = \frac{1}{2} (\bm{A}_i\cdot\bm{U}_{,j} + \bm{A}_j\cdot\bm{U}_{,i}) \;,
\end{equation}
where the derivatives of the displacement field with respect to the curvilinear coordinates, upon considering the kinematic hypothesis, simplify to
\begin{subequations} \label{eq:RM volume displacement derivatives}
\begin{align}
    & \bm{U}_{,\alpha} = \bm{u}_{,\alpha} + \xi_3 \bm{\theta}_{,\alpha} \;,\\
    & \bm{U}_{,3} = \bm{\theta} \;.
\end{align}
\end{subequations}
Introducing Eqs.\eqref{eq:DIFF volume covariant basis} and \eqref{eq:RM volume displacement derivatives} in Eq.\eqref{eq:RM volume strain} leads to the the following expressions for the components of the strain field
\begin{subequations}
\begin{align}
    \epsilon_{\alpha\beta} &= \frac{1}{2} \left[ (\bm{a}_\alpha + \xi_3 \bm{a}_{3,\alpha})\cdot(\bm{u}_{,\beta} + \xi_3 \bm{\theta}_{,\beta}) + (\bm{a}_\beta + \xi_3 \bm{a}_{3,\beta})\cdot(\bm{u}_{,\alpha} + \xi_3 \bm{\theta}_{,\alpha}) \right] =  \varepsilon_{\alpha\beta} + \xi_3\kappa_{\alpha\beta}+\xi_3^2\chi_{\alpha\beta} \;,\\
    \epsilon_{3\alpha}&= \frac{1}{2}\left[\bm{a}_{3}\cdot(\bm{u}_{,\alpha} + \xi_3 \bm{\theta}_{,\alpha}) + (\bm{a}_\alpha + \xi_3 \bm{a}_{3,\alpha})\cdot\bm{\theta}\right] = \frac{\gamma_\alpha}{2}\;,\\
    \epsilon_{33} &= 0  \;.
\end{align}
\end{subequations}
Note that in the previous equation $\epsilon_{3\alpha}=\epsilon_{\alpha 3}$. As typical in Reissner-Mindlin formulations, the tensor $\bm{\chi}$ is discarded since $\xi_3^2$ can be neglected for sufficiently small thickness values. The tensors $\bm{\varepsilon}$,  $\bm{\kappa}$, and  $\bm{\gamma}$ are referred to as the generalized strain tensors, specifically representing the membrane strain, bending strain, and shear strain, respectively. Their covariant components are defined as
\begin{subequations} \label{eq:RM - strain provv}
\begin{align}
    2\varepsilon_{\alpha\beta} &= \bm{a}_\alpha\cdot\bm{u}_{,\beta}+\bm{a}_\beta\cdot\bm{u}_{,\alpha} \;,\\
    2\kappa_{\alpha\beta} &= \bm{a}_\alpha\cdot\bm{\theta}_{,\beta}+\bm{a}_\beta\cdot\bm{\theta}_{,\alpha} + \bm{a}_{3,\alpha}\cdot\bm{u}_{,\beta} + \bm{a}_{3,\beta}\cdot\bm{u}_{,\alpha} \;,\\
    \gamma_{\alpha} &= \bm{a}_3\cdot\bm{u}_{,\alpha} + \bm{a}_{\alpha}\cdot\bm{\theta} \;,
\end{align}
\end{subequations}
in such way that $\bm{\varepsilon}=\varepsilon_{\alpha\beta}\,\bm{a}^\alpha\otimes\bm{a}^\beta$, $\bm{\kappa}=\kappa_{\alpha\beta}\,\bm{a}^\alpha\otimes\bm{a}^\beta$, $\bm{\gamma}=\gamma_{\alpha}\,\bm{a}^3\otimes\bm{a}^\alpha$. It is worth noting that in the previous equation it was used the property $\bm{a}_{3,\alpha}\cdot\bm{\theta} = - \bm{a}_3\cdot\bm{\theta}_{,\alpha}$, that follows from considering $\bm{a}_3\cdot\bm{\theta}=0$ in the equality $(\bm{a}_3\cdot\bm{\theta})_{,\alpha}= \bm{a}_{3,\alpha}\cdot\bm{\theta} + \bm{a}_3\cdot\bm{\theta}_{,\alpha}$. 
The work-conjugate quantities are known as generalized stresses and are denoted by as $\bm{N}$, $\bm{M}$, and $\bm{Q}$, representing the membrane force, bending moment, and shear force, respectively. The contravariant components of these tensors are derived using the constitutive relationships, as detailed in Appendix \ref{sec:APP constitutive}. These relationships are summarized as follows:

\begin{subequations}
    \begin{align}
        N^{\alpha \beta} &=  \mathbb{A}^{\alpha \beta \gamma \delta}\varepsilon_{\gamma \delta} + \mathbb{B}^{\alpha \beta \gamma \delta}\kappa_{\gamma \delta} \;,\\
        M^{\alpha \beta} &=  \mathbb{B}^{\alpha \beta \gamma\delta}\varepsilon_{\gamma \delta} + \mathbb{D}^{\alpha \beta \gamma \delta}\kappa_{\gamma \delta} \;,\\
        Q^\alpha &= \mathbb{S}^{\alpha\beta} \gamma_\beta \;.
    \end{align}
\end{subequations}

\subsection{The Reissner-Mindlin shell equation}
In the formulation adopted here, the Reissner-Mindlin equations use the Cartesian coordinates of the displacement vector $u_i$ (i.e., $\bm{u} = u_i \bm{e}_i$), and the covariant coordinates of the rotation vector $\theta_\alpha$ (i.e., $\bm{\theta} = \theta_\alpha \bm{a}^\alpha$).While alternative choices for expressing rotation are explored in the literature, they often lead to an increased number of degrees of freedom and necessitate additional drilling stabilization. Whereas in the present formulation, expressing $\bm{\theta}$ in a local basis tangent to the mid-surface allows for the use of only two coordinates, with $\theta_3 = 0$ inherently satisfying the non-extensibility condition. Consequently, Eq.\eqref{eq:RM - strain provv} can be simplified further as
\begin{subequations} \label{eq:RM - strain}
\begin{align}
    2\varepsilon_{\alpha\beta} &= \bm{a}_\alpha\cdot\bm{u}_{,\beta}+\bm{a}_\beta\cdot\bm{u}_{,\alpha} \;,\\
    2\kappa_{\alpha\beta} &= 
    (\theta_{\alpha,\beta}+\theta_{\beta,\alpha}) - (\bm{a}_{\alpha,\beta}\cdot\bm{\theta} + \bm{a}_{\beta,\alpha}\cdot\bm{\theta}) + (\bm{a}_{3,\alpha}\cdot\bm{u}_{,\beta} + \bm{a}_{3,\beta}\cdot\bm{u}_{,\alpha}) \;,\\
    \gamma_{\alpha} &= \bm{a}_3\cdot\bm{u}_{,\alpha} + \theta_\alpha \;.
\end{align}
\end{subequations}
Where it was used $\bm{a}_\alpha\cdot\bm{\theta}_{,\beta}=(\bm{a}_\alpha\cdot\bm{\theta})_{,\beta} - \bm{a}_{\alpha,\beta}\cdot\bm{\theta}=\theta_{\alpha,\beta} - \bm{a}_{\alpha,\beta}\cdot\bm{\theta}$. The weak form of the governing equation derives from the three-dimensional principle of virtual displacements. This involves integrating through the thickness, incorporating the shell kinematics, and assuming $\dd V=\dd \Omega \dd\xi_3$. While the detailed derivations are omitted here for brevity, they can be found in various references on shell mechanics, e.g., \cite{chapelle2011finite}. The resulting two-dimensional variational statement for the Reissner-Mindlin shell reads: find $(\bm{u},\bm{\theta})\in[\mathcal{H}^1]^5_{(\bm{\bar{u}},\bm{\bar{\theta}})}$ such that
\begin{equation} \label{eq:RM variational con}
    \mathcal{L}_{int}^{RM}(\delta\bm{u},\delta\bm{\theta},\bm{u},\bm{\theta}) = \mathcal{L}_{ext}^{RM}(\delta\bm{u},\delta\bm{\theta}) \;,\quad\quad\forall\delta\bm{u},\delta\bm{\theta}\in[\mathcal{H}^1]^5_{(\bm{0},\bm{0})}\;,
\end{equation}
where $[\mathcal{H}^1]^5_{(\bar{\bm{u}},\bar{\bm{\theta}})}$ denotes the subspace of $[\mathcal{H}^1]^5$ in which Dirichlet boundary conditions are satisfied, and $[\mathcal{H}^1]^5_{(\bm{0},\bm{0})}$ represents the subspace of $[\mathcal{H}^1]^5$ where the test functions are null on the corresponding Dirichlet boundary. The bilinear form $\mathcal{L}_{int}^{RM}$ and the linear form $\mathcal{L}_{ext}^{RM}$ represent the virtual work of the internal forces and external forces, respectively. Their expressions are given by
\begin{subequations} \label{eq:RM Lint Lext con}
\begin{align}
    \mathcal{L}_{int}^{RM} &= \int_{\Omega}(\delta\bm{\varepsilon}\colon\bm{N} + \delta\bm{\kappa}\colon\bm{M} + \delta\bm{\gamma}\colon\bm{Q})\dd\Omega \;,\\
    \mathcal{L}_{ext}^{RM} &= \int_\Omega(\delta\bm{u}\cdot\bar{\bm{f}} +\delta\bm{\theta}\cdot\bar{\bm{m}})\dd\Omega + \int_{\pd\Omega_{N_1}}\delta\bm{u}\cdot\bar{\bm{F}}\dd\pd\Omega + \int_{\pd\Omega_{N_2}}\delta\bm{\theta}\cdot\bar{\bm{M}}\dd\pd\Omega  \;,
\end{align}
\end{subequations}
where the operator $\colon$ denotes the contraction operation between tensors.

\subsection{The Kirchhoff-Love shell equation}
In Kirchhoff-Love theory, it is additionally assumed that the unit normal segment remains perpendicular to the deformed mid-surface. This implies that shear strain are zero, and thus
\begin{equation}
    \theta_\alpha = -\bm{a}_3\cdot\bm{u}_{,\alpha} \;.
\end{equation}
Substituting $\theta_\alpha$ into Eq.\eqref{eq:RM - strain} leads to the expressions for the generalized strain provided in \cite{guarino2024interior,li2018geometrically,benzaken2021} and not reported here for the sake of conciseness. Consequently, for the Kirchhoff-Love theory, the only variables are the Cartesian coordinates of the displacement vector. However, this approach results in a fourth-order differential equation, necessitating additional considerations. Firstly, the approximation space must be $C^1$ continuous to ensure that the displacement's derivatives, and thus the rotation, are continuous. Secondly, the Neumann boundary condition are adjusted as follows. The component $\bar{M}_n=\bar{\bm{M}}\cdot\bm{n}$ is applied on $\pd\Omega_{N_2}$, while the component $\bar{M}_t=\bar{\bm{M}}\cdot\bm{t}$ contributes to the definition of the Ersatz force $\bar{\bm{T}}$
\begin{equation}
    \bar{\bm{T}} = \left(\bar{F}_\alpha - \bar{M}_{t} b_{\alpha\beta} t^{\beta}\right) \bm{a}^\alpha +\left(\bar{F}_3 + \fpd{\bar{M}_{t}}{t}\right)\bm{a}^3 \;,
\end{equation}
where $\bar{F}_\alpha=\bm{a}_\alpha\cdot\bar{\bm{F}}$, $\bar{F}_3=\bm{a}_3\cdot\bar{\bm{F}}$, $t^\alpha=\bm{a}^\alpha\cdot\bm{t}$, and $\pd\bar{M}_{t}/\pd t$ denotes the arc-lenght derivative of $\bar{M}_t$. Here, $\bm{n}$ is the outer unit vector locally perpendicular to $\pd\Omega$ and lying on the plane locally tangent to $\Omega$, while $\bm{t}=\bm{a}_3\times\bm{n}$ is the unit vector locally tangent to $\pd\Omega$. Additionally, Kirchhoff-Love theory introduces the so-called corner forces at points of $\chi_N$, defined as
\begin{equation}
    \bar{R} = \lim_{\epsilon\to 0}{\left(\bar{M}_{t}(\bm{x}+\epsilon\bm{t}^+)-\bar{M}_{t}(\bm{x}-\epsilon\bm{t}^-)\right)} \;,
\end{equation}
where $\bm{t}^+$ and $\bm{t}^-$ are the tangent vectors at the two edges intersecting in a corner. The two-dimensional variational statement for the Kirchhoff-Love shell then reads: find $\bm{u}\in[\mathcal{H}^2]^3_{(\bm{\bar{u}},\bar{\theta}_n)}$ such that
\begin{equation} \label{eq:KL variational con}
    \mathcal{L}_{int}^{KL}(\delta\bm{u},\bm{u}) = \mathcal{L}_{ext}^{KL}(\delta\bm{u})\;,\quad\quad\forall\delta\bm{u}\in[\mathcal{H}^2]^3_{\bm{0},0}\;,
\end{equation}
where $[\mathcal{H}^2]^3_{(\bar{\bm{u}},\bar{\theta}_n)}$ is the subspace of $[\mathcal{H}^2]^3$ where Dirichlet boundary conditions are satisfied and $[\mathcal{H}^2]^3_{(\bm{0},0)}$ is the subspace of $[\mathcal{H}^2]^3$ where both $\bm{u}$ and $\theta_n$ are null on the corresponding Dirichlet boundary. The bilinear form $\mathcal{L}_{int}^{KL}$ and the linear form $\mathcal{L}_{ext}^{KL}$ are expressed in this case as
\begin{subequations} \label{eq:KL Lint Lext con}
\begin{align}
    \mathcal{L}_{int}^{KL} &= \int_{\Omega}(\delta\bm{\varepsilon}\colon\bm{N} + \delta\bm{\kappa}\colon\bm{M})\dd\Omega \;,\\
    \mathcal{L}_{ext}^{KL} &= \int_\Omega{\delta\bm{u}\cdot\bar{\bm{f}}\dd\Omega} + \int_{\pd\Omega_{N_1}}\delta\bm{u}\cdot\bar{\bm{T}}\dd\pd\Omega + \int_{\pd\Omega_{N_2}}\theta_n\bar{M}_n\dd\pd\Omega + \sum_{C\in\chi_N}{\left.\left(v_{3} \bar{R}\right)\right|_C} \;.
\end{align}
\end{subequations}

\section{Discretized shell variational statements} \label{sec:IP}

This section focuses on the numerical approximation of the variational statements in Eqs.\eqref{eq:RM variational con} and \eqref{eq:KL variational con}. First, the domain is discretized following the domain sub-structuring and meshing strategy as described in Section \ref{ssec: discretization}. Then, bi-variate spline spaces of arbitrary degree are constructed over each discretized patch and summed to construct the approximation space for the whole domain $\mathcal{S}_h$ as described in Section \ref{ssec:splines}. From here, the vector approximation spaces for the Reissner-Mindlin and Kirchhoff-Love theories are obtained, respectively, as $\mathcal{V}_h^{RM}=[\mathcal{S}_{h}]^5$ and $\mathcal{V}_h^{KL}=[\mathcal{S}_{h}]^3$. It is worth mentioning that in $\mathcal{V}_h^{KL}$, the spline bases are taken with a degree of at least 2 within the domain of each patch to satisfy the $C^1$ continuity requirement imposed by the Kirchhoff-Love shell equation.

The discretized variational statement for the Reissner-Mindlin shell equation reads: find $(\bm{u}_h,\bm{\theta}_h)\in \mathcal{V}_{h(\bar{\bm{u}},\bar{\bm{\theta}})}^{RM}$ such that
\begin{equation}\label{eq:RM variational}
    \mathcal{L}_{int}^{RM}(\delta\bm{u}_h,\delta\bm{\theta}_h,\bm{u}_h,\bm{\theta}_h) +  \mathcal{L}_{nit}^{RM}(\delta\bm{u}_h,\delta\bm{\theta}_h,\bm{u}_h,\bm{\theta}_h)    
    = \mathcal{L}_{ext}^{RM}(\delta\bm{u}_h,\delta\bm{\theta}_h)\;,\quad\quad\forall(\delta\bm{u}_h,\delta\bm{\theta}_h)\in\mathcal{V}_{h({\bm{0}},{\bm{0}})}^{RM}   
    \;,
\end{equation}
where $\mathcal{L}_{int}^{RM}$ and $\mathcal{L}_{ext}^{RM}$ are defined as in Eq.\eqref{eq:RM Lint Lext con}, $\mathcal{L}_{nit}^{RM}$ denotes the contribute arising from the Nitsche-based coupling method and is discussed in detail in the reminder of this section, $\mathcal{V}_{h(\bar{\bm{u}},\bar{\bm{\theta}})}^{RM}$ represents the subspace of $\mathcal{V}_h^{RM}$ where $\bm{u}_h=\bar{\bm{u}}$ and $\bm{\theta}_h=\bar{\bm{\theta}}$ on the corresponding Dirichlet boundaries, and $\mathcal{V}_{h(\bm{0},\bm{0})}^{RM}$ represents the subspace of $\mathcal{V}_h^{RM}$ where $\delta\bm{u}_h=\bm{0}$ and $\delta\bm{\theta}_h=\bm{0}$ on the corresponding Dirichlet boundaries. It is worth noting that boundary conditions on trial and test functions are imposed strongly by directly assigning values to the degrees of freedom associated with the basis functions that are non-zero on $\pd\Omega$ after an $L^2$ projection of the applied displacement and rotation, as customary when adopting spline-based approximation spaces \cite{hughes2005isogeometric,cottrell2009}. This is made possible by the adoption of the IBCM, which ensures the conformal nature of each boundary within the domain.

Similarly, the discretized variational statement for the Kirchhoff-Love shell equation is formulated as: find $\bm{u}_h\in\mathcal{V}_{h (\bm{\bar{u}})}^{KL}$ such that
\begin{equation}\label{eq:KL variational}
    \mathcal{L}_{int}^{KL}(\delta\bm{u}_h,\bm{u}_h) + \mathcal{L}_{nit}^{KL}(\delta\bm{u}_h,\bm{u}_h)= \mathcal{L}_{ext}^{KL}(\delta\bm{u}_h)\;,\quad\quad\forall\delta\bm{u}_h\in\mathcal{V}_{h(\bm{0})}^{KL}\;,
\end{equation}
where $\mathcal{L}_{int}^{KL}$ and $\mathcal{L}_{ext}^{KL}$ are defined as in Eq.\eqref{eq:KL Lint Lext con}, $\mathcal{L}_{nit}^{KL}$ denotes the Nitsche contribute, $\mathcal{V}_{h(\bar{\bm{u}})}^{KL}$ and $\mathcal{V}_{h({\bm{0}})}^{KL}$ represent the subspace of $\mathcal{V}_h^{KL}$ that satisfy in a strong sense the Dirichlet conditions on trial and test functions, respectively. It is worth mentioning that for the Kirchhoff-Love equation the external rotation on Dirichlet boundary still needs to be imposed in a weak sense following the Nitsche's method as described in \cite{benzaken2021}. However, this procedure is not reported here for the sake of conciseness.

\subsection{Nitsche method for Reissner-Mindlin shell coupling}
The Nitsche term in Eq.\eqref{eq:RM variational} is based on integrals along the interfaces between adjacent patches. This term is defined as
\begin{multline} \label{eq:IP RM nitsche}
    \mathcal{L}_{nit}^{RM} = \sum_{i=1}^{N_\Gamma}\int_{\Gamma_i}\biggl( -
    \jump{\delta\bm{u}_h}\cdot\ave{\bm{N}_n} -\jump{\delta\bm{\theta}_h}\cdot\ave{\bm{M}_n}  + \\
      -\gamma_1 \ave{\delta\bm{N}_n}\cdot\jump{\bm{u}_h} -\gamma_1 \ave{\delta\bm{M}_n}\cdot\jump{\bm{\theta}_h}+ \\
    +\mu^{RM}_u\jump{\delta\bm{u}_h}\cdot\jump{\bm{u}_h} + \mu^{RM}_\theta\jump{\delta\bm{\theta}_h}\cdot\jump{\bm{\theta}_h}
    \biggr)\dd\Gamma \;,
\end{multline}
where it is reminded that $\Gamma_i$ is the $i$-th of the $N_\Gamma$ interfaces. $\bm{N}_n=\bm{N}_n(\bm{u}_h,\bm{\theta}_h)$ and $\bm{M}_n=\bm{M}_n(\bm{u}_h,\bm{\theta}_h)$ are the  Reissner-Mindlin fluxes for the trial functions space, and $\delta\bm{N}_n=\delta\bm{N}_n(\delta\bm{u}_h,\delta\bm{\theta}_h)$ and $\delta\bm{M}_n=\delta\bm{M}_n(\delta\bm{u}_h,\delta\bm{\theta}_h)$ are the  Reissner-Mindlin fluxes for the test functions space. $\gamma_1$, $\mu_u^{RM}$, and  $\mu_\theta^{RM}$ are parameters of the Nitsche method and their choice determine its type and properties, as discussed in Section \ref{ssec:IP parameters}. The operators $\ave{\bullet}$ and $\jump{\bullet}$ are the average and jump operators defined in Section \ref{ssec:IP parameters}.

Following \cite{schollhammer2019reissner}, the fluxes for the Reissner-Mindlin shell in the outer unit direction $\bm{n}$ introduced in Eq.\eqref{eq:IP RM nitsche} are of two types: the force $\bm{N}_n$ and moment $\bm{M}_n$, defined as follows
\begin{subequations}
\begin{align}
    \bm{N}_n &= N^\alpha\bm{a}_\alpha + Q\bm{a}_3 \;,\\
    \bm{M}_n &= M_{nn}\bm{n} + M_{nt} \bm{t} \;,
\end{align}
\end{subequations}
where $N^\alpha$ and $Q$ are computed as
\begin{subequations}
\begin{align}
    N^\alpha &= (N^{\alpha\beta} - b^\alpha_\gamma M^{\gamma\beta})n_\beta \;, \\
    Q        &= Q^{\alpha}   n_\alpha \;.
\end{align}
\end{subequations}
Finally, the components for the moments $\bm{M}_n$ are obtained from
\begin{subequations}
\begin{align}
    M_{nn} &= M^{\alpha\beta} n_\alpha n_\beta \;, \\
    M_{nt} &= M^{\alpha\beta} n_\alpha t_\beta \;.
\end{align}
\end{subequations}

\subsection{Nitsche method for Kirchhoff-Love shell coupling}
Similarly to the Reissner-Mindlin case, the Nitsche term in Eq.\eqref{eq:KL variational} is defined as
\begin{multline} \label{eq:IP KL nitsche}
    \mathcal{L}_{nit}^{KL} = \sum_{i=1}^{N_\Gamma}\int_{\Gamma_i}\biggl( -
    \jump{\delta\bm{u}_h}\cdot\ave{\bm{T}_n} -\jump{\delta{\theta}_n}\cdot\ave{{M}_{nn}}  + \\
     -\gamma_1 \ave{\delta\bm{T}_n}\cdot\jump{\bm{u}_h} -\gamma_1 \ave{\delta{M}_{nn}}\cdot\jump{{\theta}_n}  + \\
    +\mu^{KL}_u\jump{\delta\bm{u}_h}\cdot\jump{\bm{u}_h} + \mu^{KL}_\theta\jump{\delta{\theta}_n}\cdot\jump{{\theta}_n}
    \biggr)\dd\Gamma \;,
\end{multline}
where $\bm{T}_n=\bm{T}_n(\bm{u}_h)$ and $\bm{M}_n=\bm{M}_n(\bm{u}_h)$ are the  Kirchhoff-Love fluxes for the trial functions space, and $\delta\bm{T}_n=\delta\bm{T}_n(\delta\bm{u}_h)$ and $\delta\bm{M}_n=\delta\bm{M}_n(\delta\bm{u}_h)$ are the  Kirchhoff-Love fluxes for the test functions space, as defined in the reminder of this section. It is worth reminding that in Kirchhoff-Love theory
\begin{equation}
    \theta_n=\theta_n(\bm{u}_h)=\theta_\alpha(\bm{u}_h) n^\alpha \;,
\end{equation}
where $n^\alpha=\bm{n}\cdot\bm{a}^\alpha$. Similarly, for the test functions, $\delta\theta_n=\delta\theta_n(\delta\bm{u}_h)$.

The correct expression of the fluxes for the Kirchhoff-Love formulation was only recently found in \cite{benzaken2021} and it's partially reported here. Also for the Kirchhoff-Love theory the fluxes are of two types: the so-called ersatz forces $\bm{T}_n$ and the bending moment $M_{nn}$. While the bending moment is defined in the same way as in the Reissner-Mindlin formulation, the expression for the ersatz forces is
\begin{equation}
    \bm{T}_n = T^\alpha\bm{a}_\alpha+T^3\bm{a}_3 \;,
\end{equation}
where $T^\alpha$ and $T^3$ are given by 
\begin{subequations}
    \begin{align}
        &T^\alpha = N^{\alpha\beta}n_\beta - b^\alpha_\gamma M^{\gamma\beta}n_\beta - M_{nt} b^\alpha_\gamma t^\gamma \;,\\
        &T^3 = M^{\alpha\beta}_{|\beta}n_\alpha + (M^{\alpha\beta}n_{\alpha}t_\beta)_{,t} \;.
    \end{align}
\end{subequations}
In the previous equation the notation $M^{\alpha\beta}_{|\beta}$ denotes the covariant derivative in the direction $\beta$ of the $\alpha\beta$-th component of the bending moment tensor, whereas the notation $(\bullet)_{,t}$ refers to the arc-length derivative along $\Gamma_i$. The details on how to compute these two derivatives are not reported here for brevity but can be found in \cite{benzaken2021,guarino2024interior}. 

\subsection{Choice of the method's parameters} \label{ssec:IP parameters}
In both Eqs.\eqref{eq:IP RM nitsche} and \eqref{eq:IP KL nitsche}, the Nitsche term relies on integrals along the interfaces between patches. These integrals consist of six additive terms each. The first two terms are the consistency ones introduced to restore consistency in the variational formulation. The average and jump operators that appear in these terms are defined as
\begin{align} \label{eq:IP ave - jump}
    \ave {\bullet} &= \gamma_2\bullet^++(1-\gamma_2)\bullet^- \;,\\
    \jump{\bullet} &= \bullet^+-\bullet^- \;,
\end{align}
where $\bullet^+$ and $\bullet^-$ denote a generic quantity computed from each side of the interface. It is worth mentioning that to compute the average of the same quantity, the fluxes from the second patch are computed relatively to a unit vector $\bm{n}$ entering the patch domain. The value of $\gamma_2$ ranges between 0 and 1 and determines how the average operator leans toward the first or the second patch.

The second and third terms in Eqs.\eqref{eq:IP RM nitsche} and \eqref{eq:IP KL nitsche} are called symmetric terms. Typical choices for the value of $\gamma_1$ that multiplies these terms are $\{0, +1, -1\}$, producing non-symmetrical, symmetrical, and anti-symmetrical contributions to the solving linear system, respectively. The last two terms are called the stabilization terms and are introduced to ensure coercivity of the bilinear form and, therefore, the stability of the method. The multiplicative terms $\mu_u$ and $\mu_\theta$ are called penalty parameters, and their choice balances two opposite needs: ensuring stability, and limiting the condition number of the linear system. It is worth mentioning that the pure penalty method relies only on these terms to restore coupling, leading to a symmetric and stable formulation \cite{breitenberger2015, herrema2019penalty, leonetti2020robust, zhao2022opensource, prosepio2022penalty}. However the values of the penalty parameters to achieve a given error need to be considerably higher given the lack of consistency \cite{babuska1973penalty, pash2021priori}.

Accordingly to the choice of $\gamma_1$, $\gamma_2$, and the penalty values, different Nitsche-type formulations can be constructed. Selecting $\gamma_1=1$ leads to the symmetric Nitsche method \cite{douglas2008interior, guo2015, nguyen2017isogeometric, guo2018variationally, guo2021isogeometric, chasapi2023fast,guarino2024interior}, also known as interior penalty method. This formulation strictly requires the presence of the penalty terms, and therefore $\mu>0$, to ensure coercivity.  With $\gamma_1=-1$ and $\mu=0$ the Nitsche contribute is anti-symmetrical and leads to a more stable formulation \cite{bauman1999discontinuous, hu2018skew, guo2017parameter}, where penalty terms can still be adopted to further increase stability, loosing however every type of symmetry \cite{riviere1999improved, wang2022isogeometric, yu2023isogeometric}. The value of $\gamma_2$ can be used to balance the computation of the average of the fluxes accordingly to the different properties of the two side of the interface \cite{annavarapu2012, jiang2015robust}. 

Adequately choosing these parameters to construct a stable method is even more difficult when the elements in the patches to be coupled are trimmed by the interface. In limit situations in order to ensure coercivity, in the symmetric Nitsche method, a local eigenvalue problem needs to be addressed and the value of the penalty parameter can grow unbounded \cite{griebel2003,ruess2014weak,guo2017parameter}.

In this work, the symmetric Nitsche method is adopted for coupling internal patches and boundary layers. The stability of the method relies on the definition of the average operator and, therefore, the choice for the parameter $\gamma_2$. Here, the fluxes are computed uniquely from the boundary layer, setting $\gamma_2=1$ in Eq.\eqref{eq:IP ave - jump} and leveraging the fact that the elements in this layer are not trimmed, thus leading to Nitsche contributes seamless to stabilize \cite{antolin2021overlapping, wei2021}. Regarding the penalty values, these are chosen as
\begin{subequations}
\begin{align}
    &\mu_u^{RM} = \mu_u^{KL} = \beta E \tau   /h  \;,\\
    &\mu_\theta^{RM}= \mu_\theta^{KL} = \beta E \tau^3 /h  \;,
\end{align}
\end{subequations}
where $E$ represent the maximum Young modulus of the layers of the laminate, $\beta$ is an arbitrary value that for Nitsche based formulations typically ranges in $[10,10^3]$, and $h$ is a measure of the dimension of the elements of the boundary layer.

\section{Results} \label{sec:RESULTS}

The proposed method is tested in this section through a series of scenarios involving shell structures. Specifically, five tests are conducted: a laminated plate with circular cut-outs, a hyperbolic-paraboloid shell with two internal cut-outs, a structure comprising two intersecting cylinders, a generally curved spline-based surface with a cut-out and a trimmed external boundary, and a cylindrical shell with an axial through the thickness fracture. The constitutive laws in these tests involve either laminated materials with orthotropic layers or single-layer isotropic materials.

In these tests, Dirichlet boundary conditions are applied using different strategies depending on whether they pertain to displacements or rotations, the shell theory in question, and whether the boundary is conformal or trimmed. Trimmed boundaries are included solely for comparison purposes with the present IBCM approach. For conformal boundaries, essential boundary conditions on displacement are applied strongly, while on rotation this applies only for the Reissner-Mindlin shell theory, where the components of the rotation vector are explicitly available as degrees of freedom. Conversely, for the Kirchhoff-Love theory, rotation boundary conditions are enforced weakly using Nitsche's method, which is also adopted for all types of Dirichlet conditions on trimmed boundaries. As standard for spline-based approximation spaces, imposing Dirichlet conditions strongly requires an $L^2$ projection of the applied external function into the B-spline space of the corresponding edge in order to determine the values of the associated degrees of freedom. Coupling conditions between boundary layers and the internal patch are imposed weakly using Nitsche's formulation, as shown in Eqs.\eqref{eq:RM variational} and \eqref{eq:KL variational}.

A common issue with Nitsche-based coupling and boundary conditions is that the resulting bilinear form might lose coercivity if not properly stabilized. In the presented tests, coercivity is verified by assessing the positive definiteness of the stiffness matrix. Specifically, this is established through a Cholesky decomposition algorithm, which fails if the stiffness matrix is severely ill-conditioned or not positive definite.

\begin{figure}	\centering
    \includegraphics[width = 0.60\textwidth]{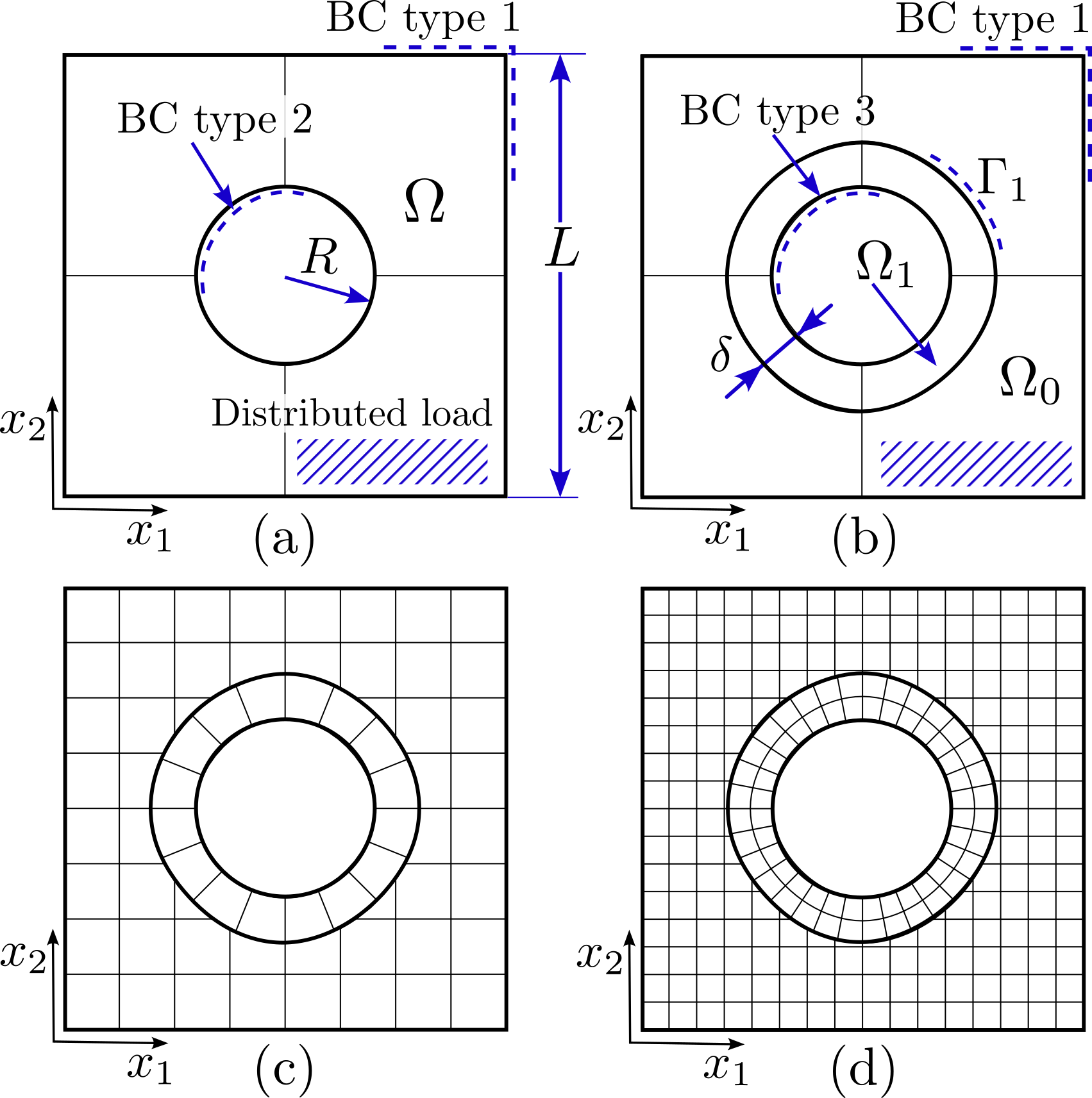}
    \begin{subfigure}{0\textwidth}\phantomcaption\label{fig:RES - LamGeom a}\end{subfigure}
    \begin{subfigure}{0\textwidth}\phantomcaption\label{fig:RES - LamGeom b}\end{subfigure}
    \begin{subfigure}{0\textwidth}\phantomcaption\label{fig:RES - LamGeom c}\end{subfigure}
    \begin{subfigure}{0\textwidth}\phantomcaption\label{fig:RES - LamGeom d}\end{subfigure}
    \caption{Geometry and discretization for the laminated plate in the first set of tests. (a) First refinement level for the single patch discretization  with corresponding types of boundary conditions and applied loads. (b) First refinement level for the IBCM-based discretization with corresponding types of boundary conditions and applied loads. Third (c) and fourth (d) refinement levels for the IBCM-based discretization.} \label{fig:RES - LamGeo}
\end{figure}

\subsection{Laminated plate with cut-out}\label{ssec:Res - lam}
The first set of tests focuses on a square laminated plate with an internal circular cut-out, as shown in Fig.\figref{fig:RES - LamGeo}. The plate has a length of $L=1$ [m], and the internal hole is centered on the square with a radius of $R=0.2$ [m]. Regarding the constitutive properties, the layup consists of four layers of an orthotropic material having Young moduli $E_1=25$ [GPa], $E_2=1$ [GPa], Poisson's ratio $\nu_{12}=0.25$, and Shear moduli  $G_{12}=G_{31}=G_{32}=0.4$ [GPa]. The layers are oriented in a sequence of $[0,90,90,0]$, with the angles measured relative to the $\xi_1$ axis. 

The external applied force and boundary conditions are chosen to reproduce the manufactured solution
\begin{equation}
    \bm{u}_{ex} = \mathrm{U}_{i} \sin(\pi n_1\xi_1) \sin(\pi n_2\xi_2) \bm{e}_i \;,
    \end{equation}
where $n_1=n_2=2$ are the number of half-waves in the $\xi_1$ and $\xi_2$ directions, respectively. The amplitude of the displacement functions are $\mathrm{U}_1=\mathrm{U}_2=\mathrm{U}_3=0.1$ [m]. For the Kirchhoff-Love theory, no additional information is required to construct the manufactured solution. However, for Reissner-Mindlin kinematic, additional fields are introduced for the rotation vector
\begin{equation}
    \bm{\theta}_{ex} = \Theta_{\alpha} \sin(\pi n_1\xi_1) \sin(\pi n_2\xi_2) \bm{e}_\alpha \;,
\end{equation}
where $\Theta_1=\Theta_2=0.1$. 

Due to the conformal nature of the external boundary of the plate, boundary conditions on the displacements are applied strongly (denoted as boundary conditions of type 1 in Fig.\figref{fig:RES - LamGeom a}) for both Reissner-Mindlin and Kirchhoff-Love elements. For the rotation, this is imposed strongly only for the Reissner-Mindlin theory, whereas, for the Kirchhoff-Love theory, homogeneous Neumann moment boundary conditions arise for the selected displacement function. The boundary conditions on the internal hole depend on the discretization strategy. Two discretization strategies are compared: a single patch with non-conformal trimmed boundaries (shown in Fig.\figref{fig:RES - LamGeom a}) and an IBCM-based discretization (shown in Fig.\figref{fig:RES - LamGeom b}) with an additional conformal boundary layer characterized by an offset distance $\delta=0.1$ [m]. 

Boundary conditions of type 2 are applied on the internal boundary of the first discretization, meaning that non-homogeneous Dirichlet boundary conditions are imposed weakly through the Nitsche's method. Due to the conformal nature of the internal boundary in the second discretization, Dirichlet conditions on the displacements can be imposed strongly. However, Dirichlet conditions for the rotation in Kirchhoff-Love elements still need to be imposed weakly.

Starting from the discretizations shown in Figs.\figref{fig:RES - LamGeom a} and \figref{fig:RES - LamGeom b}, subsequent meshing levels are obtained through dyadic refinements. However, to generate elements with comparable size in both directions, in the coordinate $\eta_2$, which refers to the offset direction of the boundary layer, the first dyadic refinement occurs when the dimension of the elements in the direction $\xi_1$ has reached values comparable to the offset distance. To exemplify, Figs.\figref{fig:RES - LamGeom c} and \figref{fig:RES - LamGeom d} show the third and fourth refinement levels for the IBCM discretization.

In computing the penalty values, the arbitrary parameter $\beta$ is chosen as $\beta=10$. The characteristic mesh size for coupling conditions in the IBCM discretization is computed as $h=\mathrm{min}(1/2^{n_{ref}},\delta)$, where $n_{ref}$ is the refinement level. For the Nitsche boundary conditions in the trimmed single patch case, $h=1/2^{n_{ref}}$.

\begin{figure}	\centering
    \begin{subfigure}[b]{0.32\textwidth}  \includegraphics[width = 0.99\textwidth]{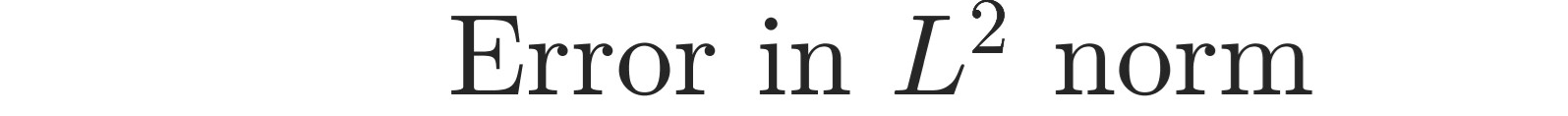}
    \end{subfigure}
    \begin{subfigure}[b]{0.32\textwidth}  \includegraphics[width = 0.99\textwidth]{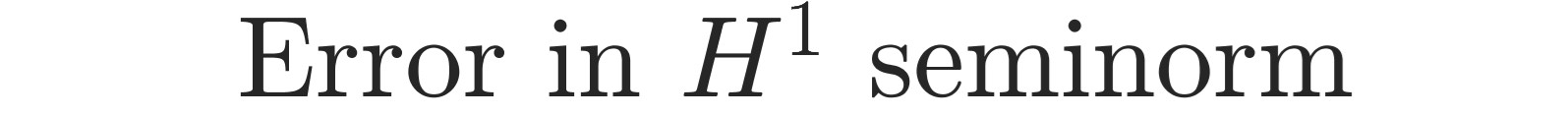}
    \end{subfigure}\\  
    \begin{subfigure}[b]{0.32\textwidth}  \includegraphics[width = 0.99\textwidth]{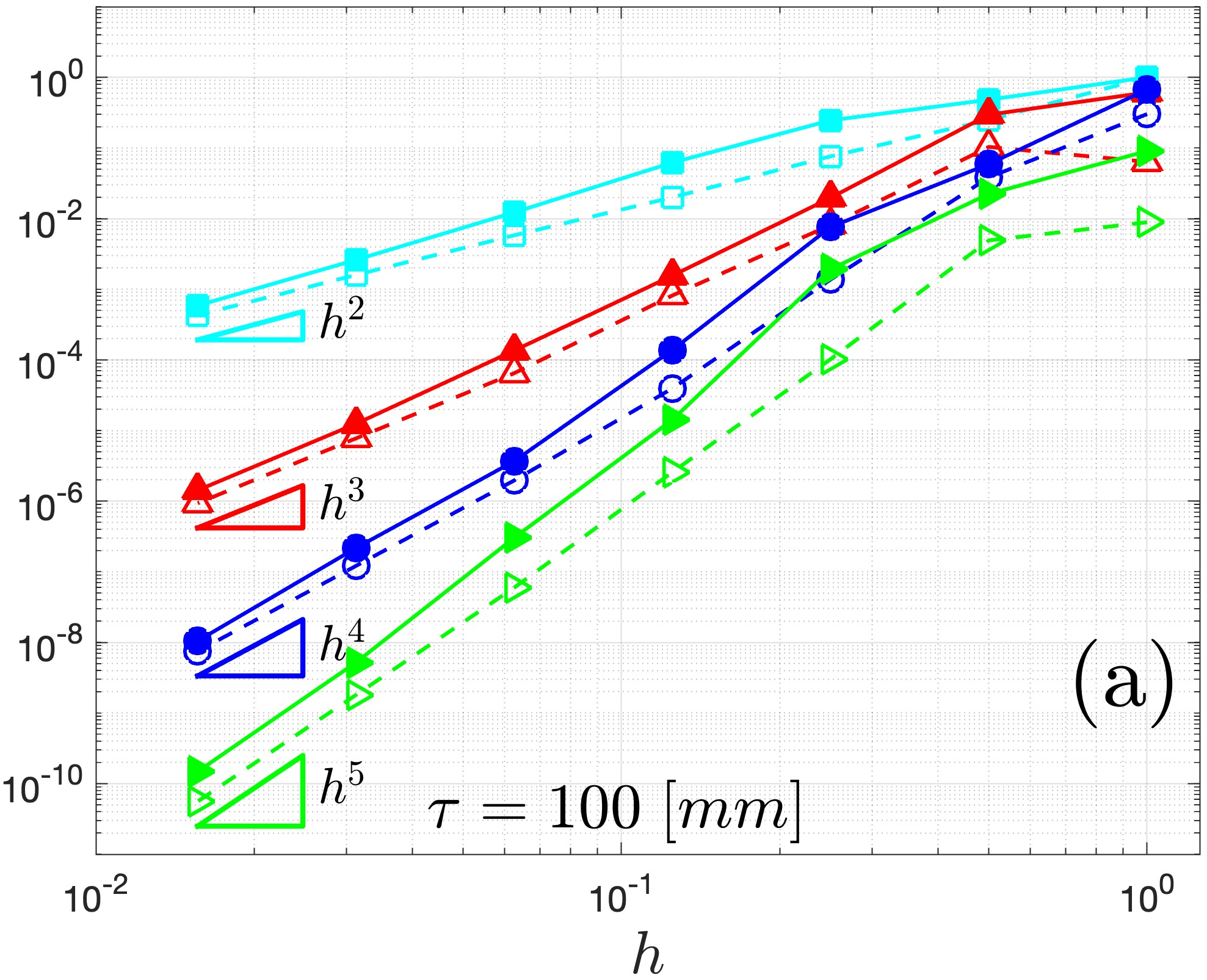}
    \end{subfigure}
    \begin{subfigure}[b]{0.32\textwidth}  \includegraphics[width = 0.99\textwidth]{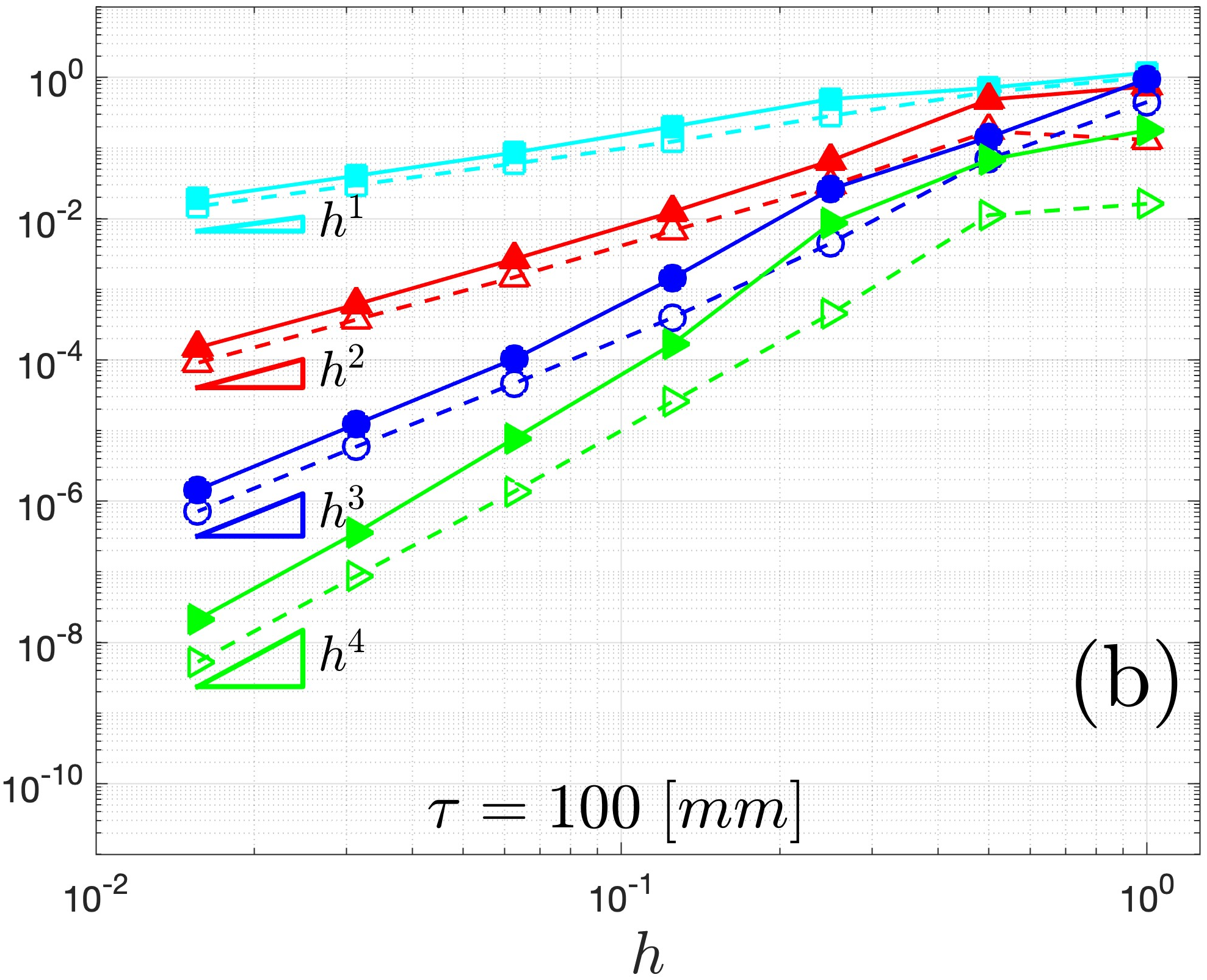}
    \end{subfigure}\\
    \begin{subfigure}[b]{0.32\textwidth}  \includegraphics[width = 0.99\textwidth]{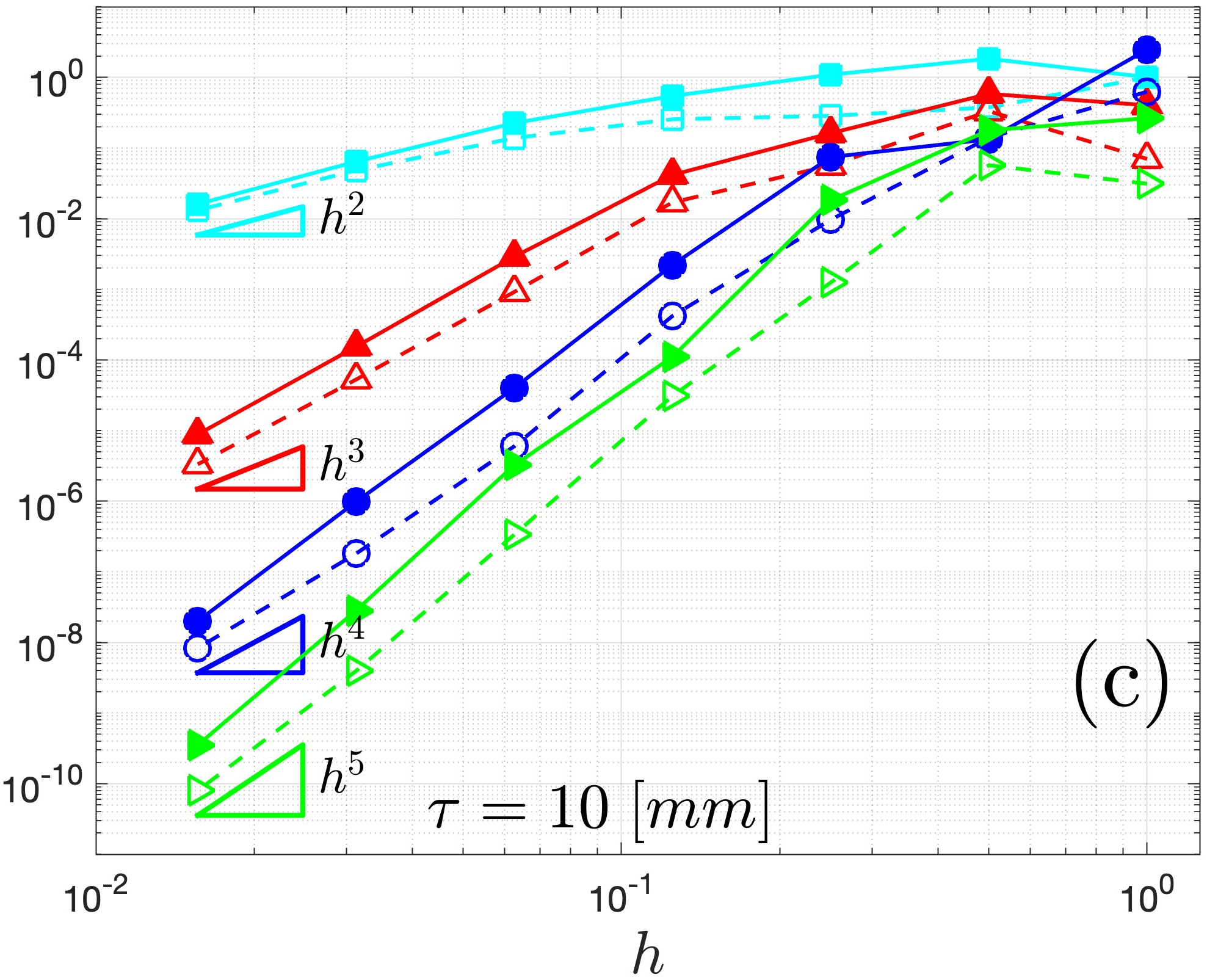}
    \end{subfigure}
    \begin{subfigure}[b]{0.32\textwidth}  \includegraphics[width = 0.99\textwidth]{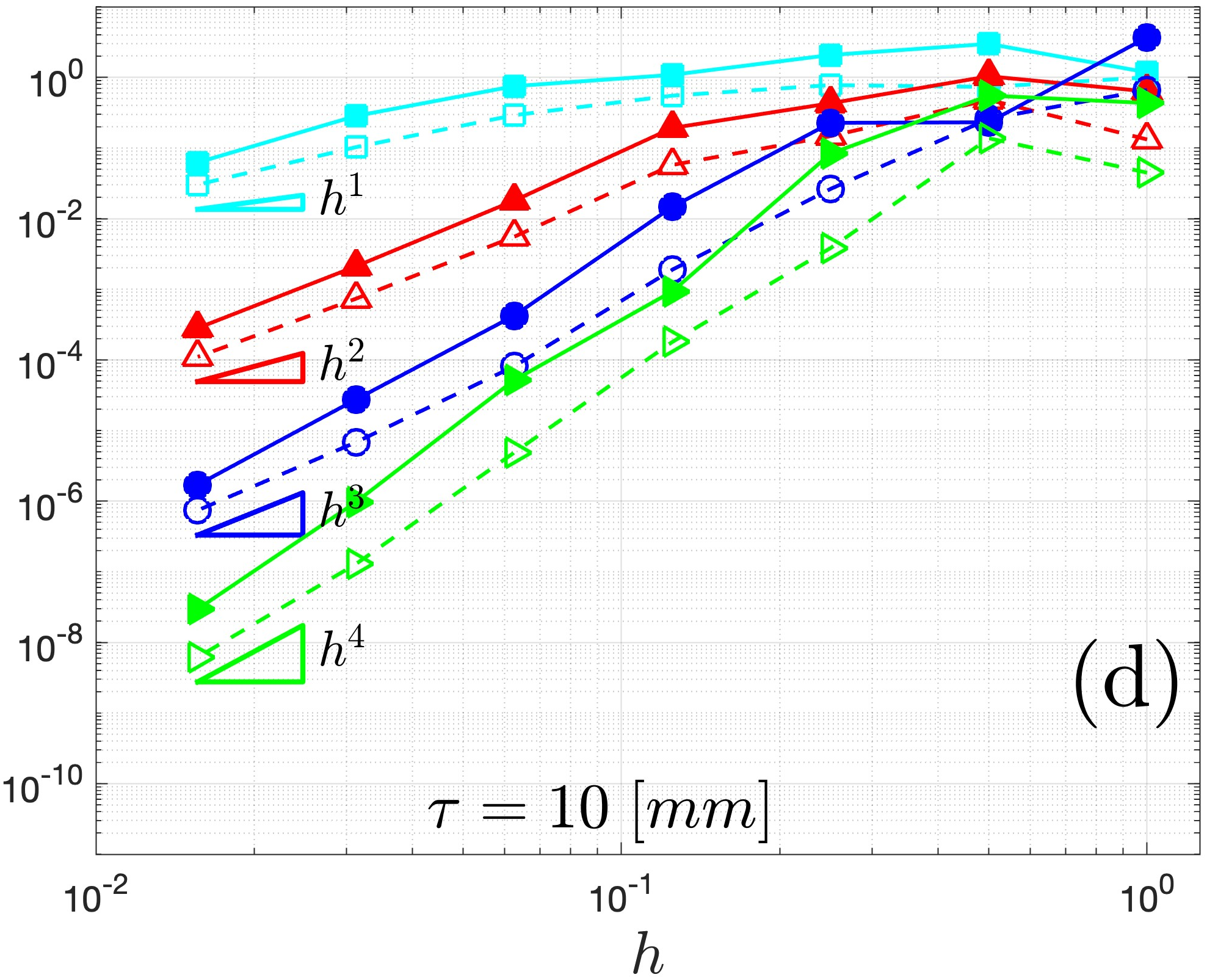}
    \end{subfigure}\\
    \begin{subfigure}[b]{0.32\textwidth}  \includegraphics[width = 0.99\textwidth]{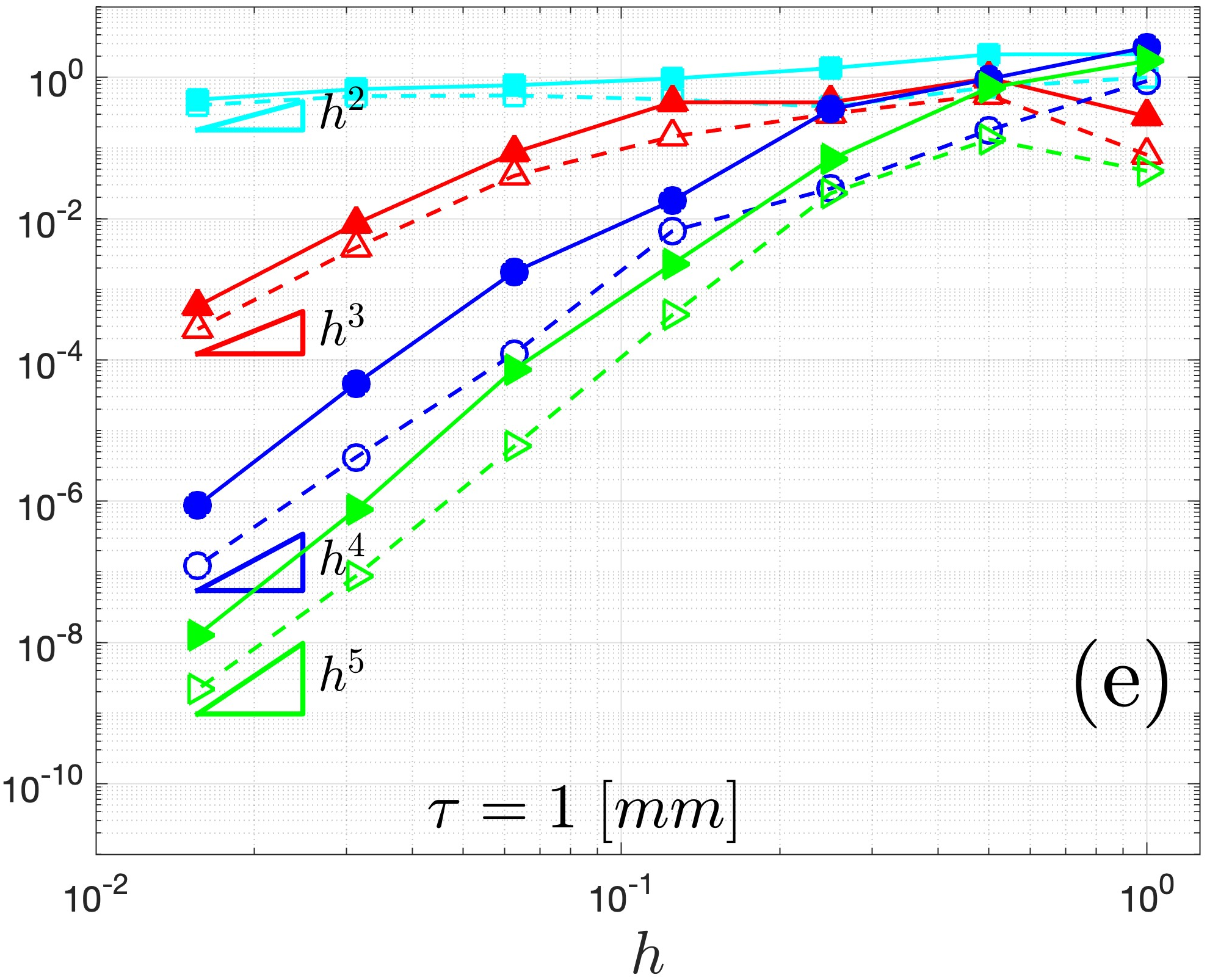}
    \end{subfigure}
    \begin{subfigure}[b]{0.32\textwidth}  \includegraphics[width = 0.99\textwidth]{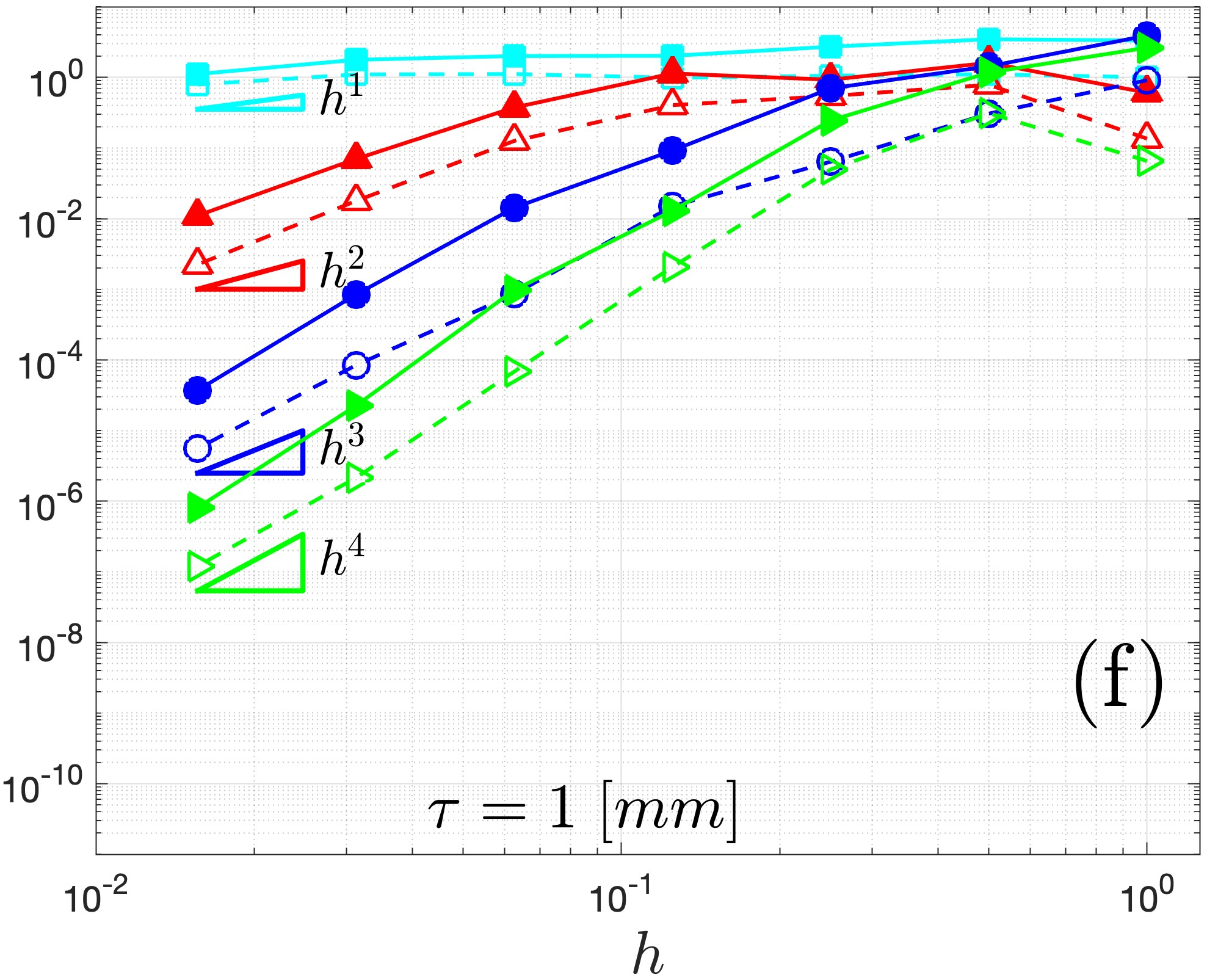}
    \end{subfigure}\\
    \begin{subfigure}[b]{0.66\textwidth}   \includegraphics[width = 0.99\textwidth]{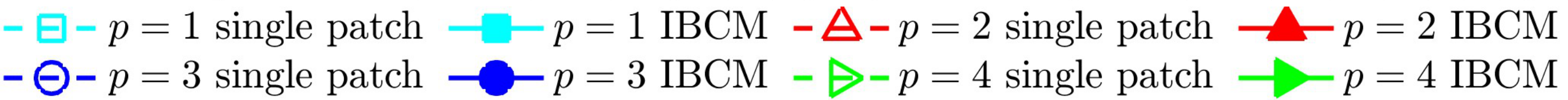}
    \end{subfigure}
	\caption{Convergence curves for the plate shown in Fig.\figref{fig:RES - LamGeo} in $L^2$ error norm and $H^1$ error seminorm for a Reissner-Mindlin theory. The curves are obtained for four different polynomial orders $p=1,2,3,4$ and three thickness values $\tau=100,10,1$ [mm]. Two discretization are taken into account, a single trimmed patch as shown in Fig.\figref{fig:RES - LamGeom a} and a IBCM-based one as shown in Fig.\figref{fig:RES - LamGeom b}.} \label{fig:RES - lam_RM_Conv}
\end{figure}

Fig.\figref{fig:RES - lam_RM_Conv} shows the convergence of the error in $L^2$ norm and $H^1$ seminorm for the Reissner-Mindlin theory. Three thickness values are considered: $\tau=100$ [mm], $\tau=10$ [mm], and $\tau=1$ [mm]. The plots show convergence curves for different polynomial orders $p=1,2,3,4$, and for the two proposed discretization strategies. The mesh size reported on the abscissa refers to a global value computed as the common edge length of the untrimmed element of the main patch. The expected optimal convergence rates are represented by auxiliary triangles below each convergence curve. The first observation is that asymptotic optimal convergence is reached in all cases except for first-degree polynomials for $\tau=1$ [mm]. In fact, as the thickness decreases, increasingly severe shear locking can be identified in the curves for all polynomials, although higher degree curves are less affected and recover proper convergence after approximately two refinements.

The second observation is that the IBCM discretization tends to have a slightly higher error value than the single patch discretization. This is expected due to the higher discrepancy between the approximation space and the manufactured target function, which is a bi-sinusoidal in $\xi_1$ and $\xi_2$. The curvilinear coordinates $\eta_1$ and $\eta_2$ of the boundary layer are not aligned with $\xi_1$ and $\xi_2$ leading to shape functions that result in higher errors. However, it is crucial to mention that for the single patch discretization, a loss of coercivity is observed for each polynomial value at the second-to-last refinement level, and for polynomials from 2 to 4 at the last refinement level. In contrast, the IBCM discretization, with the boundary layer allowing Nitsche's fluxes to be computed from the non-trimmed side, leads to a stable formulation in all refinement levels.

\begin{figure}	\centering
    \begin{subfigure}[b]{0.32\textwidth}  \includegraphics[width = 0.99\textwidth]{RES_ConvergenceTitle_L2.png}
    \end{subfigure}
    \begin{subfigure}[b]{0.32\textwidth}  \includegraphics[width = 0.99\textwidth]{RES_ConvergenceTitle_H1.png}
    \end{subfigure}
    \begin{subfigure}[b]{0.32\textwidth}  \includegraphics[width = 0.99\textwidth]{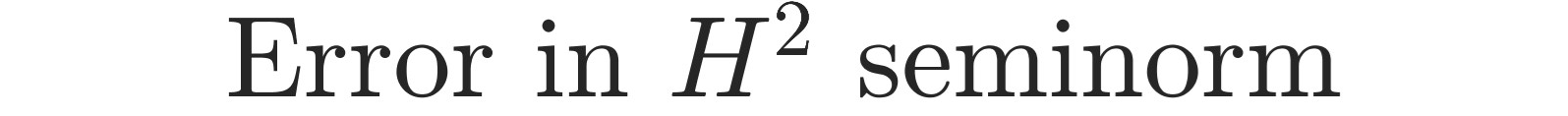}
    \end{subfigure}\\  
    \begin{subfigure}[b]{0.32\textwidth}  \includegraphics[width = 0.99\textwidth]{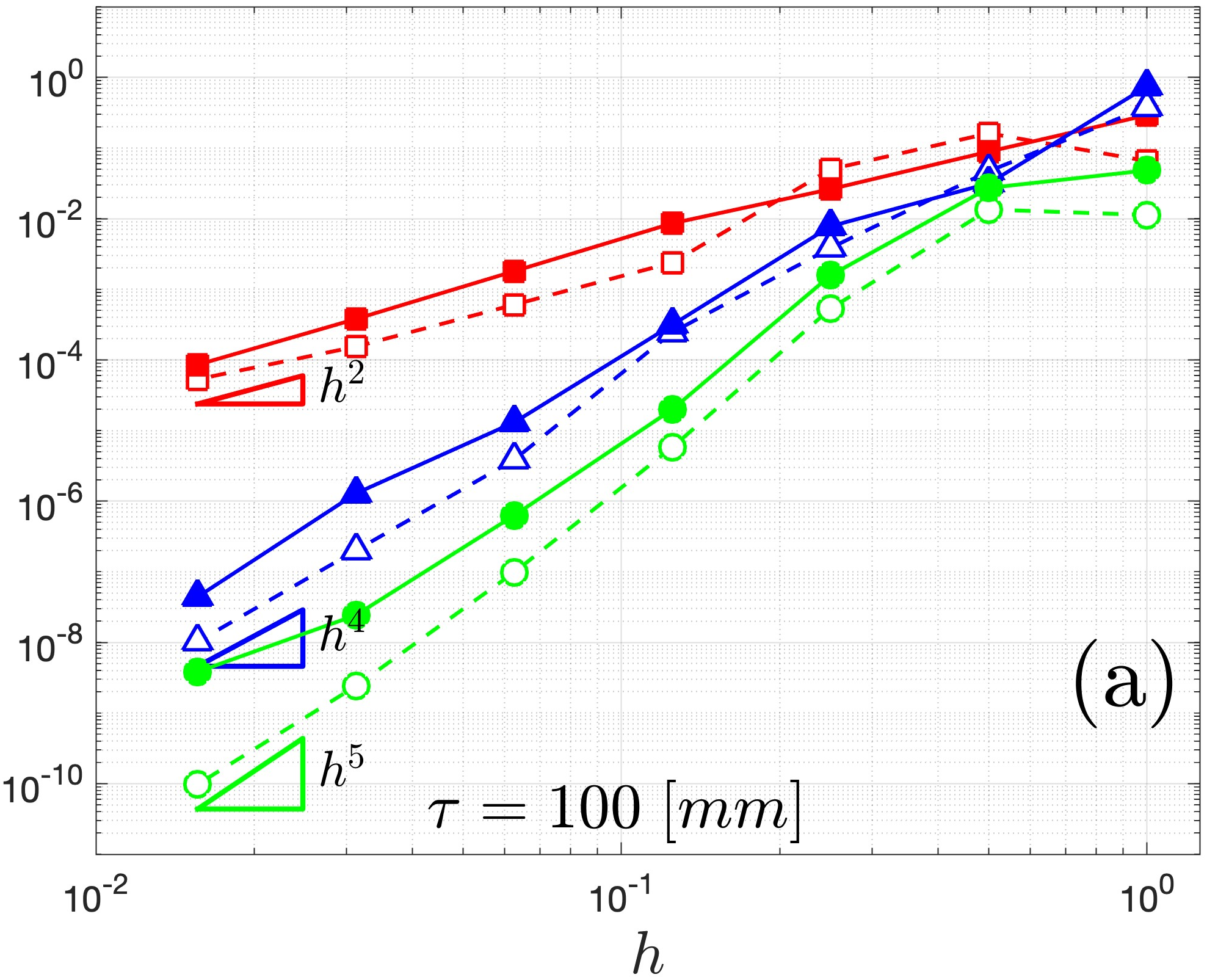}
    \end{subfigure}
    \begin{subfigure}[b]{0.32\textwidth}  \includegraphics[width = 0.99\textwidth]{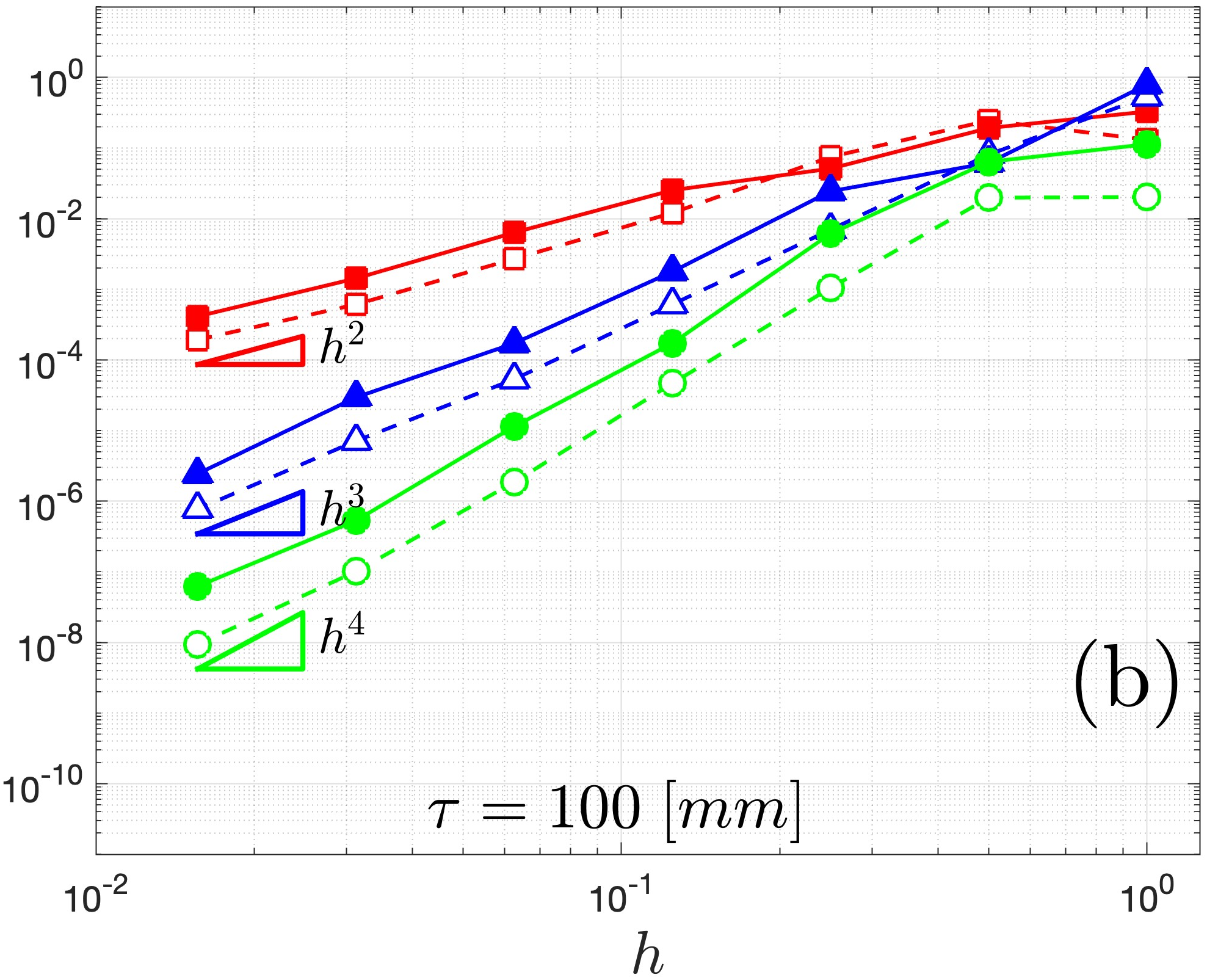}
    \end{subfigure}
    \begin{subfigure}[b]{0.32\textwidth}  \includegraphics[width = 0.99\textwidth]{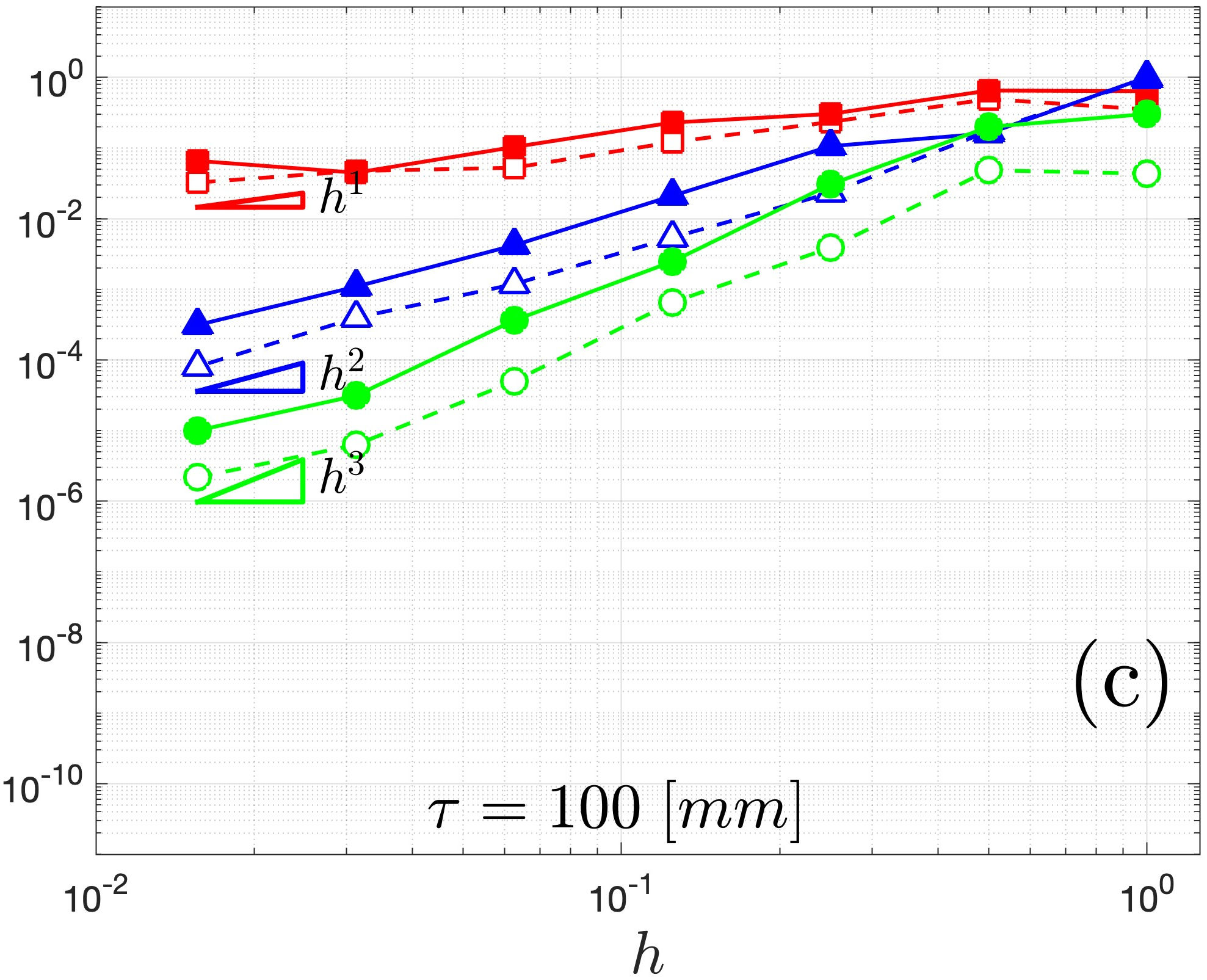}
    \end{subfigure}\\
    \begin{subfigure}[b]{0.32\textwidth}  \includegraphics[width = 0.99\textwidth]{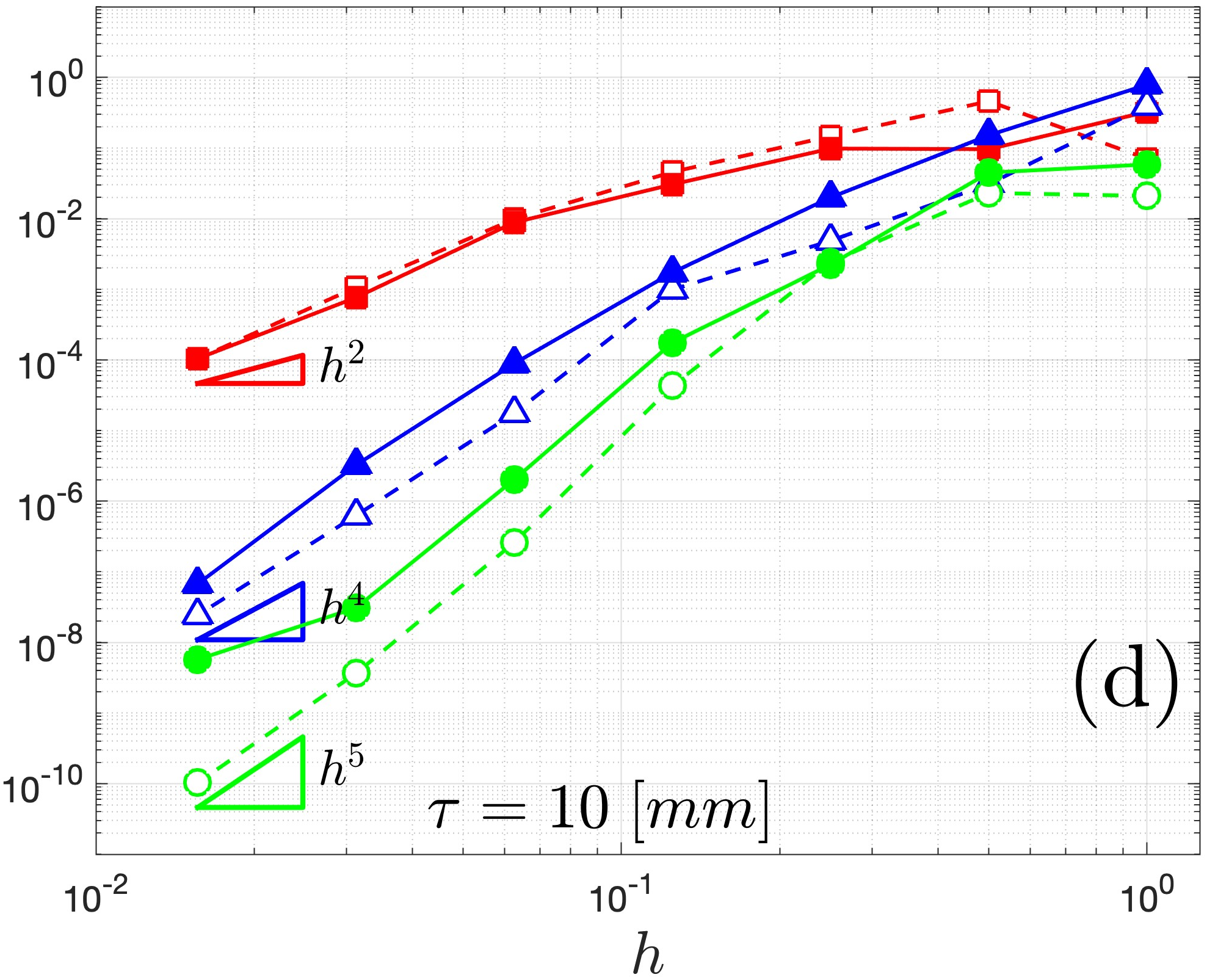}
    \end{subfigure}
    \begin{subfigure}[b]{0.32\textwidth}  \includegraphics[width = 0.99\textwidth]{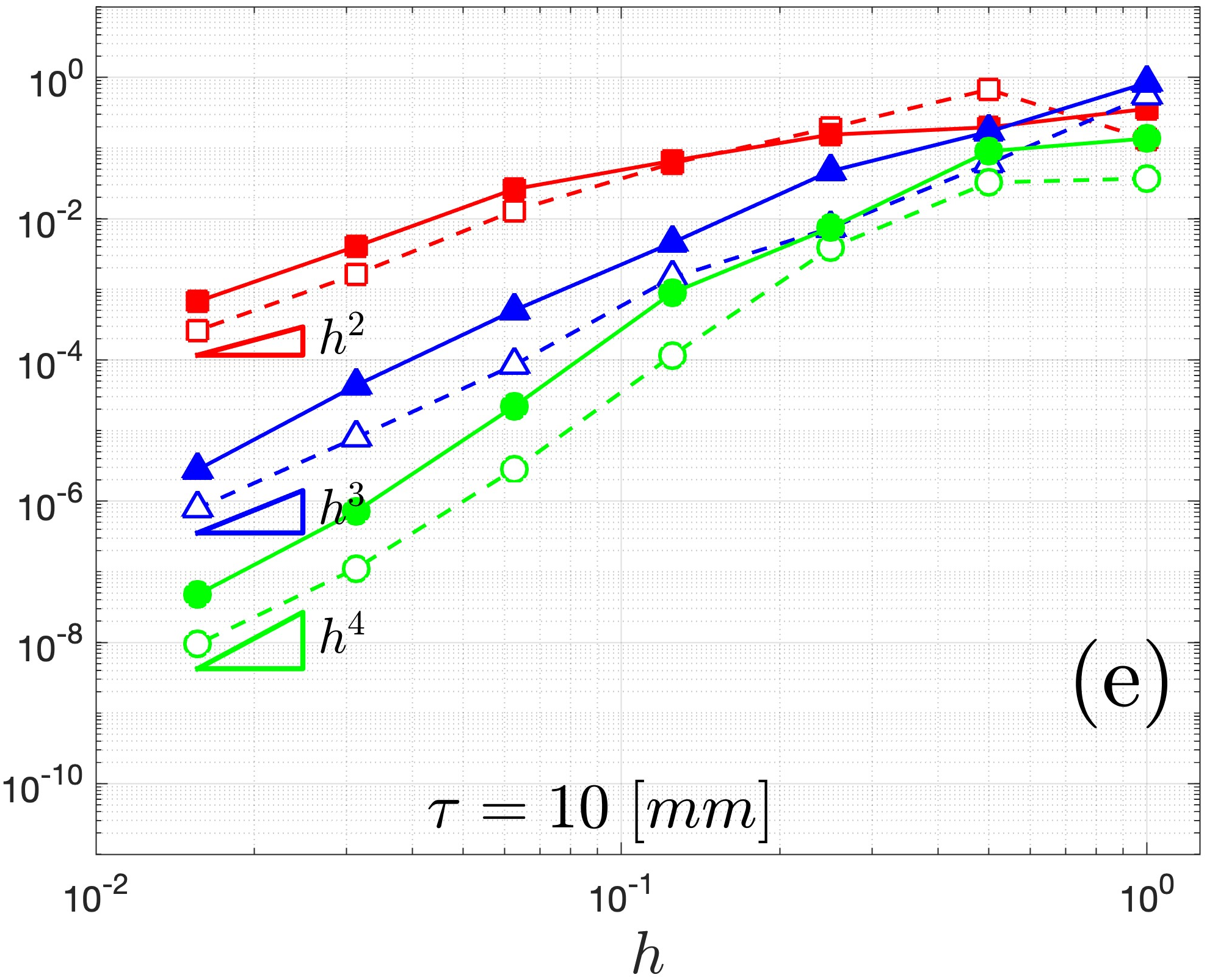}
    \end{subfigure}
    \begin{subfigure}[b]{0.32\textwidth}  \includegraphics[width = 0.99\textwidth]{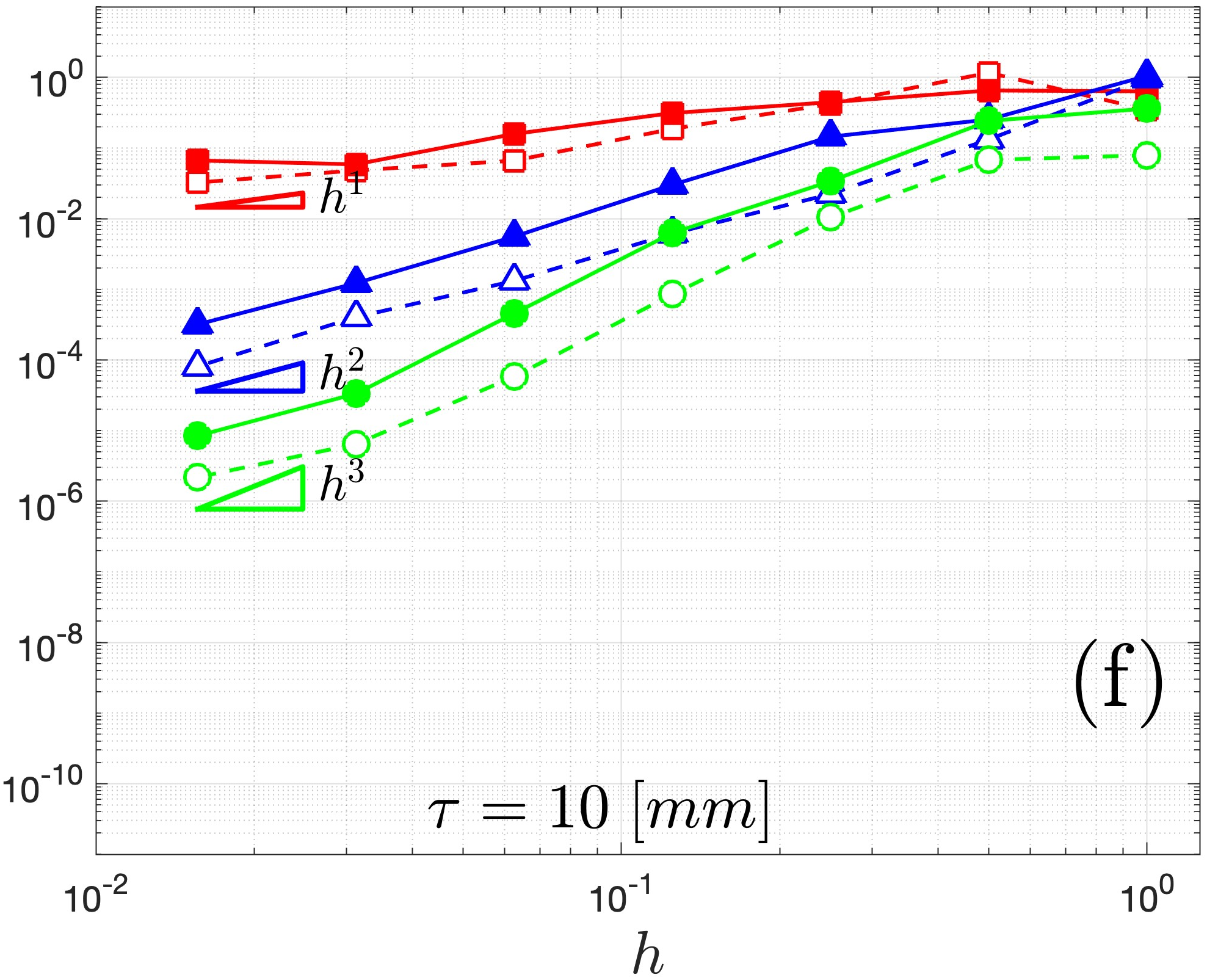}
    \end{subfigure}\\
    \begin{subfigure}[b]{0.32\textwidth}  \includegraphics[width = 0.99\textwidth]{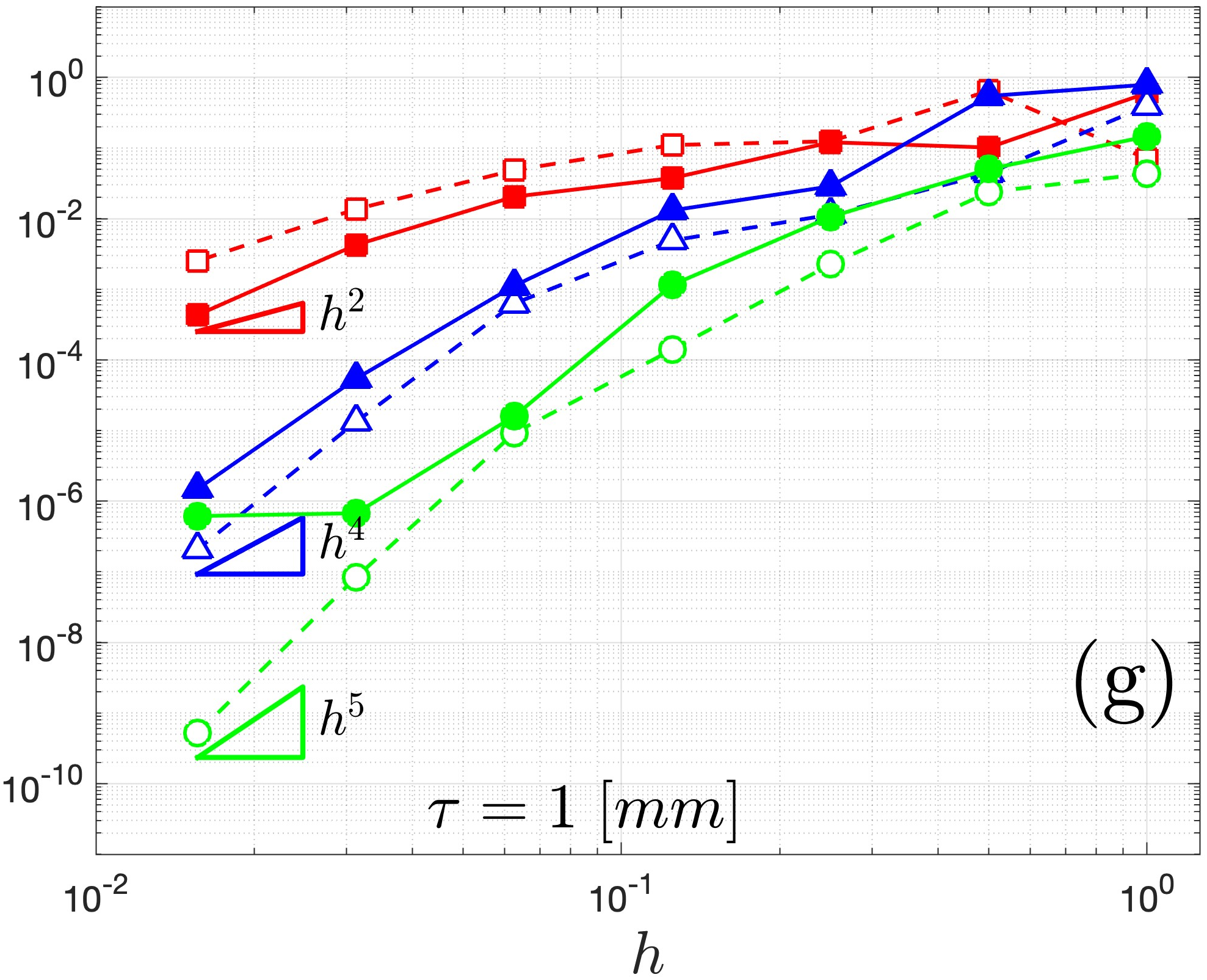}
    \end{subfigure}
    \begin{subfigure}[b]{0.32\textwidth}  \includegraphics[width = 0.99\textwidth]{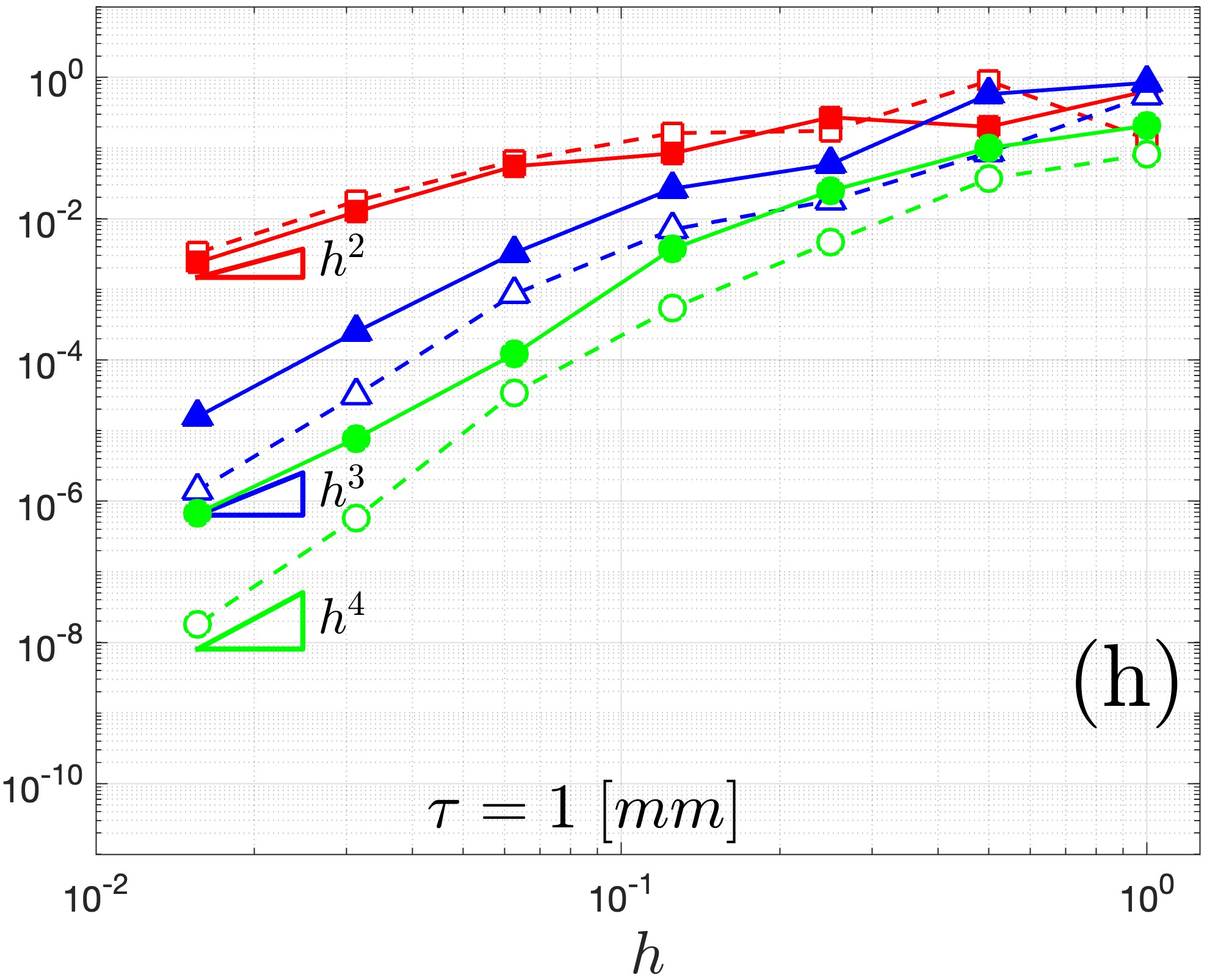}
    \end{subfigure}
    \begin{subfigure}[b]{0.32\textwidth}  \includegraphics[width = 0.99\textwidth]{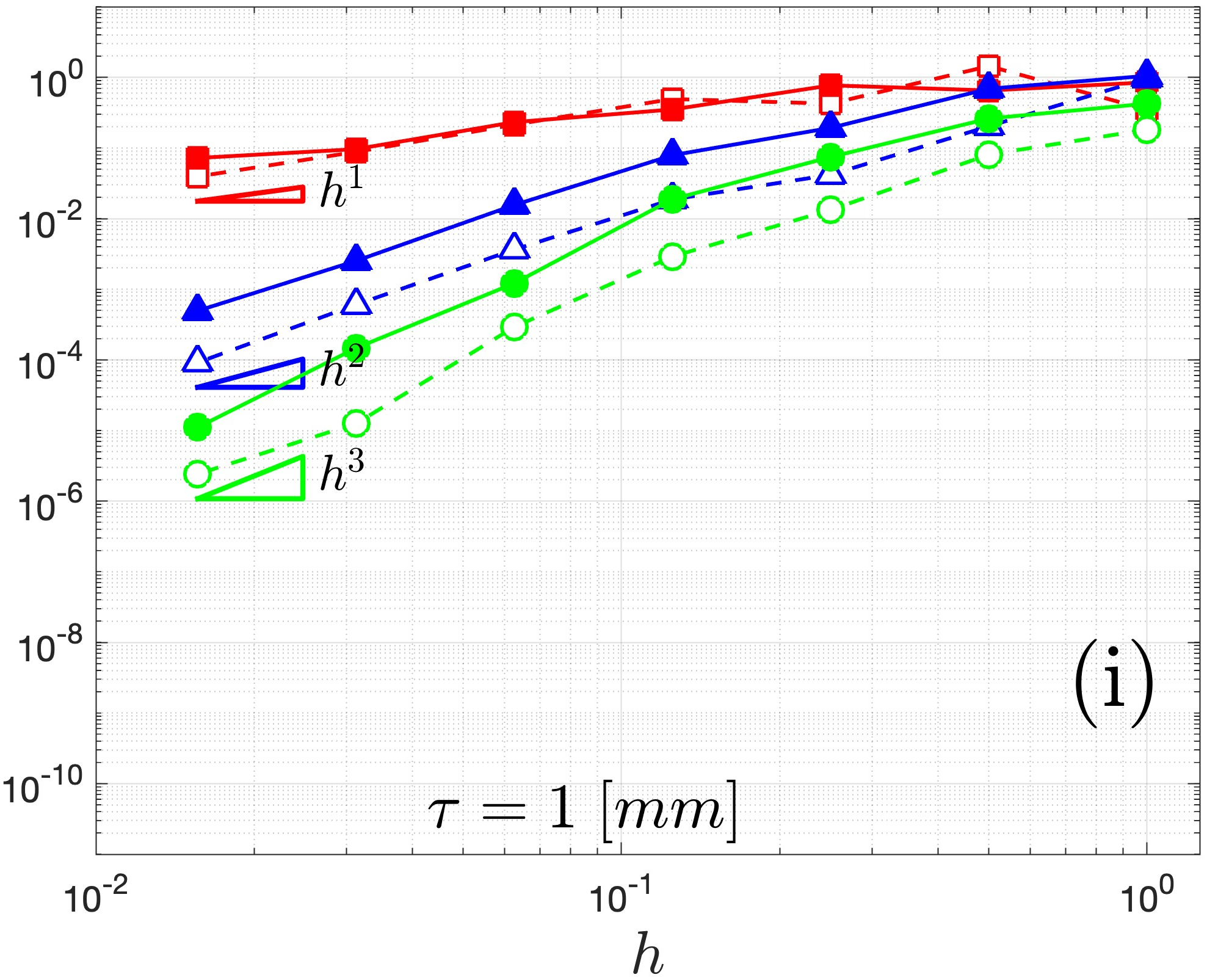}
    \end{subfigure}\\
    \begin{subfigure}[b]{0.99\textwidth}   \includegraphics[width = 0.99\textwidth]{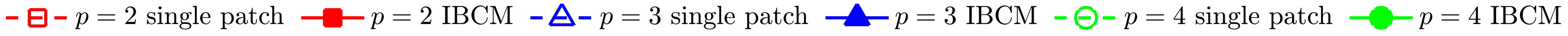}
    \end{subfigure} 
	\caption{Convergence curves for the plate shown in Fig.\figref{fig:RES - LamGeo} in $L^2$ error norm and $H^1$ and $H^2$ error seminorms for a Kirchhoff-Love theory. The curves are obtained for three different polynomial orders $p=2,3,4$ and three thickness values $\tau=100,10,1$ [mm]. Two discretization are taken into account, a single trimmed patch as shown in Fig.\figref{fig:RES - LamGeom a} and a IBCM-based one as shown in Fig.\figref{fig:RES - LamGeom b}.} \label{fig:RES - lam_KL_Conv}
\end{figure}

Fig.\figref{fig:RES - lam_KL_Conv} shows the same curves but for the Kirchhoff-Love discretization. In this case, first-order polynomials are not considered as they do not satisfy the $C^1$ requirement of the Kirchhoff-Love theory. This time, the $H^2$ seminorm error convergence is also shown since it is significant in evaluating the bending moment. Similar observations made for the previous case hold here as well. The convergence curves approach optimal rates. Additionally, Kirchhoff-Love elements do not suffer from shear locking in the early refinement stages. However, some saturation of the error appears in the final refinement level for $\tau=1$ [mm] and $p=4$ due to the ill-conditioning of the stiffness matrix typical of higher polynomial discretizations combined with weak coupling conditions. 

Regarding stability, the positive definiteness of the stiffness matrix is maintained in all tests with an IBCM-based discretization, while coercivity is lost for the single patch discretization in the last refinement level for $p=3,4$ at  $\tau=100$ [mm] and in the last two refinement levels for $p=3,4$ at $\tau=1$ [mm].

\begin{figure}	\centering
    \includegraphics[width = 0.90\textwidth]{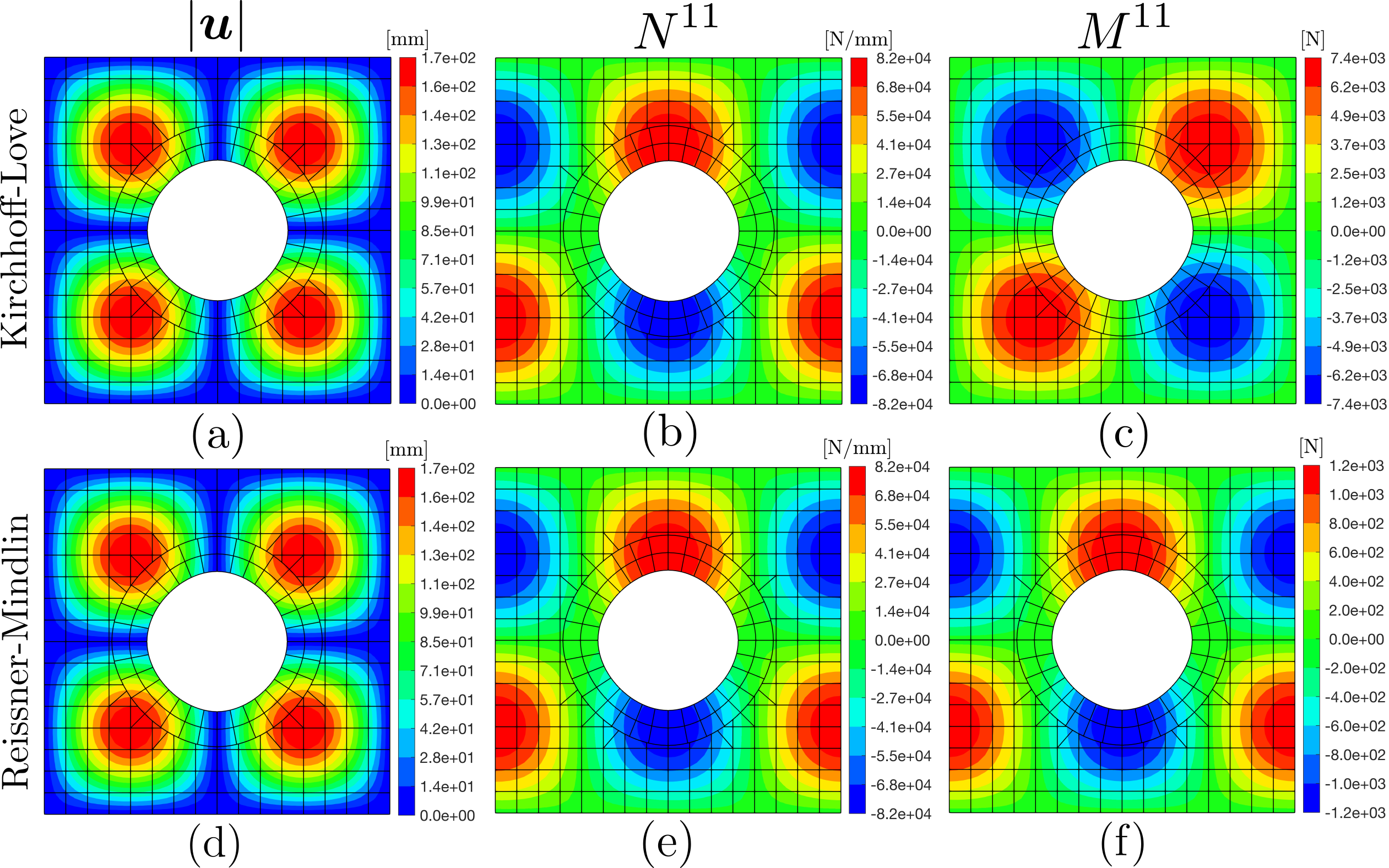}
    \caption{Contour plots of the laminated plate in Section \ref{ssec:Res - lam} with superimposed discretization. Displacement magnitude (a), representative component of the membrane force $N^{11}$ (b), representative component of the bending moment $M^{11}$ (c) for the Kirchhoff-Love theory, and correspondents contour plots for the Reissner-Mindlin theory (d), (e), (f), respectively. It is worth mentioning that the noticeable difference between (c) and (f) is due to the different manufactured solution adopted in this test. In particular assigning also a manufactured rotation field in the Reissner-Mindlin discretization affects in this case only the bending response.} \label{fig:RES - Lam_Cont}
\end{figure}

To conclude, Fig.\figref{fig:RES - Lam_Cont} shows the contour plots of the displacement and a component of the membrane force and bending moment for both the Kirchhoff-Love and Reissner-Mindlin discretizations, with superimposed mesh, demonstrating the smooth coupling between the internal patch and the boundary layer.

\subsection{Kirchhoff-Love hyperbolic paraboloid with cut-outs}\label{ssec:Res - hyp}
\begin{figure}	\centering
    \includegraphics[width = 0.60\textwidth]{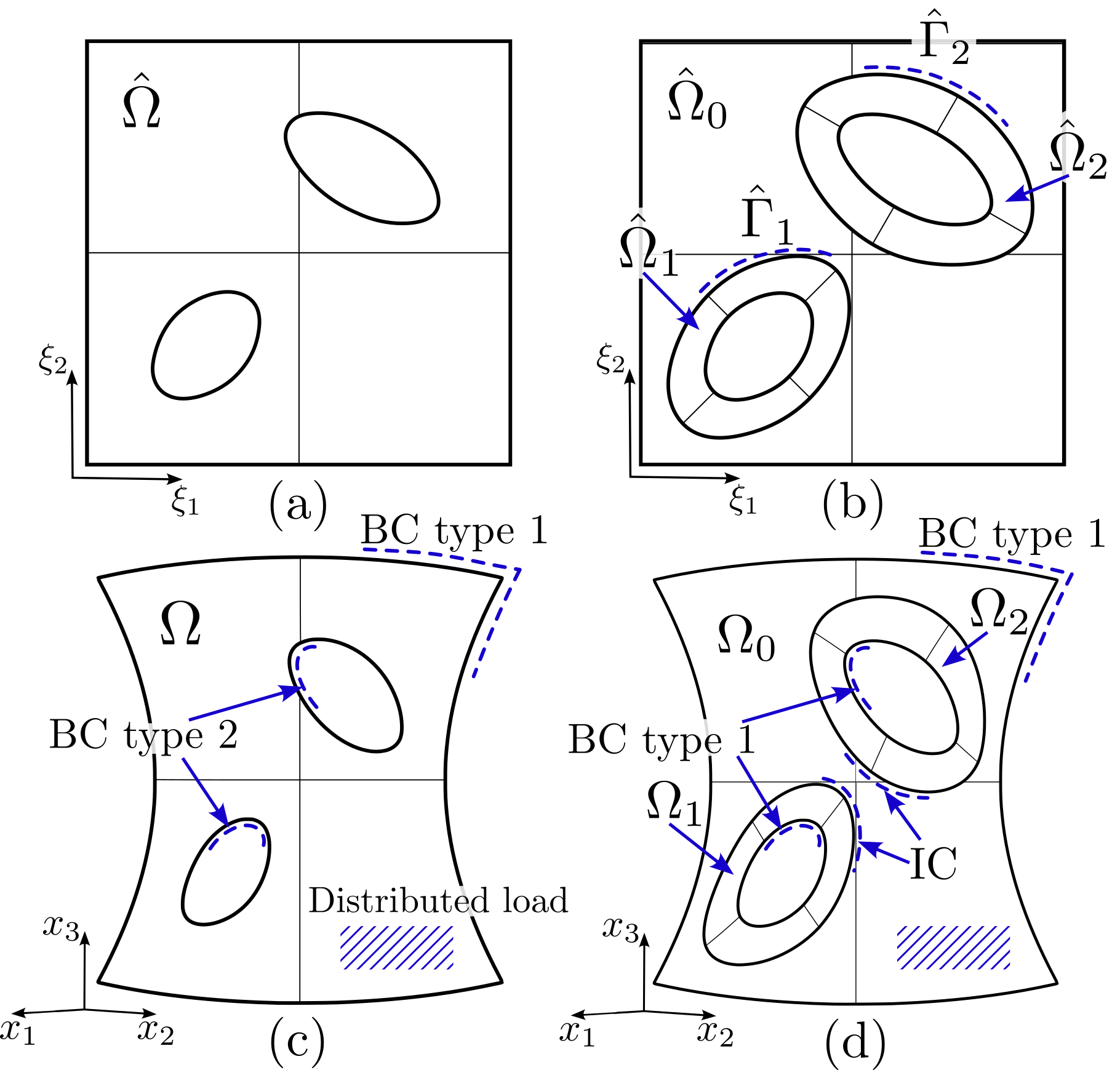}
    \begin{subfigure}{0\textwidth}\phantomcaption\label{fig:RES - HypGeom a}\end{subfigure}
    \begin{subfigure}{0\textwidth}\phantomcaption\label{fig:RES - HypGeom b}\end{subfigure}
    \begin{subfigure}{0\textwidth}\phantomcaption\label{fig:RES - HypGeom c}\end{subfigure}
    \begin{subfigure}{0\textwidth}\phantomcaption\label{fig:RES - HypGeom d}\end{subfigure}
    \caption{Geometry and discretization for the hyperbolic paraboloid shell in the second set of tests. First refinment level for the single patch discretization in parametric (a) and physical (c) domains, with corresponding types of boundary conditions and applied loads. First refinment level for the IBCM-based discretization in the parametric (b) and physical (d) domains.} \label{fig:RES - HypGeo}
\end{figure}

This next test focuses on a hyperbolic paraboloid geometry taken from the new shell obstacle course proposed in \cite{benzaken2021}. For brevity, the complete description of the geometrical map is not provided here but can be found in the aforementioned work. Additionally, two elliptical holes are constructed in the parametric domain to test the IBCM-based discretization. Figs.\figref{fig:RES - HypGeom a} and \figref{fig:RES - HypGeom c} show the resulting geometry in the parametric and physical domains, respectively. The shell section is a single-layer isotropic material with a Young's modulus of  $E=10$ [GPa] and $\nu=0.3$. The following displacement field is adopted to construct a manufactured solution:
\begin{equation}
    \bm{u}_{ex} = \mathrm{U}_1\xi_2\sin{\left(\frac{\pi}{2}\xi_2\right)}\bm{e}_1 + \mathrm{U}_2\xi_2\sin{\left(\frac{\pi}{2}\xi_2\right)}\bm{e}_2,
\end{equation}
where $\mathrm{U}_1=\mathrm{U}_2=0.1$ [m]. This displacement field results in non-homogeneous Dirichlet boundary conditions for both displacement and rotation. The shell theory adopted for this test is uniquely the Kirchhoff-Love one, as proposed in the benchmark. As in the previous test, two discretization strategies are employed: a trimmed single-patch and an IBCM-based discretization. The first refinement levels are shown in Figs.\figref{fig:RES - HypGeom c} and \figref{fig:RES - HypGeom d}.

\begin{figure}	\centering
    \begin{subfigure}[b]{0.32\textwidth}  \includegraphics[width = 0.99\textwidth]{RES_ConvergenceTitle_L2.png}
    \end{subfigure}
    \begin{subfigure}[b]{0.32\textwidth}  \includegraphics[width = 0.99\textwidth]{RES_ConvergenceTitle_H1.png}
    \end{subfigure}
    \begin{subfigure}[b]{0.32\textwidth}  \includegraphics[width = 0.99\textwidth]{RES_ConvergenceTitle_H2.png}
    \end{subfigure}\\  
    \begin{subfigure}[b]{0.32\textwidth}  \includegraphics[width = 0.99\textwidth]{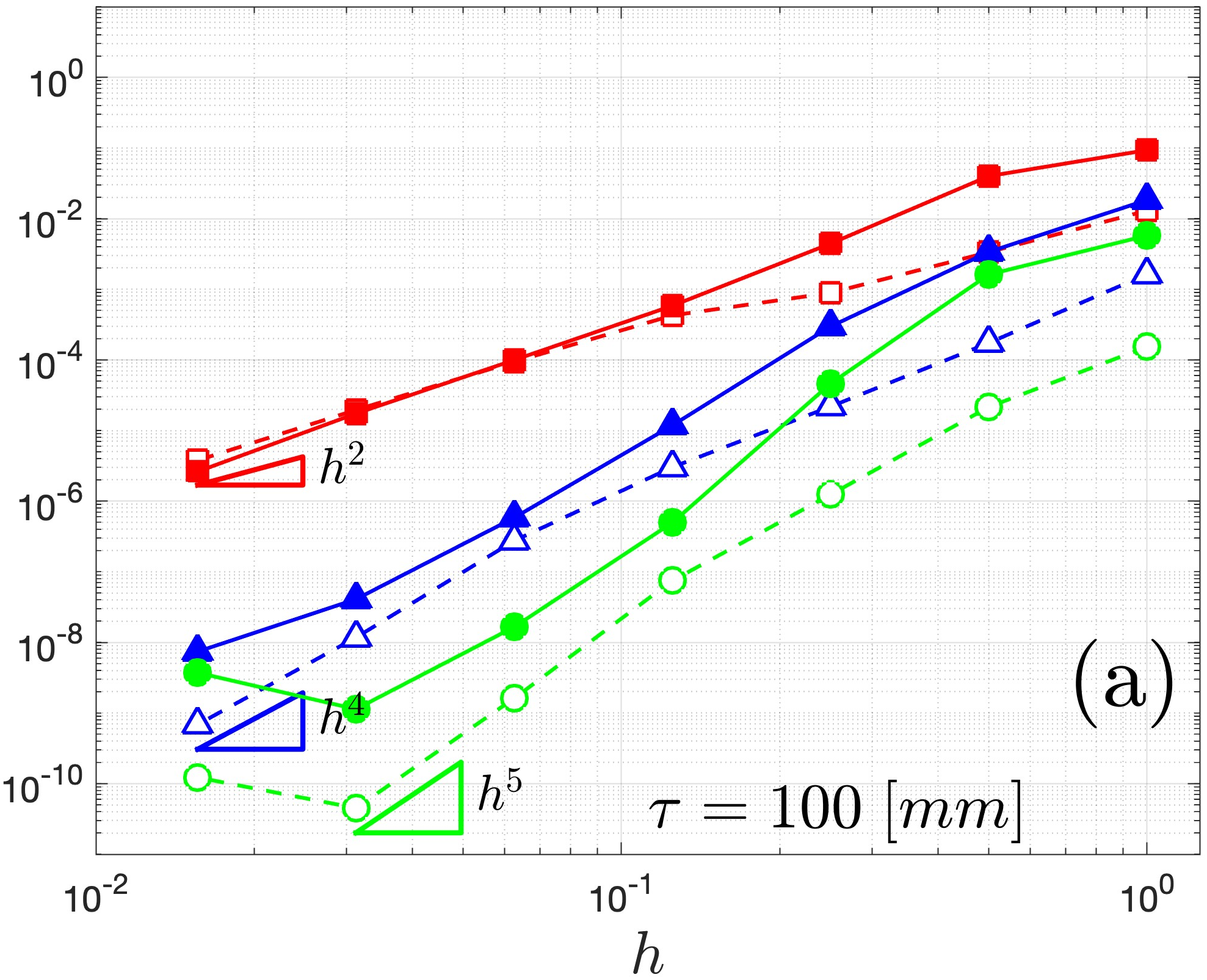}
    \end{subfigure}
    \begin{subfigure}[b]{0.32\textwidth}  \includegraphics[width = 0.99\textwidth]{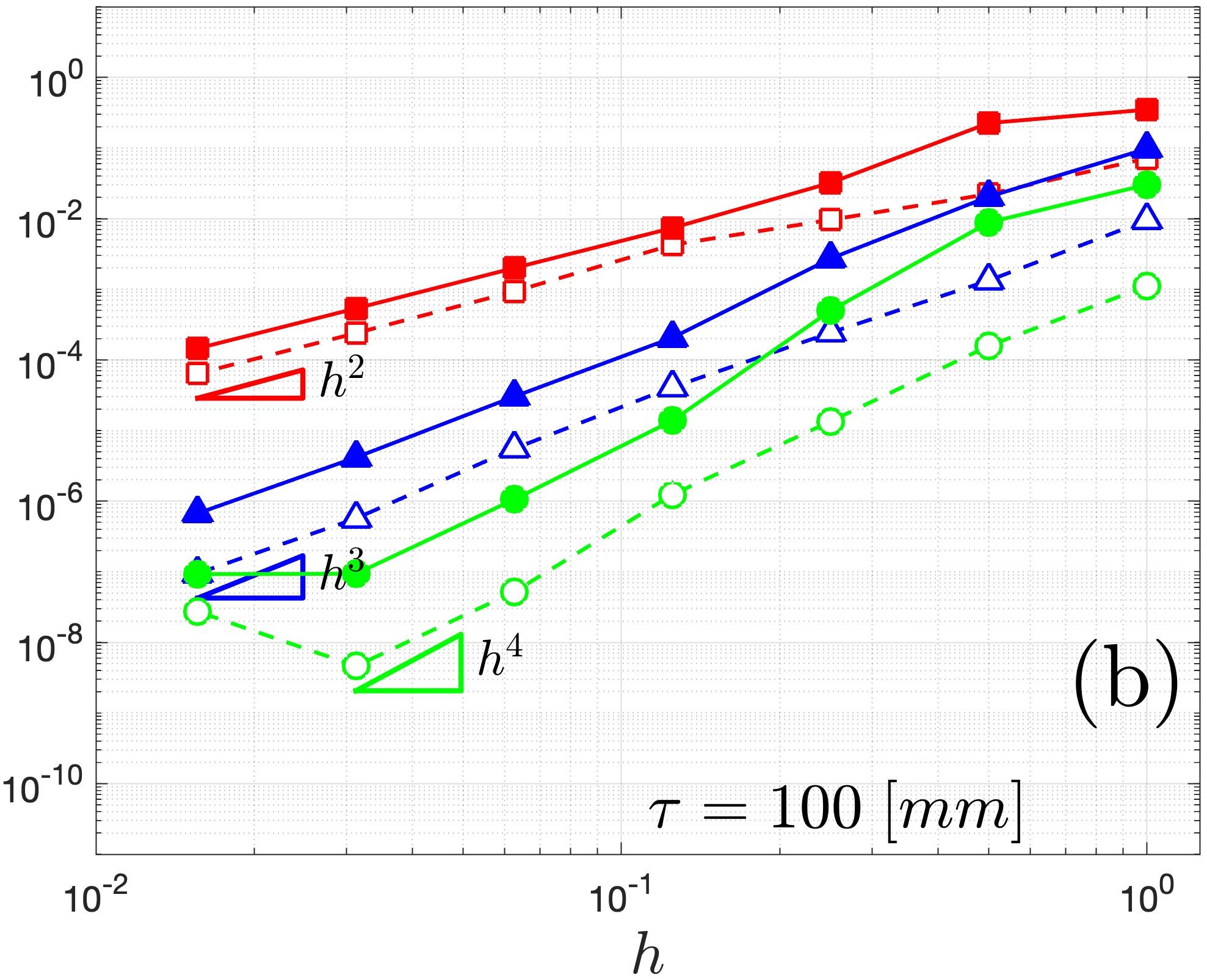}
    \end{subfigure}
    \begin{subfigure}[b]{0.32\textwidth}  \includegraphics[width = 0.99\textwidth]{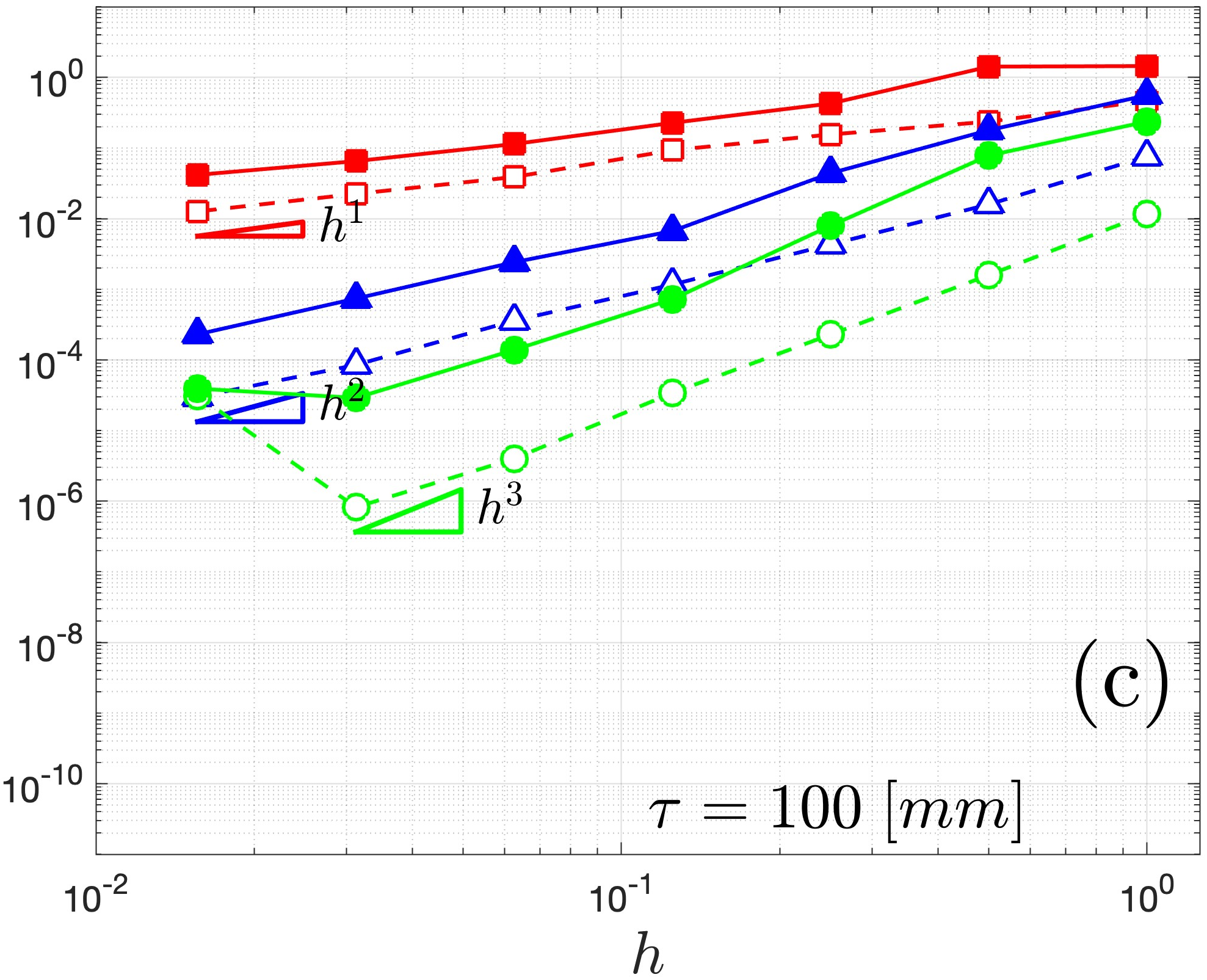}
    \end{subfigure}\\
    \begin{subfigure}[b]{0.32\textwidth}  \includegraphics[width = 0.99\textwidth]{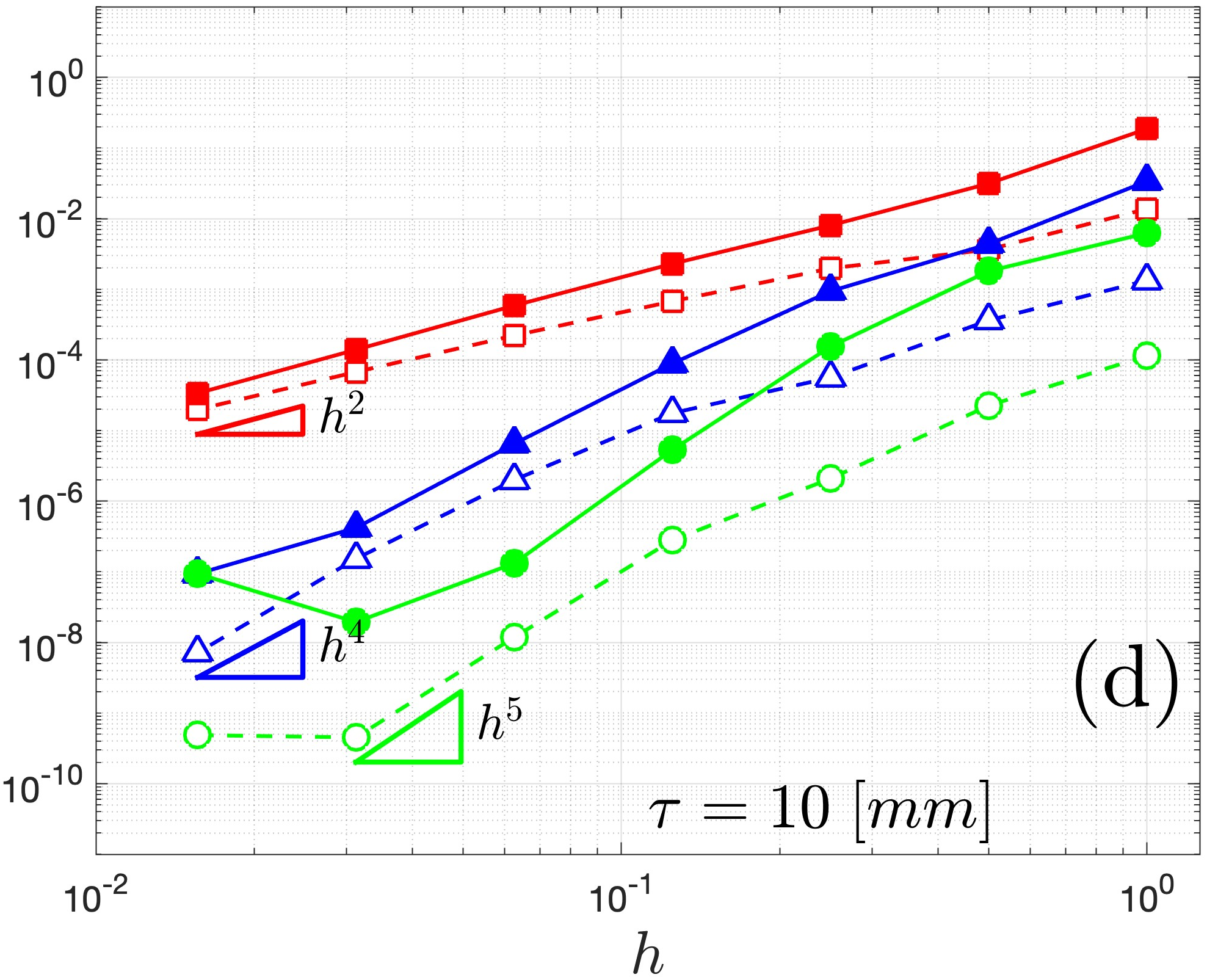}
    \end{subfigure}
    \begin{subfigure}[b]{0.32\textwidth}  \includegraphics[width = 0.99\textwidth]{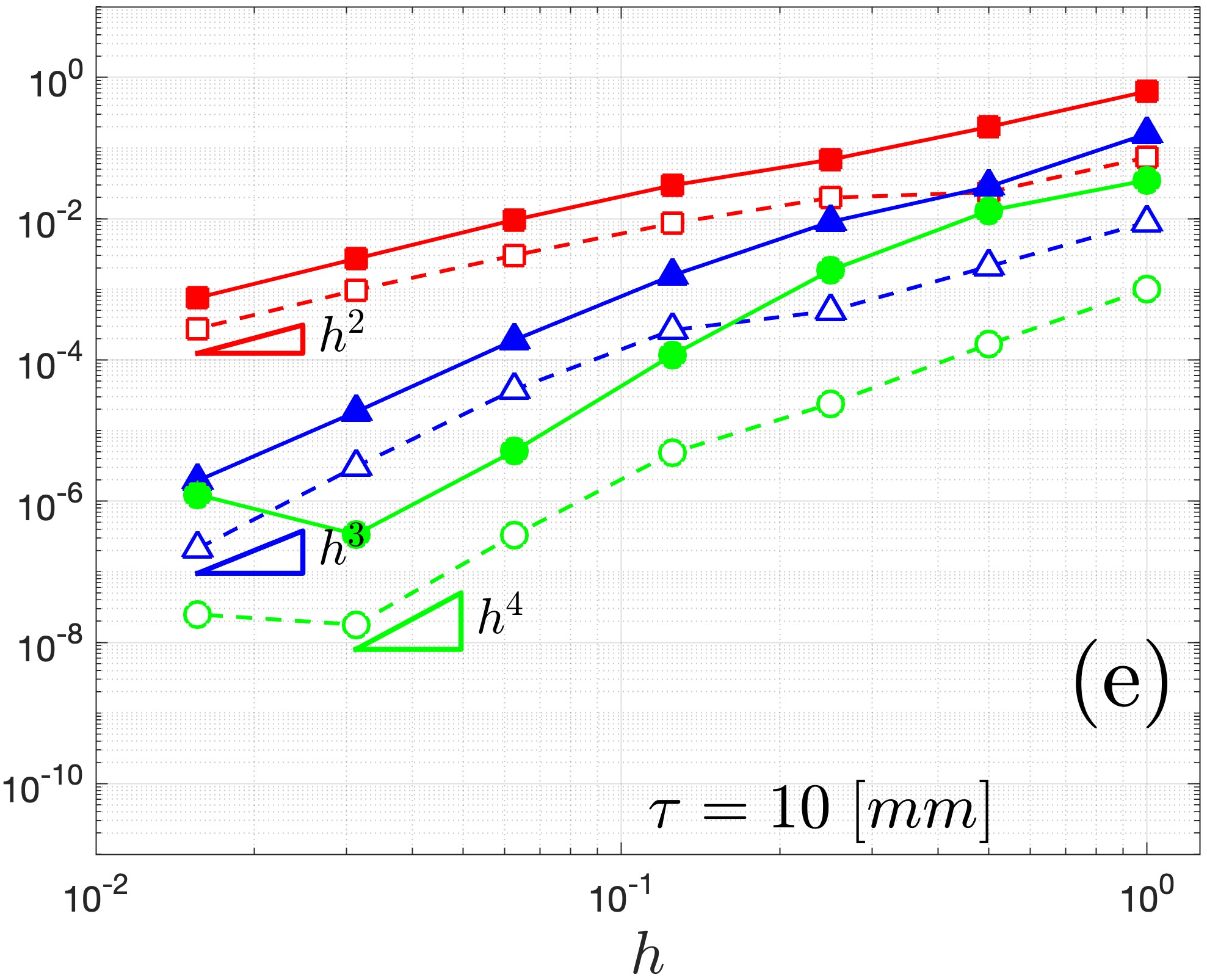}
    \end{subfigure}
    \begin{subfigure}[b]{0.32\textwidth}  \includegraphics[width = 0.99\textwidth]{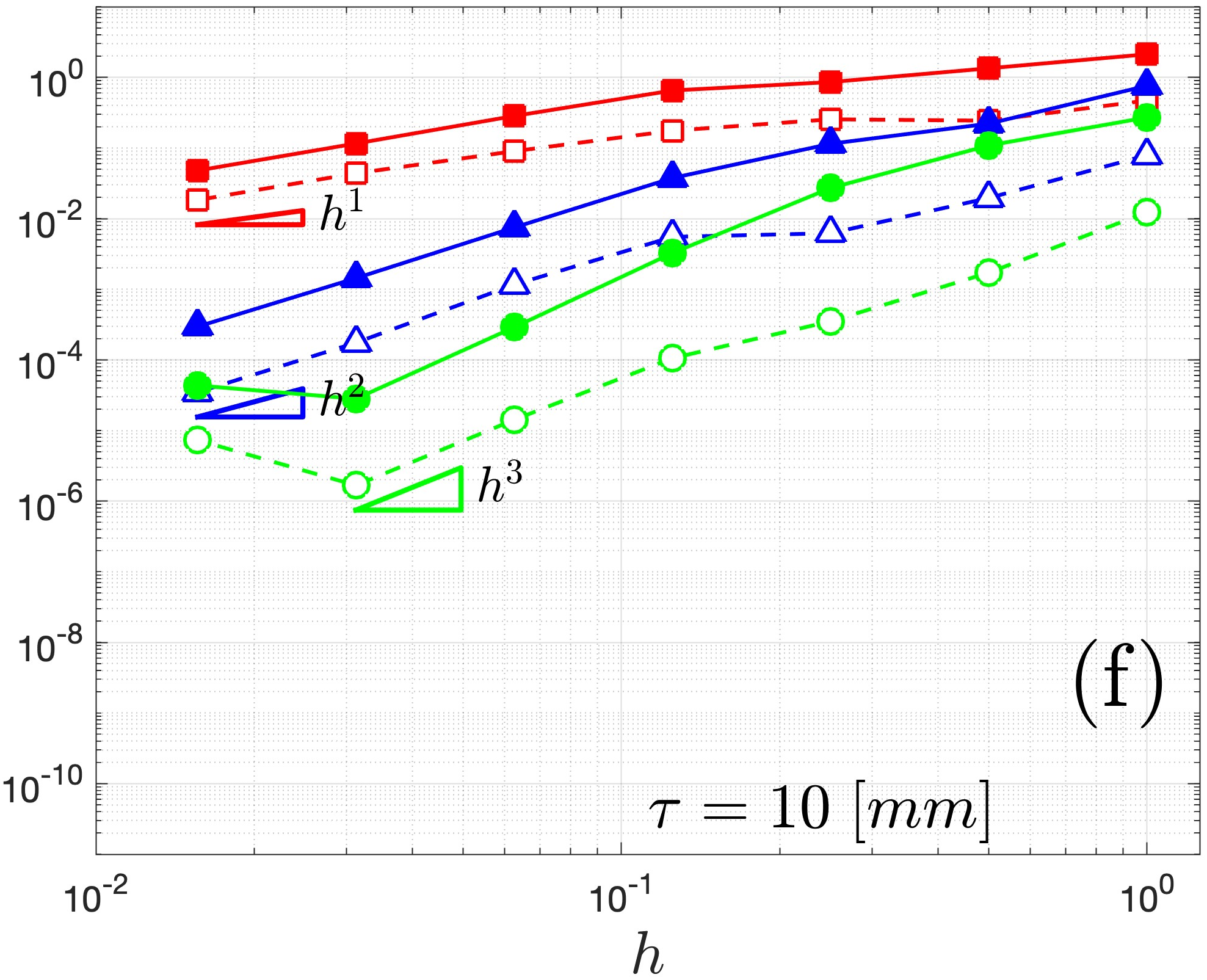}
    \end{subfigure}\\
    \begin{subfigure}[b]{0.32\textwidth}  \includegraphics[width = 0.99\textwidth]{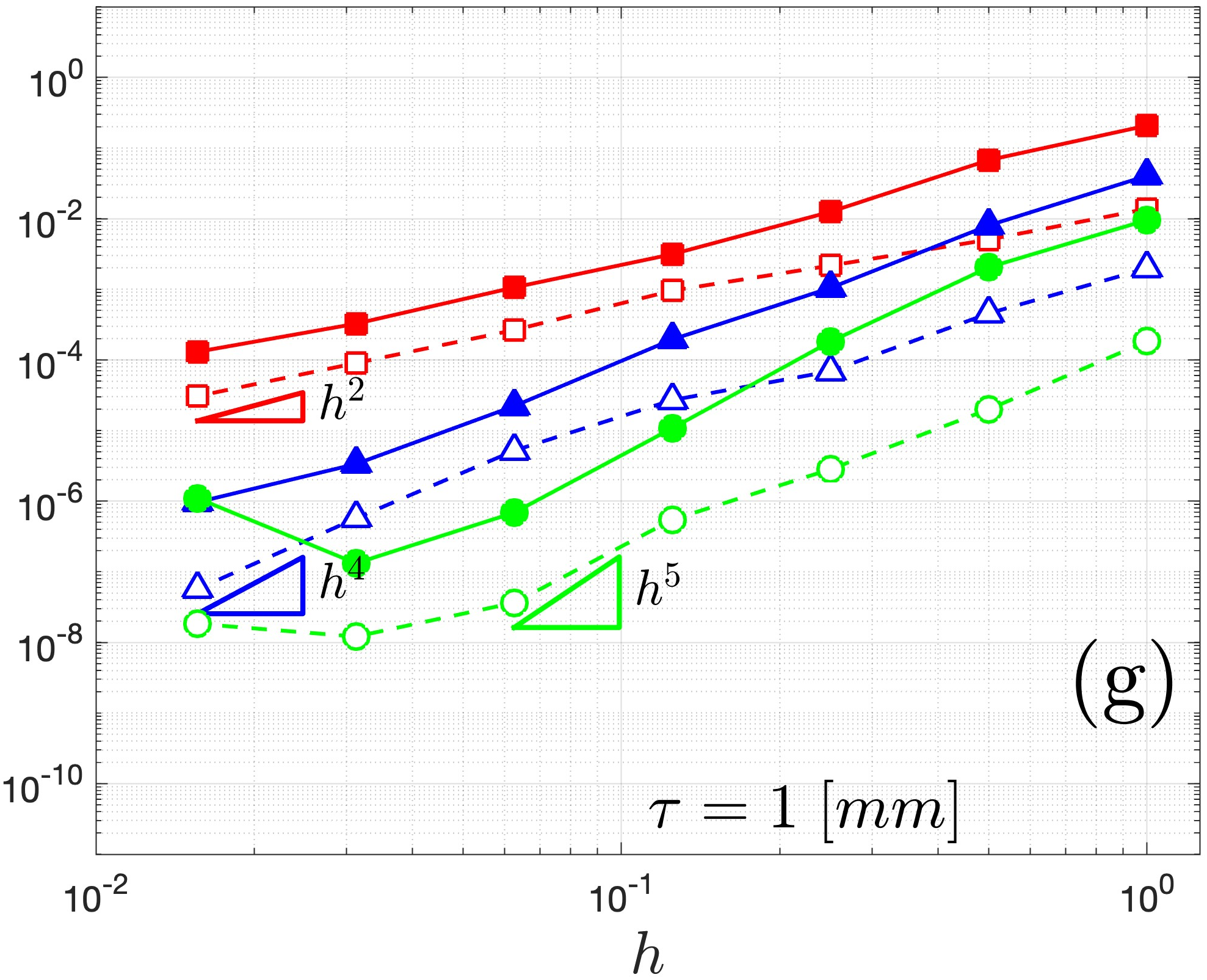}
    \end{subfigure}
    \begin{subfigure}[b]{0.32\textwidth}  \includegraphics[width = 0.99\textwidth]{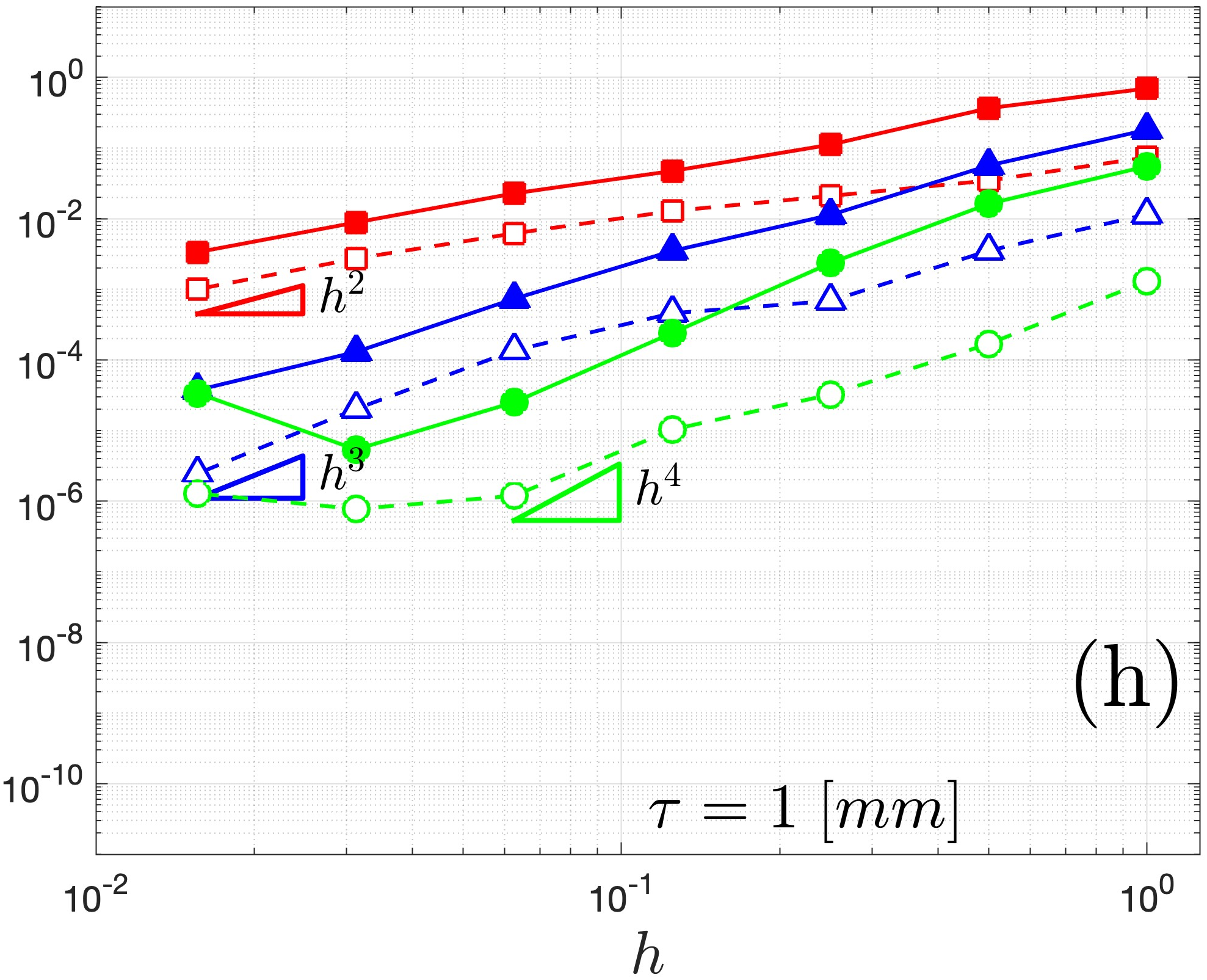}
    \end{subfigure}
    \begin{subfigure}[b]{0.32\textwidth}  \includegraphics[width = 0.99\textwidth]{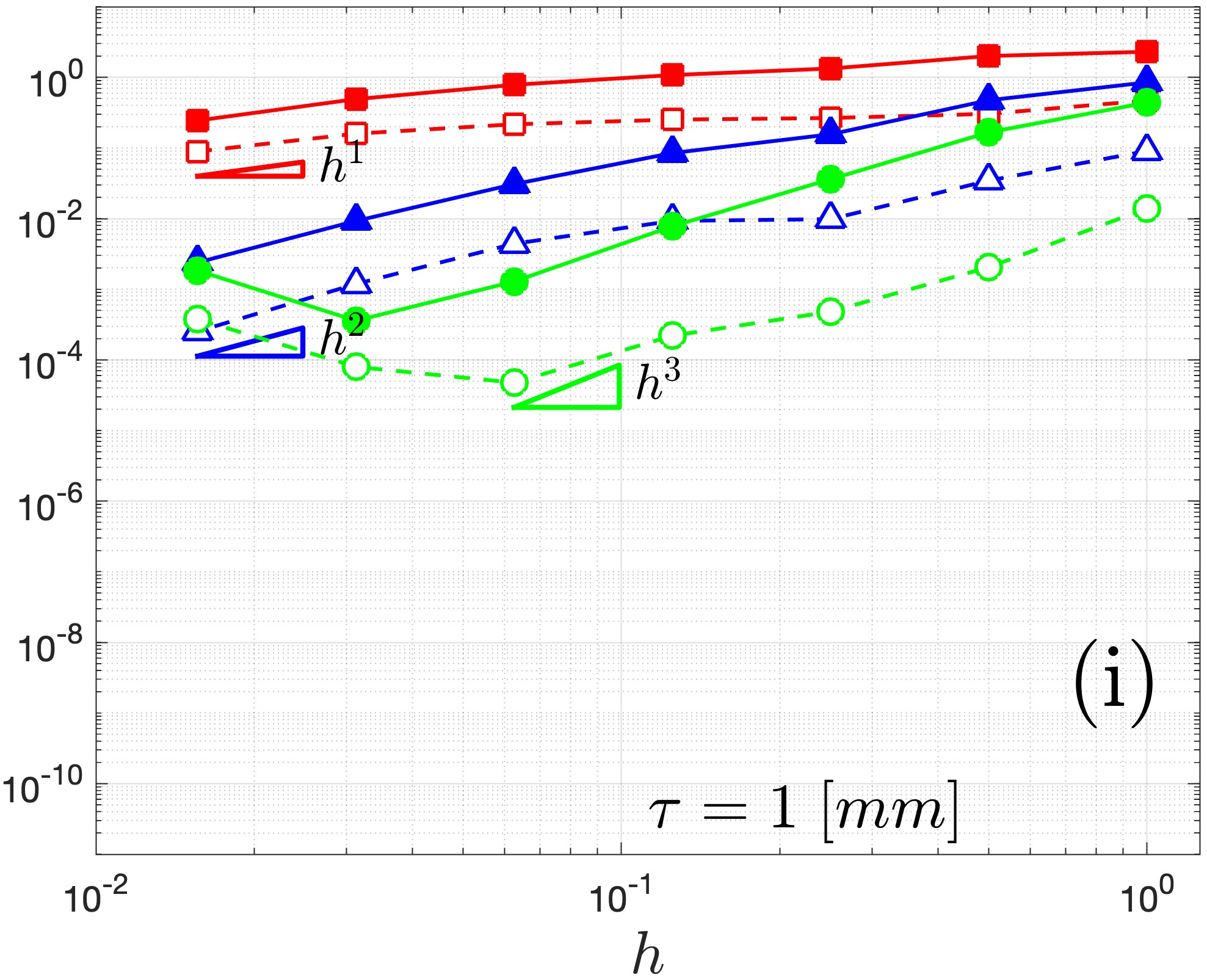}
    \end{subfigure}\\
    \begin{subfigure}[b]{0.99\textwidth}   \includegraphics[width = 0.99\textwidth]{RES_ConvergenceLegend_KL.png}
    \end{subfigure} 
    \caption{Convergence curves for the shell shown in Fig.\figref{fig:RES - HypGeo} in $L^2$ error norm and $H^1$ and $H^2$ error seminorms for a Kirchhoff-Love theory. The curves are obtained for three different polynomial orders $p=2,3,4$ and three thickness values $\tau=100,10,1$ [mm]. Two discretization are taken into account, a single trimmed patch as shown in Fig.\figref{fig:RES - HypGeom a} and a IBCM-based one as shown in Fig.\figref{fig:RES - HypGeom b}.} \label{fig:RES - hyp_KL_Conv}
\end{figure}

The convergence of the error in $L^2$ norm and $H^1$ and $H^2$ seminorms are shown in the graphs collected in Fig.\figref{fig:RES - hyp_KL_Conv} for different thicknesses, specifically $\tau=100$ [mm], $\tau=10$ [mm], and $\tau=1$ [mm], and for the polynomial values  $p=2,3,4$. The parameter $\beta=100$, and the choice for the mesh size to scale the penalty is performed as in the previous test. It can be noticed that, in general, the error convergence curves corresponding to the IBCM discretization follow the corresponding ones related to the single patch discretization with a slightly higher error, which is due to the same reason as in the previous case. In the case of  $p=4$, this difference reaches one order of magnitude. However, the optimal convergence rates are approached except for the last refinement level for  $p=4$ where the stiffness matrix starts suffering from ill-conditioning. This effect becomes even more noticeable in the $H^1$ and $H^2$ norms.

The positive definiteness of the stiffness matrix is achieved in all tests with an IBCM-based discretization, while this is not the case for the single patch discretization for $\tau=100$ [mm] in the last two refinement levels for $p=3,4$ and in the last refinement level for $p=2$, and for $\tau=10$ [mm] in the last refinement level for $p=2,3,4$.

\begin{figure}	\centering
    \includegraphics[width = 0.90\textwidth]{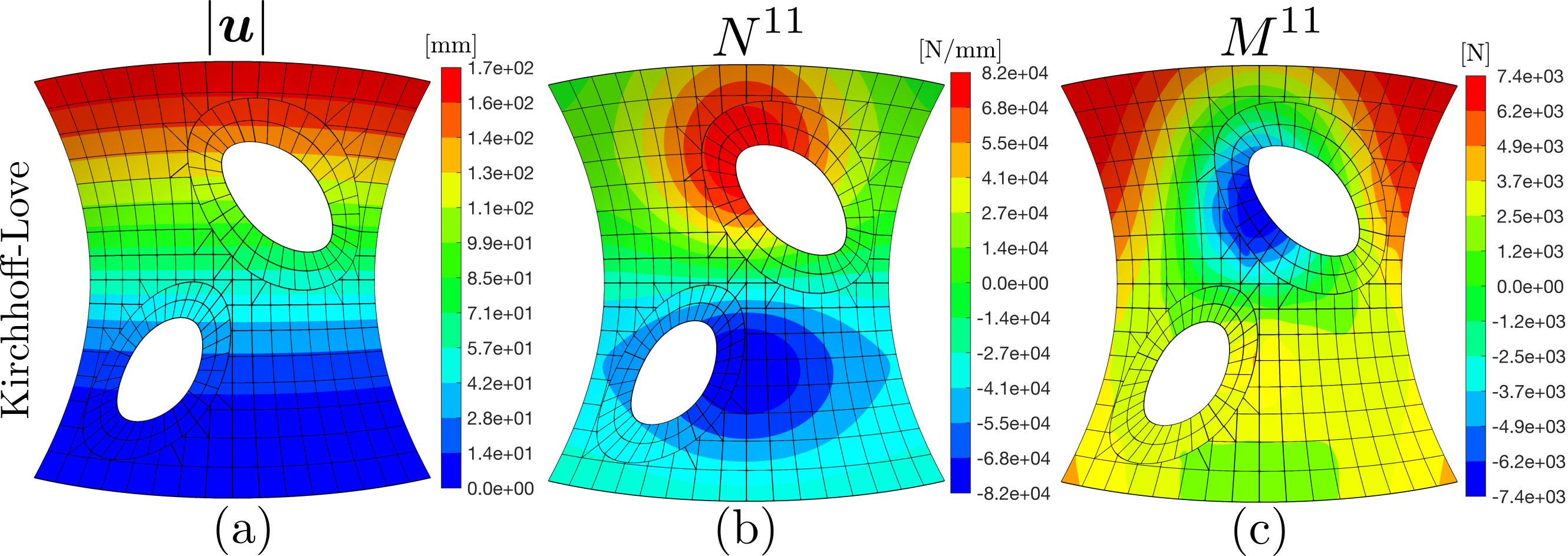}
    \caption{Contourplots of the shell in Section \ref{ssec:Res - hyp} with superimposed discretization. Displacement magnitude (a), representative component of the membrane force $N^{11}$, representative component of the bending moment $M^{11}$ for the Kirchhoff-Love theory.}
\label{fig:RES - Hyp_Cont}
\end{figure}

For completeness, some contour plots representative of the displacement, membrane force, and bending moment fields are shown in Fig.\figref{fig:RES - Hyp_Cont}, with superimposed mesh. The smoothness can be observed in the first two graphs, whereas in the last one some artifacts can be noticed at the interface between the inner patch and the boundary layers. These artifacts tend to disappear with further mesh refinements.

\subsection{Conformal coupling of trimmed cylindrical shells}\label{ssec:Res - cyl}
\begin{figure}	\centering
    \includegraphics[width = 0.80\textwidth]{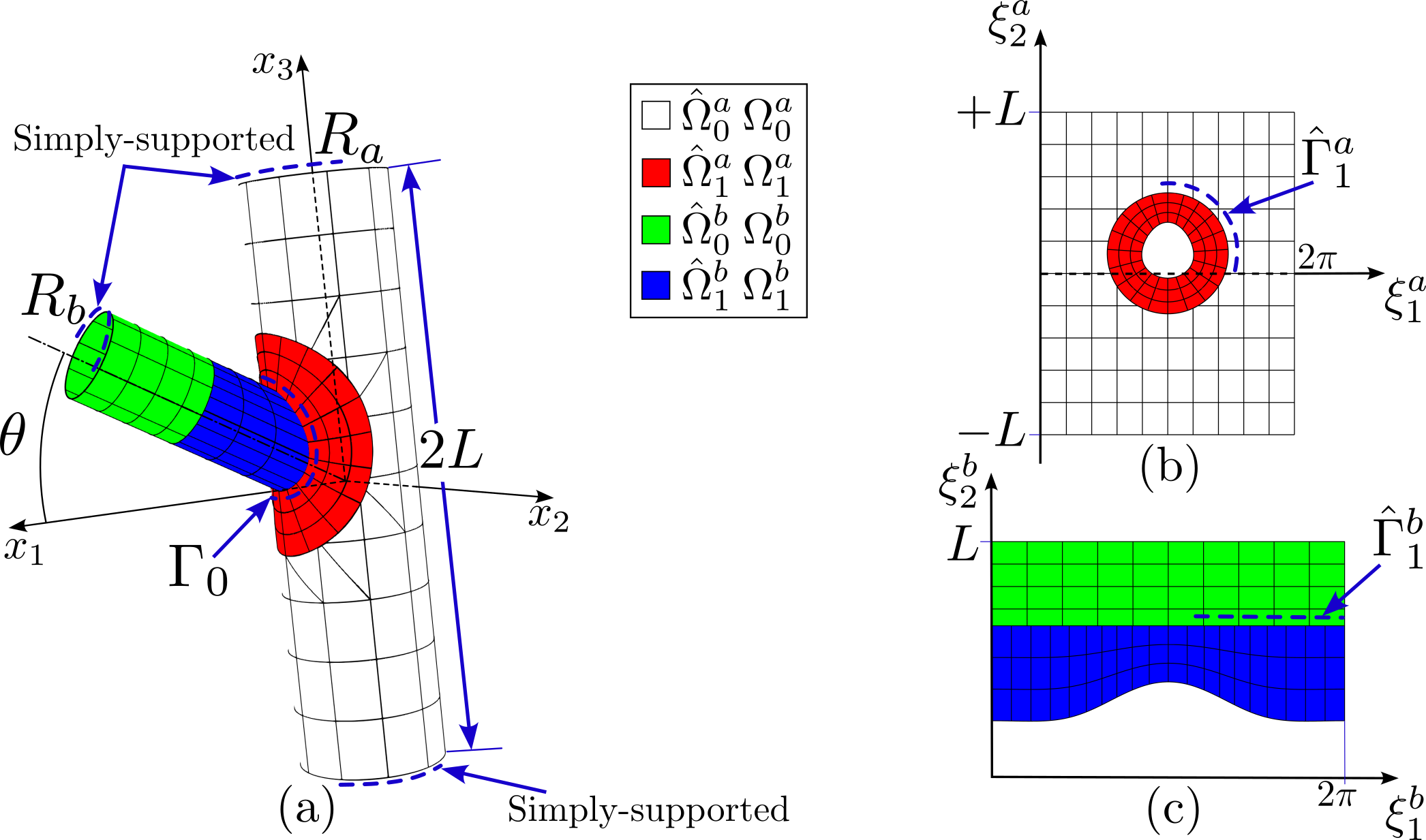}
    \begin{subfigure}{0\textwidth}\phantomcaption\label{fig:RES - CylGeom a}\end{subfigure}
    \begin{subfigure}{0\textwidth}\phantomcaption\label{fig:RES - CylGeom b}\end{subfigure}
    \begin{subfigure}{0\textwidth}\phantomcaption\label{fig:RES - CylGeom c}\end{subfigure}
    \caption{Geometry and discretization of the structure in Section \ref{ssec:Res - cyl} in the Euclidean space (a). Parametric domains of each cylinder, (a) and (b), showing the inner patches and the boundary layers.} \label{fig:RES - CylGeo}
\end{figure}

In this test, two boundary layers are constructed at the intersection between two cylindrical patches to allow for a conformal discretization at the interface and a strong coupling of the displacements, as shown in Fig.\figref{fig:RES - CylGeo}. The geometries of the cylindrical shells are defined through the analytic maps
\begin{equation}
    \bm{x}^a (\xi_1^a,\xi_2^a) = \hat{\bm{\mathcal{F}}}^a(\xi_1^a,\xi_2^a) = \begin{bmatrix} -R_a \cos{(\xi_1^a)}\\ -R_a \sin{(\xi_1^a)}\\ \xi_2 \end{bmatrix} \;,
\end{equation}
where $(\xi_1^a,\xi_2^a)\in\hat{\Omega}^a=[0,2\pi]\times[-L,L]$, and 
\begin{equation}
    \bm{x}^b (\xi_1^a,\xi_2^a) = \hat{\bm{\mathcal{F}}}^b(\xi_1^b,\xi_2^b) = \begin{bmatrix}\cos{(\theta)}&0&-\sin{(\theta)}\\0&1&0\\ \sin{(\theta)}&0&\cos{(\theta)}\end{bmatrix} \begin{bmatrix}-R_b\cos{(\xi_1^b)}\\-R_b\sin{(\xi_1^b)}\\\xi_2^b\end{bmatrix} \;,
\end{equation}
where $(\xi_1^b,\xi_2^b)\in\hat{\Omega}^b=[0,2\pi]\times[0,L]$. The parameters for the geometry are set as follows: $L=4$ [m], $R_a=1$ [m], $R_b=0.6$ [m], and $\theta=-\pi/3$ [rad].

Fig.\figref{fig:RES - CylGeo} displays the geometry of the structure, highlighting the intersection between the two cylinders, which creates an interface $\Gamma_0$. The shell section has a thickness $\tau=1$ [cm] and is composed of a single isotropic layer with a Young's modulus of $E=100$ [GPa] and a Poisson's ratio of $\nu=0.3$. A uniform distributed force of $\bar{\bm{f}}=\{1,1,1\}^\Tr$ [MPa] is applied on the surfaces of the shells, while the external boundaries are simply-supported, meaning $\bm{u}=0$ at $\xi_2^a=\pm L/2$ and $\xi_2^b=L$. Using the procedure described in Section \ref{sec:IBCM}, a boundary layer is constructed for each patch in such a way that they result in a conformal interface. Furthermore, if the two boundary layers have the same discretization in the interface direction and if the two approximation spaces are constructed with the same properties, the coupling of the displacement degrees of freedom can be realized in a strong sense. It is worth noting, however, that the rotation still requires a weak coupling. In the Kirchhoff-Love theory, this necessity arises from the absence of rotation as a main variable, whereas in the Reissner-Mindlin theory, this is motivated by the different local contravariant bases for the rotation variables between the two patches.

Fig.\figref{fig:RES - CylGeo} illustrates the discretization of the structure in the Euclidean space as well as the parametric domains of the cylindrical surfaces, each constituted by a inner patch and a boundary layer. It can be noticed how, in the Euclidean space, the discretization from both sides is conformal at the interface, allowing the strong coupling.

\begin{figure}	\centering
    \includegraphics[width = 0.95\textwidth]{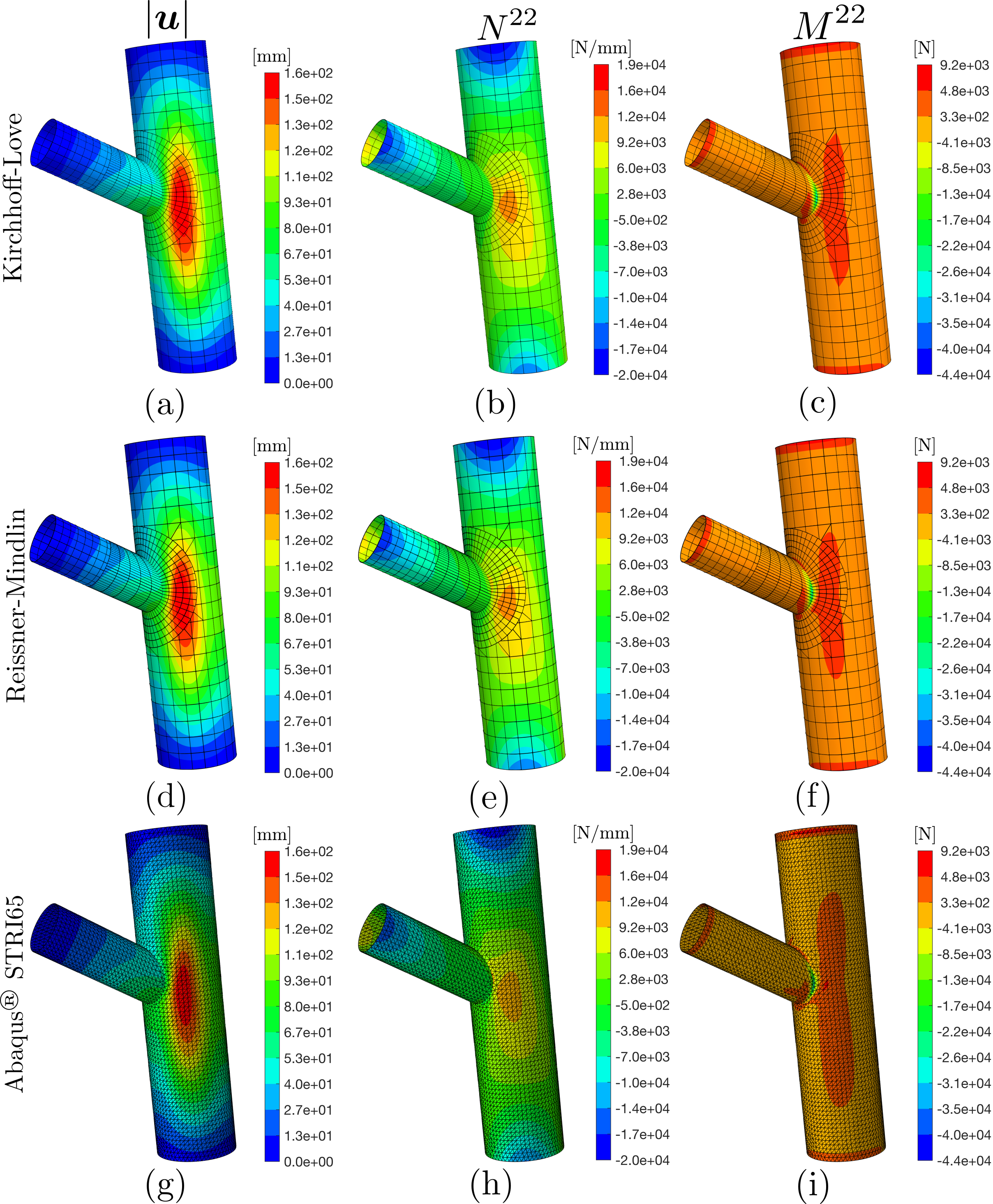}
    \caption{Contour plots for the intersecting cylindrical shells in Section \ref{ssec:Res - cyl} with superimposed discretization. Displacement magnitude (a), representative component of the membrane force $N^{22}$ (b), representative component of the bending moment $M^{22}$ (c) for the Kirchhoff-Love theory, and corresponding contour plots for the Reissner-Mindlin theory (d), (e), (f), respectively. As a reference, the same contour plots obtained with Abaqus\textsuperscript{\textregistered} quadratic triangular elements STRI65 are shown in (g), (h), (i), with superimposed mesh.} \label{fig:RES - CylCont}
\end{figure}

The contour plots of the solution are shown in Fig.\figref{fig:RES - CylCont}. In particular, the magnitude of the displacement, a representative component of the membrane force, and a representative component of the bending moment are shown. Due to the small thickness of the shells, the corresponding plots for the Kirchhoff-Love theory and the Reissner-Mindlin one are very similar. For reference, in (g), (h), (i), the same contour plots are replicated using Abaqus\textsuperscript{\textregistered} quadratic triangular elements STRI65 \cite{abaqus}. It is observed that a denser mesh is necessary to achieve a comparable level of refinement, particularly near the interface, where the mesh becomes unstructured.

\subsection{Boundary layers over a B-spline surface with KL-RM coupling} \label{ssec:Res - mix}
\begin{figure}	\centering
    \includegraphics[width = 0.80\textwidth]{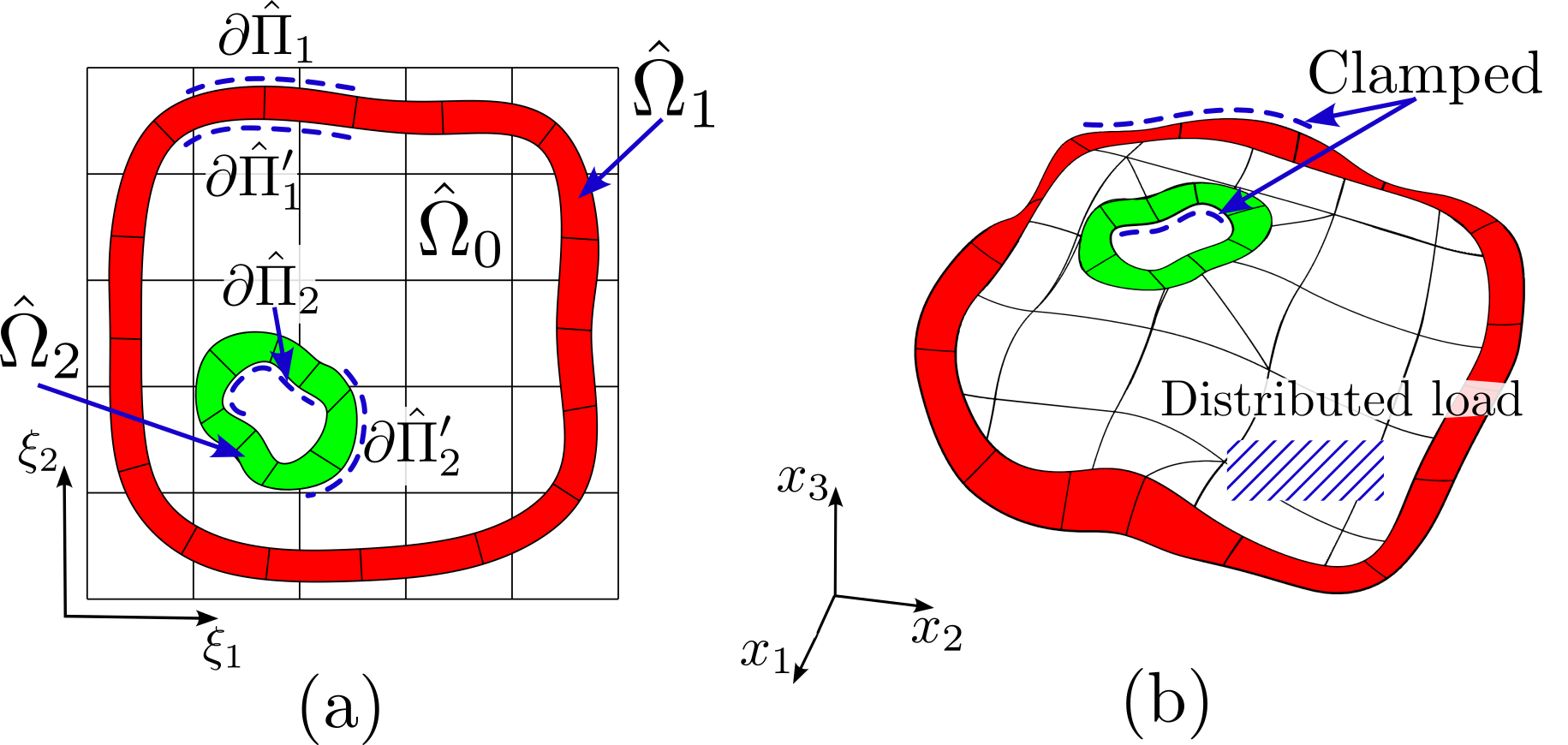}
    \begin{subfigure}{0\textwidth}\phantomcaption\label{fig:RES - MixGeom a}\end{subfigure}
    \begin{subfigure}{0\textwidth}\phantomcaption\label{fig:RES - MixGeom b}\end{subfigure}
    \caption{Geometry of the shell described in Section \ref{ssec:Res - mix} in the parametric (a) and Euclidean (b) domains.} \label{fig:RES - MixGeo}
\end{figure}

In the next test, we examine a shell whose geometry depicted in Fig.\figref{fig:RES - MixGeom b}, which features an immersed external boundary and an immersed cut-out These geometrical features are defined by the curves $\pd\hat{\Pi}_1$ and $\pd\hat{\Pi}_2$ in the parametric domain, as illustrated in Fig.\figref{fig:RES - MixGeom a}. Conformal boundary layers are created corresponding to both curves.

The material used is a laminate with four identical layers having $E_1=250$ [GPa], $E_2=10$ [GPa], $\nu_{12}=0.25$, lamination sequence [0,90,90,0], and a total shell thickness $\tau=10$ [mm]. The shell boundaries are clamped, and a distributed load $\bar{\bm{f}}=f_0 \bm{e}_3$ is applied on the shell surface, where $f_0=10^5$ [kN]. The analysis considers three different shell theory settings based on the kinematics adopted for the inner patch and the boundary layers. In the first case, both are modeled as Kirchhoff-Love shells; in the second case, both use Reissner-Mindlin theory; and in the last case, the boundary layers employ Reissner-Mindlin theory while the internal patch adopts Kirchhoff-Love theory.

\begin{figure}	\centering
    \includegraphics[width = 0.95\textwidth]{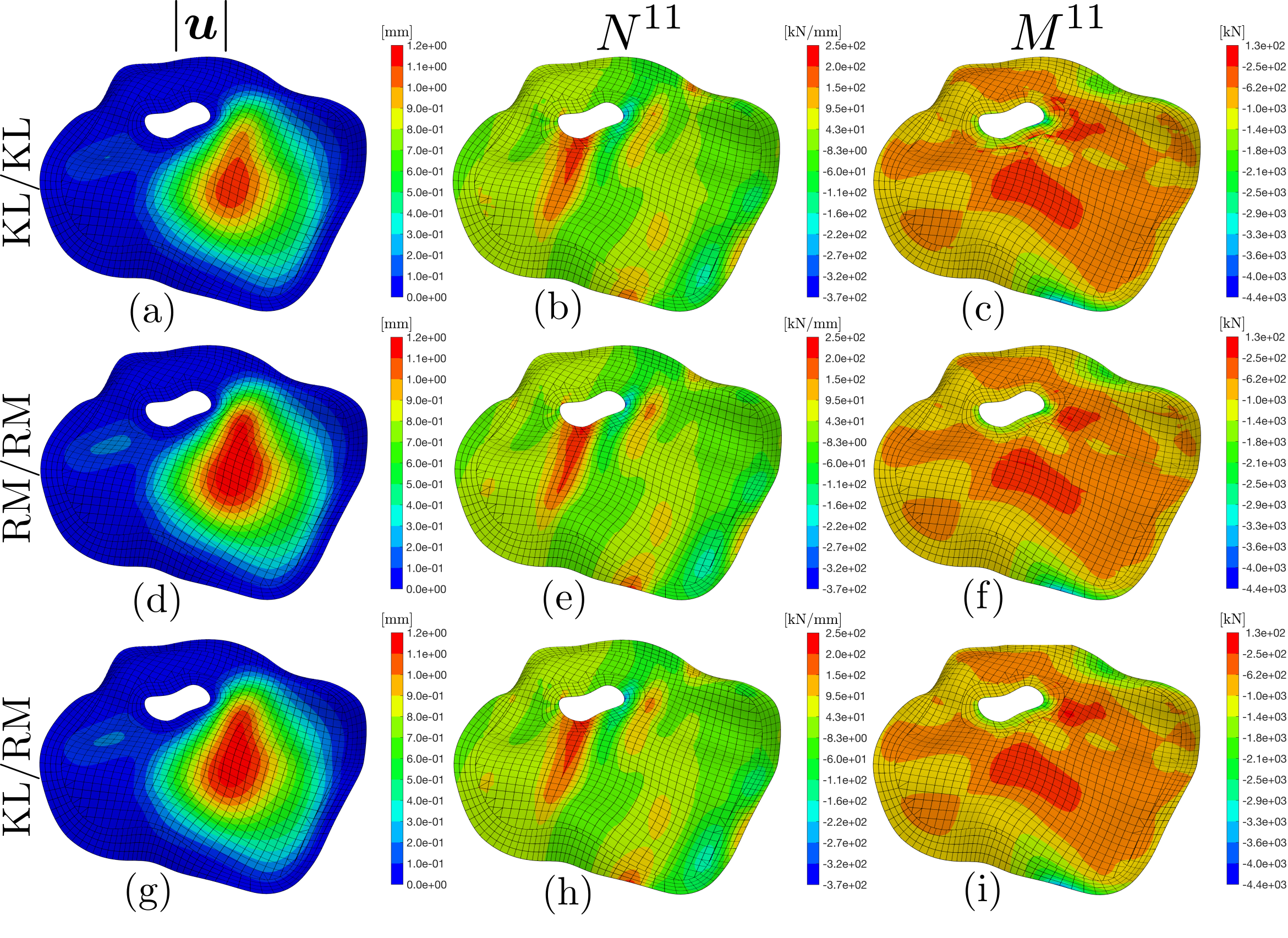}
    \caption{Contour plots of the magnitude of the displacement $|\bm{u}|$, first component of the membrane force $N^{11}$, and the first component of the bending moment $M^{11}$ for the shell described in Section \ref{ssec:Res - mix} with superimposed discretization. The plots are categorized based on the shell theory used: Kirchhoff-Love/Kirchhoff-Love based contours are shown in (a), (b), and (c); Reissner-Mindlin/Reissner-Mindlin based contours are shown in (d), (e), and (f); Kirchhoff-Love/Reissner-Mindlin based contours are shown in (g), (h), and (i).} \label{fig:RES - MixCont}
\end{figure}

Fig.\figref{fig:RES - MixCont} shows contour plots of the displacement magnitude, the first component of the membrane force, and the first components of the bending moment. The figures illustrate the smooth coupling between the boundary layers and the internal patch for all represented quantities. There is minimal difference between the Kirchhoff-Love and Reissner-Mindlin solutions due to the shell's small thickness. However, artifacts in the bending moment are noticeable at the interface between the boundary layer of the cut-out and the internal patch. These artifacts are mitigated by employing Reissner-Mindlin theory for the boundary layer.

\begin{figure}	\centering
    \includegraphics[width = 0.70\textwidth]{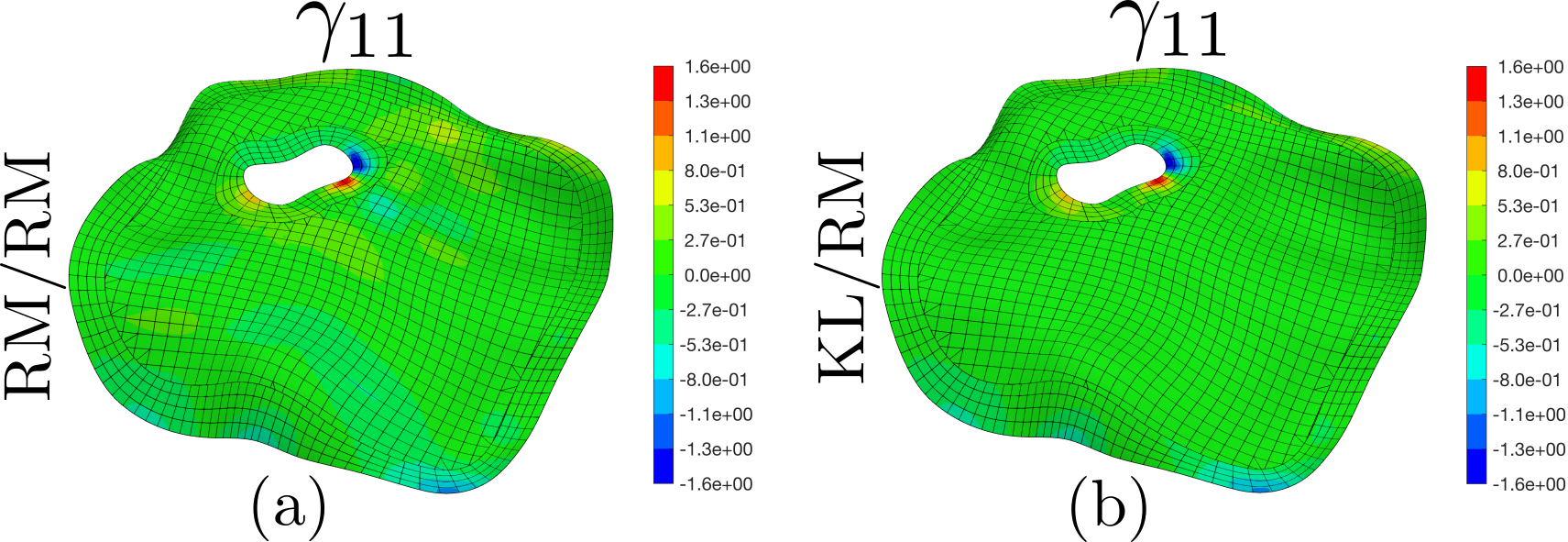}
    \caption{Contour plots of the first component of the shear strain tensor for a Reissner-Mindlin/Reissner-Mindlin discretization (a) and for a Kirchhoff-Love/Reissner-Mindlin discretization (b). For the mixed discretization, the Kirchhoff-Love theory cannot capture the shear strain in the internal patch. However, at the boundaries, where concentrations occur and delamination can initiate, the Reissner-Mindlin theory effectively models these potentially dangerous localizations.} \label{fig:RES - MixGamm}
\end{figure}

Additionally, Fig.\figref{fig:RES - MixGamm} shows the contour plots for the first component of the shear strain tensor for a Reissner-Mindlin/Reissner-Mindlin discretization and for a Kirchhoff-Love/Reissner-Mindlin one. The Kirchhoff-Love theory assumes the shear strain to be null. Therefore, to accurately capture this field at the boundaries, where 
local concentrations occur, a more refined shell theory, such as the Reissner-Mindlin one, is required. This refinement is particularly important for laminates, where the associated shear stress can initiate delamination phenomena.

The adoption of different formulations for boundary layers and the internal patch offers several advantages. Firstly, it allows for local enhancement of analysis accuracy where needed, depending on the application requirements. Additionally, Reissner-Mindlin boundary layers enable strong enforcement of rotation boundary conditions, and at the same time the overall number of degrees of freedom is kept limited thanks to the Kirchhoff-Love kinematics of the internal patch.

It is important to note that the boundary layer map results from the composition of two maps as detailed in Section \ref{ssec:IBCM boundary layers}, where the underlying map is constructed using a B-spline. Consequently, lines corresponding to its knots exhibit reduced continuity and pass through the internal parts of boundary layer elements, leading to certain challenges. Higher-order integration on maps with non-$C^\infty$ continuity cannot be achieved using standard Gaussian schemes. Potential solutions include reparametrizing elements into tiles where the map is locally $C^\infty$ or using different integration schemes suitable for reduced continuity maps, such as Newton–Cotes based methods.

Another issue arises when using elements with orders higher than the continuity of the underlying map. In such instances, these elements are expected to behave asymptotically as if their order conformed to the continuity of the map itself. Adopting Reissner-Mindlin kinematics for boundary layers helps mitigate this issue in the preasymptotic regime, as this is a smaller order problem when compared with Kirchhoff-Love one. However, a thorough investigation of this behavior is not within the scope of this study and is reserved for future research.

\subsection{Damaged cylindrical shell under internal pressure} \label{ssec:RES - dam}
\begin{figure}	\centering
    \includegraphics[width = 0.80\textwidth]{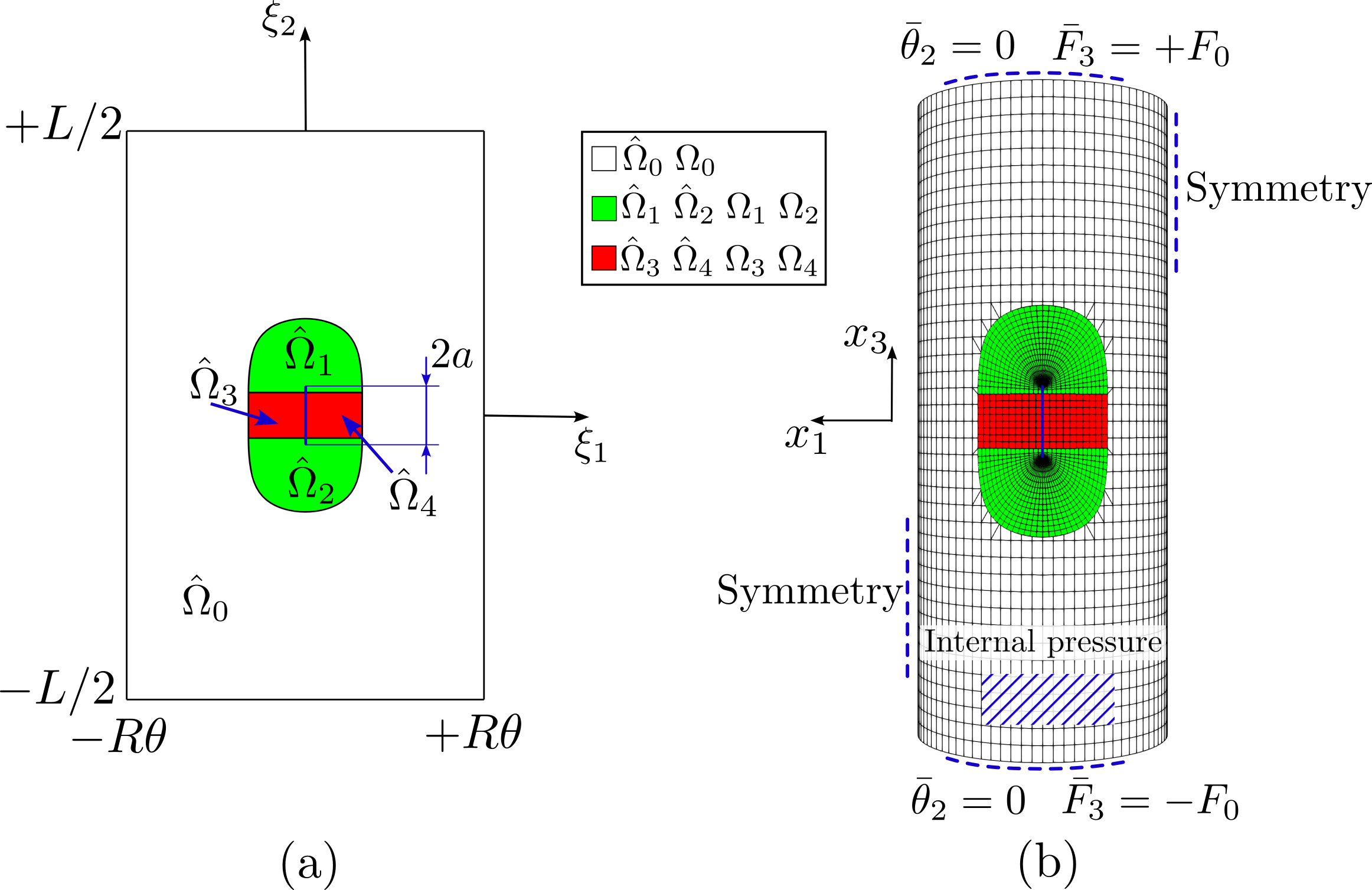}
    \begin{subfigure}{0\textwidth}\phantomcaption\label{fig:RES - DamGeom a}\end{subfigure}
    \begin{subfigure}{0\textwidth}\phantomcaption\label{fig:RES - DamGeom b}\end{subfigure}
    \caption{Geometry for the test described in Section \ref{ssec:RES - dam} in parametric domain (a) and Euclidean domain (b) with superimposed mesh.} \label{fig:RES - DamGeom}
\end{figure}

This test aims to demonstrate the capabilities of the IBCM in modeling damaged structures. Fig.\figref{fig:RES - DamGeom} illustrates the setup for this test: a cylinder with an axial crack subjected to an internal pressure $p_0$. The shell's analytic mapping is defined as:
\begin{equation}
    \bm{x} (\xi_1,\xi_2) = \begin{bmatrix} -R \sin{(\xi_1/R)}\\ +R \cos{(\xi_1)/R}\\ \xi_2 \end{bmatrix} \;,
\end{equation}
where $(\xi_1,\xi_2)\in[-R\theta,+R\theta]\times[-L/2,+L/2]$, with $L=100$ [mm], $\theta=\pi/2$ [rad], and $R=20$ [mm]. The axial crack, centered within the cylinder, has a length $2a=10$ [mm].

The material of the shell is isotropic having Young modulus $E$, Poisson's ratio $\nu=1/3$, and the shell section has a thickness $\tau=1$ [mm]. The shell is subjected to boundary conditions approximating those of an indefinitely long cylinder with closed top and bottom parts. Specifically: an external distributed force $\bar{\bm{F}} = \pm F_0 \bm{e}_3$ is applied on the boundaries at $\xi_2=\pm L/2$, with $F_0=p_0R/2$. The rotation component $\theta_2$ , referring to the local contravariant basis vector $\bm{a}^2$, is constrained on these boundaries. Symmetry boundary conditions $u_2=0$ and $\theta_1=0$ are applied on the boundaries $\xi_1=\pm R\theta$. An internal pressure $p_0$ results in a domain traction $\bar{\bm{t}}=p_0\bm{a}_3$.

The IBCM is employed to model the region surrounding the axial crack using four auxiliary patches: $\hat{\Omega}_1$, $\hat{\Omega}_2$, $\hat{\Omega}_3$, and $\hat{\Omega}_4$. $\hat{\Omega}_1$ and $\hat{\Omega}_2$ correspond to the top and bottom regions at the tips of the crack, while $\hat{\Omega}_3$ and $\hat{\Omega}_4$ represent the left and right side regions of the crack. These boundary patches are constructed in the parametric domain as depicted in Fig.\figref{fig:RES - DamGeom a}. Fig.\figref{fig:RES - DamGeom b} shows the same partitioning in the Euclidean domain with superimposed discretization, illustrating how using four patches for the boundary layer allows for tuning accuracy in regions where different stress concentrations are expected around the crack.

Furthermore, it is important to note that this configuration effectively resolves the issue of cross-talk. In this case, there is no non-active region that needs to be trimmed away. If only a single patch were used, there would be no discernible difference in the analysis between the scenarios with and without the crack. This is because disconnecting degrees of freedom for B-splines is not straightforward and becomes impossible if the crack intersects inner elements' domains. However, by employing a boundary layer approach, the necessary discontinuity is introduced into the approximation space, enabling accurate modeling of the crack opening.

\begin{figure}\centering \begin{tikzpicture} 
    \node[anchor=south west, inner sep=0] (image_ref) at (0,0) {\includegraphics[width = 0.5\textwidth]{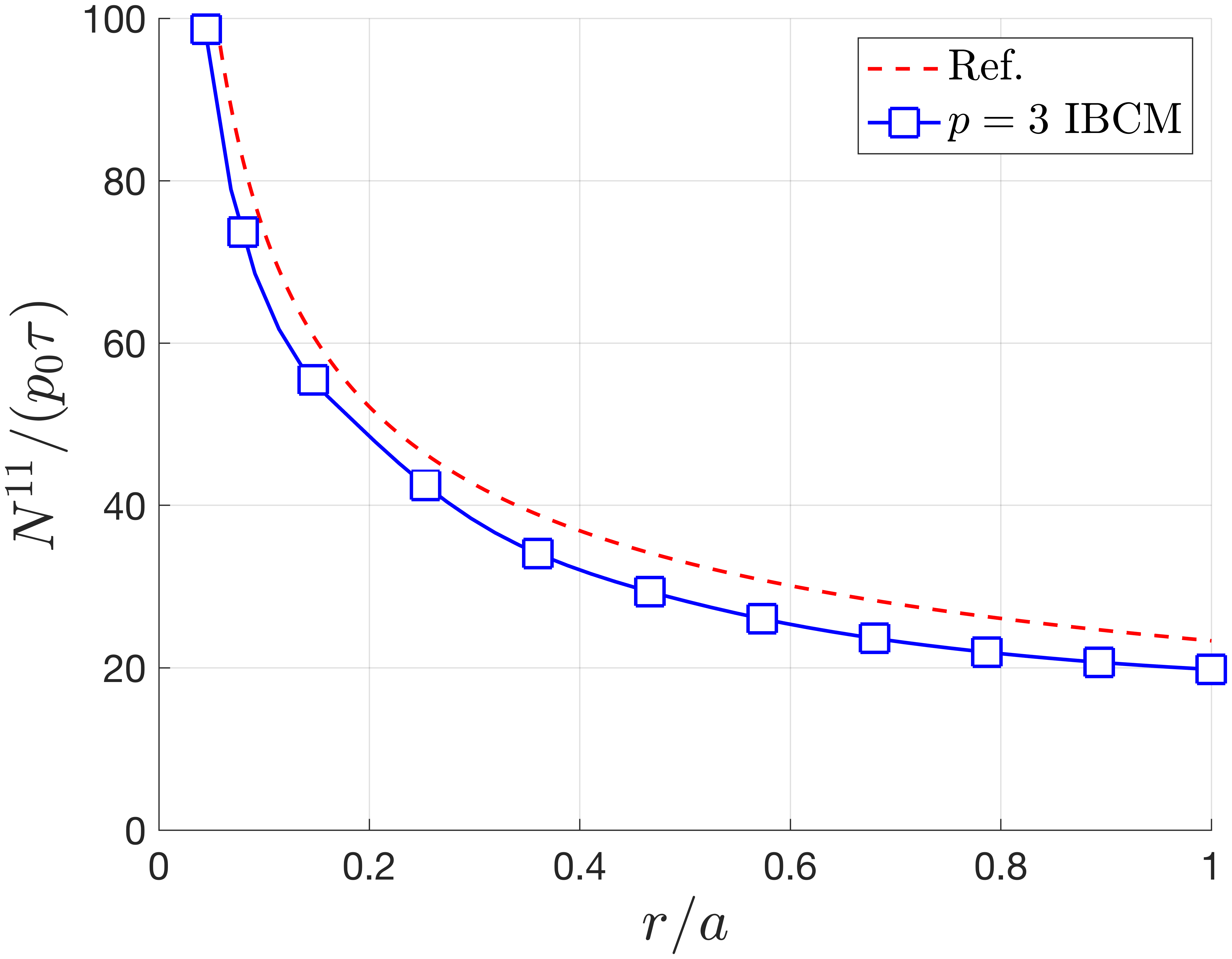}};
    \node (x1) at (7.25,6.05) [rectangle,rotate=0,anchor=center,color=black] {\small \cite{folias1969}};
    \end{tikzpicture}\caption{Comparison between the analytical and numerical solution to the problem in Section \ref{ssec:RES - dam}, regarding the non-dimensional membrane force $N^{11}/(p_0\tau)$ as a function of the distance with the crack tip on an axially-oriented segment.}\label{fig:RES - DamGrap}
\end{figure}

\begin{figure}	\centering
    \includegraphics[width = 0.80\textwidth]{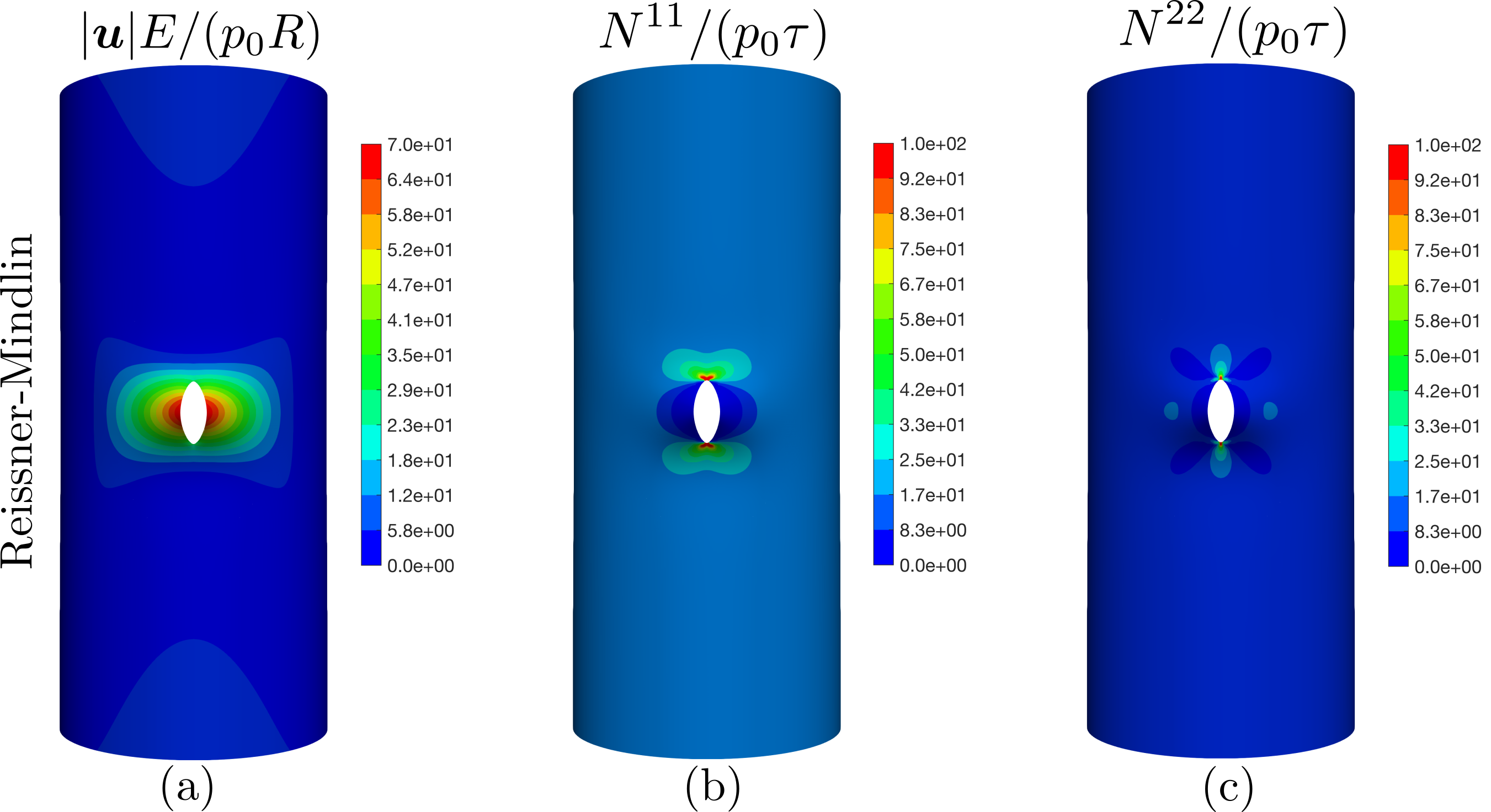}
    \caption{Contour of the non-dimensional magnitude of the displacement, and two non-dimensional components of the membrane force for the damaged cylindrical shell described in Section \ref{ssec:RES - dam}.} \label{fig:RES - DamCont}
\end{figure}

For a shallow and infinitely long cylinder, an analytical solution exists for stress at the crack tip, as demonstrated in \cite{folias1969}, and reported here for the sake of completeness for $\nu=1/3$
\begin{equation}
    N^{11}/(p_0\tau) = \sqrt{\frac{c}{2r}} \left(1+(0.37-0.30\ln{\lambda})\lambda^2 \vphantom{\frac{R}{\tau}}\right) \left(\frac{R}{\tau}\right) \;, \quad \mathrm{where}\; \lambda^4=\frac{12(1-\nu^2)a^4}{R^2\tau^2}\;.
\end{equation}
Fig.\figref{fig:RES - DamGrap} compares this analytical solution with the numerical results, representing the non-dimensional membrane force $N^{11}/(p_0\tau)$ as a function of the local distance from the crack tip along an axial segment. Despite some noticeable errors near the crack tip due to the singular solution, good agreement is observed. Finally, Fig.\figref{fig:RES - DamCont} presents the non-dimensional magnitude of displacement $|\bm{u}|E/(\sigma R)$, and the non-dimensional components of the membrane force $N^{11}/(p_0\tau)$ and $N^{22}/(p_0\tau)$.

\section{Conclusions}\label{sec:CONCLUSIONS}
This work has explored the application of the Immersed Boundary Conformal Method (IBCM) to analyze Kirchhoff-Love and Reissner-Mindlin shell structures within an immersed domain framework. By employing this innovative method, we successfully addressed several challenges inherent to traditional immersed boundary approaches. Specifically, the IBCM method offers the following benefits:

\begin{enumerate}
\item Due to the conformal nature of otherwise trimmed boundaries, Dirichlet conditions on primary variables can be enforced strongly by acting directly on the degrees of freedom associated with the spline corresponding to the specific edge. 
\item IBCM facilitates targeted mesh refinement, which is crucial for accurately capturing stress concentrations and other localized gradients in the solution.

\item Continuity between the boundary layer and the inner patch is enforced weakly through Nitsche-based coupling. However, the average operator is chosen non-symmetrically by computing the fluxes only from the conformal boundary layer, allowing for straightforward stable contributes with bounded penalty parameters. 
\item In structures comprised of multiple patches intersecting in a non-conformal manner, IBCM can be adopted to create boundary layers from each patch, restoring a conformal interface and allowing for seamless strong coupling by directly matching degrees of freedom pairs.
\item Boundary layers can be constructed to address the cross-talk phenomenon that arises in traditional spline-based methods. These auxiliary patches introduce the necessary discontinuities in the approximation space.
\item Local physical phenomena can be modeled accurately through refined equations in the boundary layer, provided adequate interface conditions to link it with the inner patch.
\end{enumerate}

Through comprehensive numerical experiments, the proposed formulation has demonstrated its capability to accurately and efficiently handle shell structures. Optimal convergence rates were achieved in tests where analytical reference solutions are available. The coercivity of the bilinear form was tested numerically through the positive definiteness of the stiffness matrix and was verified in all the tests investigated, unlike in the standard trimmed approach adopted for comparison. A mixed Kirchhoff-Love/Reissner-Mindlin approach was employed to locally refine the shell theory at the boundary of an embedded domain. Finally, the IBCM introduced the necessary discontinuity to model an axial crack in a cylindrical shell, enabling at the same time straightforward local refinement at its apex.

In summary, the integration of the IBCM method, B-splines basis functions, and Nitsche-based coupling techniques shows great potential in the field of shell analysis. This approach improves the accuracy and stability of numerical solutions while retaining the versatility of an embedded geometrical description of the mid-surface.

\section*{Acknowledgements}
This research was supported by the European Union Horizon 2020 research and innovation program, under grant agreement No 862025 (ADAM2), and by the Swiss National Science Foundation through the project No 200021\_214987 (FLASh). 

\section*{Declaration of competing interest}
The authors declare that they have no known competing financial interests or personal relationships that could have appeared to influence the work reported in this paper.

\appendix
\section{Composition of maps} \label{sec:APP - MAP COMP}
The map of the boundary layers in the Euclidean space is given as a composition of maps
\begin{equation}
    \bm{x}\left(\eta_1^i,\eta_2^i\right) = \hat{\bm{\mathcal{F}}}\circ\tilde{\bm{\mathcal{F}}}\left(\eta_1^i,\eta_2^i\right)\;.
\end{equation}
According to the problem of interest built over the surface, some derivatives in the domain of the boundary layers might be necessary. For completeness, these derivatives are reported here up to the third order
\begin{subequations}
\begin{align}
    &\fpd{  \bm{x}}{\eta_\alpha}= \fpd{\bm{x}}{\xi_\lambda}\fpd{\xi_\lambda}{\eta_\alpha}\;,\\
    &\fpd{^2\bm{x}}{\eta_\alpha\pd\eta_\beta}= \fpd{^2\bm{x}}{\xi_\lambda\pd\xi_\mu}\fpd{\xi_\lambda}{\eta_\alpha}\fpd{\xi_\mu}{\eta_\beta}+ \fpd{\bm{x}}{\xi_\lambda}\fpd{^2\xi_\lambda}{\eta_\alpha\pd\eta_\beta}\;,\\
    &\fpd{^3\bm{x}}{\eta_\alpha\pd\eta_\beta\pd\eta_\gamma}= 
    \fpd{^3\bm{x}}{\xi_\lambda\pd\xi_\mu\pd\xi_\nu}\fpd{\xi_\lambda}{\eta_\alpha}\fpd{\xi_\mu}{\eta_\beta}\fpd{\xi_\nu}{\eta_\gamma}+
    \fpd{^2\bm{x}}{\xi_\lambda\pd\xi_\mu}\fpd{^2\xi_\lambda}{\eta_\alpha\pd\eta_\gamma}\fpd{\xi_\mu}{\eta_\beta}+
    \fpd{^2\bm{x}}{\xi_\lambda\pd\xi_\mu}\fpd{\xi_\lambda}{\eta_\alpha}\fpd{^2\xi_\mu}{\eta_\beta\eta_\gamma}+\nonumber\\
    &\hspace{2cm}+\fpd{^2\bm{x}}{\xi_\lambda}\fpd{^2\xi_\lambda}{\eta_\alpha\pd\eta_\beta}\fpd{\xi_\mu}{\eta_\gamma}+
     \fpd{\bm{x}}{\xi_\lambda}\fpd{^3\xi_\lambda}{\eta_\alpha\pd\eta_\beta\pd\eta_\mu}  \;,
\end{align}
\end{subequations}
where the superscript ``$i$'' in the auxiliary curvilinear coordinates $\eta_1^i$ and $\eta_2^i$ has been discarded to enhance readability. 

\section{Laminate layout and constitutive behaviour}\label{sec:APP constitutive}
In order to derive the constitutive relationship between strain and stress, let us introduce a local orthonormal reference system $\bm{n}_1 \bm{n}_1 \bm{n}_3$, whose vectors are defined as
\begin{equation}
    \begin{aligned}
         &\bm{n}_1 = {\bm{a}_1} / {|\bm{a}_1|}  \;, \\
         &\bm{n}_2 = {\bm{a}^2} / {|\bm{a}^2|}  \;, \\
         &\bm{n}_3 = {\bm{a}_3}                 \;.
    \end{aligned}
\end{equation}
It is assumed, here, that each layer has one of the three orthotropic directions oriented along $\bm{n}_3$, whereas the main direction, that coincides with the direction of deposition for fiber-reinforced composites, forms an angle $\theta^{\lam{\ell}}$ from $\bm{n}_1$ for a rotation along $\bm{n}_3$. In the orthotropic direction, the relationship between stress and strain is expressed as
\begin{align}
    &\tilde{\bm{\sigma}}_L = \tilde{\bm{c}}^\lam{\ell}_L \tilde{\bm{\epsilon}}_L \;, \\
    &\tilde{\bm{\sigma}}_T = \tilde{\bm{c}}^\lam{\ell}_T \tilde{\bm{\epsilon}}_T \;,
\end{align}
where $\tilde{\bm{\sigma}}_L=\{\tilde{\sigma}^{11},\tilde{\sigma}^{22},\tilde{\sigma}^{12}\}$ and $\tilde{\bm{\epsilon}}_L = \{\tilde{\epsilon}_{11},\tilde{\epsilon}_{22},2\tilde{\epsilon}_{12}\}$ collect the in-plane stress and strain in Voigt notation, respectively, while $\tilde{\bm{\sigma}}_T=\{ \tilde{\sigma}^{31},\tilde{\sigma}^{32}\}$ and $\tilde{\bm{\epsilon}}_T = \{ 2\tilde{\epsilon}_{31}, 2\tilde{\epsilon}_{32}\}$ collect the out-of-plane stress and strain in Voigt notation, respectively. It is reminded that, in classical shell theories, the plane stress hypothesis states that $\tilde{\sigma}_{33}=0$, uniformly over the shell surface. The matrices $\tilde{\bm{c}}^\lam{\ell}_L$ and  $\tilde{\bm{c}}^\lam{\ell}_T$ are the in-plane and out-of-plane stiffness matrices built through the engineering elastic coefficients, namely the Young's moduli $E_1^\lam{\ell}$, $E_2^\lam{\ell}$, the shear moduli $G_{12}^\lam{\ell}$, $G_{31}^\lam{\ell}$ and $G_{32}^\lam{\ell}$ and the Poisson's ratio $\nu_{12}^\lam{\ell}$ as detailed in \cite{jones1998}. Provided these quantities, the stiffness matrices are constructed as
\begin{align}
    &\tilde{\bm{c}}^\lam{\ell}_L 
    = 
    \begin{bmatrix}
                1/E_1^\lam{\ell} & -\nu_{12}^\lam{\ell}/E_1^\lam{\ell} &        0\\
        -\nu_{21}^\lam{\ell}/E_2^\lam{\ell} &         1/E_2^\lam{\ell} &        0\\
                    0 &             0 & 1/G_{12}^\lam{\ell}
    \end{bmatrix}^{-1} \;,\\
    & \tilde{\bm{c}}^\lam{\ell}_T
    =
    \begin{bmatrix}
     1/(\alpha_s G_{31}^\lam{\ell})& 0\\
     0& 1/(\alpha_s G_{32}^\lam{\ell})
    \end{bmatrix}^{-1}\;,
\end{align}
where $\alpha_s$ is the so-called shear correction factor, that is assumed equal to $5/6$, unless differently stated. The constitutive relationship can be expressed in the common orthonormal basis $\bm{n}_1\bm{n}_2\bm{n}_3$ as
\begin{align}
    &\bar{\bm{\sigma}}_L = \bar{\bm{c}}^\lam{\ell}_L \bar{\bm{\epsilon}}_L \;, \\
    &\bar{\bm{\sigma}}_T = \bar{\bm{c}}^\lam{\ell}_T \bar{\bm{\epsilon}}_T \;,
\end{align}
where $\bar{\bm{c}}^\lam{\ell}_L = \bm{T}_L\tilde{\bm{c}}_L^\lam{\ell}\bm{T}^{\Tr}_L$ and $\bar{\bm{c}}^\lam{\ell}_T = \bm{T}_T\tilde{\bm{c}}_T^\lam{\ell}\bm{T}^{\Tr}_T$ are the stiffness matrices in the basis $\bm{n}_1\bm{n}_2\bm{n}_3$, and the rotation matrices $\bm{T}_L$ and $\bm{T}_T$ are defined as
\begin{subequations}
\begin{align}
    &\bm{T}_L = 
    \begin{bmatrix}
            \cos^2\theta^\lam{\ell} &          \sin^2\theta^\lam{\ell} &    -2\sin\theta^\lam{\ell}\cos\theta\\
            \sin^2\theta^\lam{\ell} &          \cos^2\theta^\lam{\ell} &     2\sin\theta^\lam{\ell}\cos\theta^\lam{\ell}\\
    \sin\theta^\lam{\ell}\cos\theta^\lam{\ell} & -\sin\theta^\lam{\ell}\cos\theta^\lam{\ell} & \cos^2\theta^\lam{\ell}-\sin^2\theta^\lam{\ell}  
        \end{bmatrix} \;,\\
    &\bm{T}_T = 
    \begin{bmatrix}
         \cos\theta^\lam{\ell} &     -\sin\theta^\lam{\ell}\\
         \sin\theta^\lam{\ell} &      \cos\theta^\lam{\ell}
    \end{bmatrix} \;.
\end{align}
\end{subequations}
The virtual work of the internal forces in 3D is simplified integrating along the thickness direction $\xi_3$ and substituting the constitutive relationship as
\begin{equation}
    \delta\mathcal{L}_{int}=\sum_{\ell=1}^{N_\ell}{\int_{V^\lam{\ell}}{\delta\epsilon_{ij}\sigma^{ij}\dd V}} =  \int_\Omega{\left(\sum_{\ell=1}^{N_\ell}
    {\int_{\xi_{3b}^\lam{\ell}}^{\xi_{3t}^\lam{\ell}}{(\delta\bar{\bm{\epsilon}}_L^\Tr\bar{\bm{\sigma}}_L + \delta\bar{\bm{\epsilon}}_T^\Tr\bar{\bm{\sigma}}_T)}\dd\xi_3}\right)\dd\Omega} \;,
\end{equation}
where the volume increment was approximated as $\dd V =\dd\xi_3\dd\Omega$. Upon substituting the expression for the strain and the constitutive relationship in the local orthonormal basis, the previous equation leads to the definition of the generalized stiffness matrices
\begin{subequations}
\begin{align}
    &\bar{\bm{A}} = \sum_{\ell=1}^{N_\ell}\int_{\tau_b^\lam{\ell}}^{\tau_t^\lam{\ell}}{\overline{\bm{c}}_L^\lam{\ell}       \dd\xi_3}  \;,\\
    &\bar{\bm{B}} = \sum_{\ell=1}^{N_\ell}\int_{\tau_b^\lam{\ell}}^{\tau_t^\lam{\ell}}{\overline{\bm{c}}_L^\lam{\ell}\xi_3  \dd\xi_3}  \;,\\
    &\bar{\bm{D}} = \sum_{\ell=1}^{N_\ell}\int_{\tau_b^\lam{\ell}}^{\tau_t^\lam{\ell}}{\overline{\bm{c}}_L^\lam{\ell}\xi_3^2\dd\xi_3}  \;,\\
    &\bar{\bm{S}} = \sum_{\ell=1}^{N_\ell}\int_{\tau_b^\lam{\ell}}^{\tau_t^\lam{\ell}}{\overline{\bm{c}}_T^\lam{\ell}\dd\xi_3}  \;.
\end{align}
\end{subequations}
To obtain the constitutive tensor, it is sufficient to pass from Voigt notation to tensor notation. Taking $\bar{\bm{A}}$ as an example, the associated stiffness tensor notation is established as $\bar{\mathbb{A}}^{\alpha\beta\gamma\delta}=\bar{A}^{ab}$, using the correspondences $\alpha\beta\longleftrightarrow a$ and $\gamma\delta\longleftrightarrow b$, where the indices 11, 22, 12, and 21 correspond to 1, 2, 3, and 3, respectively.  Instead, the constitutive tensor $\bar{\mathbb{S}}$ is obtained directly expressing as a tensor the matrix $\bar{\bm{S}}$, i.e., $\bar{\mathbb{S}}^{\alpha\beta}=\bar{S}^{\alpha\beta}$. It is worth remarking that these constitutive tensors in the local orthonormal basis $\bm{n}_1\bm{n}_2\bm{n}_3$ depends only on the lamination sequence and not on the geometry of the mid-surface. As such, for a uniform laminate material, the previous integration needs to be performed only once. Finally, to pass from the local orthonormal basis to the local covariant basis, the following transformation laws are employed
\begin{align}
    &\mathbb{A}^{\alpha_1\beta_1\gamma_1\delta_1} = \bar{\mathbb{A}}^{\alpha_1\beta_1\gamma_1\delta_1}\;(\bm{n}_{\alpha_2}\cdot\bm{a}^{\alpha_1})(\bm{n}_{\beta_2}\cdot\bm{a}^{\beta_1})(\bm{n}_{\gamma_2}\cdot\bm{a}^{\gamma_1})(\bm{n}_{\delta_2}\cdot\bm{a}^{\delta_1})  \;, \\    &\mathbb{S}^{\alpha_1\beta_1} = \bar{\mathbb{S}}^{\alpha_1\beta_1}\;(\bm{n}_{\alpha_2}\cdot\bm{a}^{\alpha_1})(\bm{n}_{\beta_2}\cdot\bm{a}^{\beta_1})    \;.
\end{align}


\bibliography{Bibliography}

\end{document}